\documentclass[12pt]{article}
\usepackage{latexsym,amsmath,a4wide,theorem,amsfonts,eucal,amssymb}

\def \beq {\begin{eqnarray}}
\def \eeq {\end{eqnarray}}
\def \beqn {\begin{eqnarray*}}
\def \eeqn {\end{eqnarray*}}

\newcommand{\halmos}{\rule{1ex}{1.4ex}}

\newcounter{for}[section]

\newtheorem{itlemma}{Lemma}[section]
\newtheorem{itproposition}[itlemma]{Proposition}
\newtheorem{itfact}[itlemma]{Fact}
\newtheorem{theorem}[itlemma]{Theorem}
\newtheorem{itcorollary}[itlemma]{Corollary}
\newtheorem{itremark}[itlemma]{Remark}
\newtheorem{itremarks}[itlemma]{Remarks}
\newtheorem{itdefinition}[itlemma]{Definition}
\newtheorem{itexample}[itlemma]{Example}

\newenvironment{fact}{\begin{itfact}\rm}{\end{itfact}}
\newenvironment{claim}{\begin{itclaim}\rm}{\end{itclaim}}
\newenvironment{lemma}{\begin{itlemma}}{\end{itlemma}}
\newenvironment{remark}{\begin{itremark}\rm}{\end{itremark}}
\newenvironment{remarks}{\begin{itremarks} \rm}{\end{itremarks}}
\newenvironment{corollary}{\begin{itcorollary}}{\end{itcorollary}}
\newenvironment{proposition}{\begin{itproposition}}{\end{itproposition}}
\newenvironment{definition}{\begin{itdefinition}\rm}{\end{itdefinition}}
\newenvironment{example}{\begin{itexample}\rm}{\end{itexample}}
\newenvironment{proof}{\noindent {\em Proof}.\ \
}{\hspace*{\fill}$\halmos$\medskip}
\newcommand{\be}[1]{\addtocounter{for}{1} \begin{equation}\label{#1}}
\newcommand{\ee}{\end{equation}}
\newcommand{\bl}[1]{\begin{lemma}\label{#1}}
\newcommand{\br}[1]{\begin{remark}\label{#1}}
\newcommand{\brs}[1]{\begin{remarks}\label{#1}}
\newcommand{\bt}[1]{\begin{theorem}\label{#1}}
\newcommand{\bd}[1]{\begin{definition}\label{#1}}
\newcommand{\bp}[1]{\begin{proposition}\label{#1}}
\newcommand{\bfact}[1]{\begin{fact}\label{#1}}
\newcommand{\bc}[1]{\begin{corollary}\label{#1}}
\newcommand{\bex}[1]{\begin{example}\label{#1}}
\newcommand{\ec}{\end{corollary}}
\newcommand{\efact}{\end{fact}}
\newcommand{\eex}{\end{example}}
\newcommand{\el}{\end{lemma}}
\newcommand{\er}{\end{remark}}
\newcommand{\ers}{\end{remarks}}
\newcommand{\et}{\end{theorem}}
\newcommand{\ed}{\end{definition}}
\newcommand{\ep}{\end{proposition}}
\newcommand{\epr}{\end{proof}}
\newcommand{\bpr}{\begin{proof}}
\newcommand{\bcl}[1]{\begin{claim}\label{#1}}
\newcommand{\ecl}{\end{claim}}

\newcommand{\ecs}{\end{corollary}}
\newcommand{\eers}{\end{exercise}}
\newcommand{\eexs}{\end{example}}
\newcommand{\eems}{\end{example}}
\newcommand{\els}{\end{lemma}}
\newcommand{\eles}{\end{lemmaex}}
\newcommand{\ets}{\end{theorem}}
\newcommand{\eds}{\end{definition}}
\newcommand{\eps}{\end{proposition}}

\newcommand{\bi}{\begin{itemize}}
\newcommand{\ei}{\end{itemize}}
\newcommand{\ben}{\begin{enumerate}}
\newcommand{\een}{\end{enumerate}}

\def\vbar{\mathchoice{\vrule height6.3ptdepth-.5ptwidth.8pt\kern-.8pt}
   {\vrule height6.3ptdepth-.5ptwidth.8pt\kern-.8pt}
   {\vrule height4.1ptdepth-.35ptwidth.6pt\kern-.6pt}
   {\vrule height3.1ptdepth-.25ptwidth.5pt\kern-.5pt}}
\def\fudge{\mathchoice{}{}{\mkern.5mu}{\mkern.8mu}}
\def\bbc#1#2{{\rm \mkern#2mu\vbar\mkern-#2mu#1}}
\def\bbb#1{{\rm I\mkern-3.5mu #1}}
\def\bba#1#2{{\rm #1\mkern-#2mu\fudge #1}}
\def\bb#1{{\count4=`#1 \advance\count4by-64 \ifcase\count4\or\bba A{11.5}\or
   \bbb B\or\bbc C{5}\or\bbb D\or\bbb E\or\bbb F \or\bbc G{5}\or\bbb H\or
   \bbb I\or\bbc J{3}\or\bbb K\or\bbb L \or\bbb M\or\bbb N\or\bbc O{5} \or
   \bbb P\or\bbc Q{5}\or\bbb R\or\bbc S{4.2}\or\bba T{10.5}\or\bbc U{5}\or
   \bba V{12}\or\bba W{16.5}\or\bba X{11}\or\bba Y{11.7}\or\bba Z{7.5}\fi}}

\def \Z {{\mathbb Z}}
\def \R {{\mathbb R}}

\def \N {{\mathbb N}}

\def \ra {\rightarrow }

\def \O{\Omega}
\def \s {\sigma}

\def \LL {{\cal{L}}}

\def \A {{\cal{A}}}

\def \dir {{\cal{D}}}

\newcommand{\ba}[1]{\addtocounter{for}{1} \begin{eqnarray}\label{#1}}
\newcommand{\ea}{\end{eqnarray}}

\def\sqr#1#2{{\vcenter{\vbox{\hrule height .#2pt
                             \hbox{\vrule width .#2pt height#1pt \kern#1pt
                                   \vrule width .#2pt}
                             \hrule height .#2pt}}}}

\def\pmb#1{\setbox0=\hbox{#1}%
   \kern-.025em\copy0\kern-\wd0
   \kern.05em\copy0\kern-\wd0
   \kern-.025em\raise.0433em\box0 }
\def\sqr#1#2{{\vcenter{\vbox{\hrule height.#2pt
     \hbox{\vrule width.#2pt height#1pt \kern#1pt
   \vrule width.#2pt}\hrule height.#2pt}}}}

\def\A{{\mathcal A}}

\def\EE{{\mathcal E}}   %

\def\e{\epsilon}                

\def\e{\epsilon}
\def\d{\delta}
\def\l{\lambda}
\def\L{\Lambda}
\def\n{\nu_{\Lambda}}

\def\mur{\mu_{\rho}}
\def\g{\gamma}

\def\a{\alpha}

\def\nuln{\nu_{\L}^N}
\def\mur{\mu_{\rho}}

\newcommand{\ov}[1]{\overline{#1}}

\numberwithin{equation}{section}

\newtheorem{teo}{Theorem}[section]
\newtheorem{lem}[teo]{Lemma}
\newtheorem{pro}[teo]{Proposition}
\newtheorem{cor}[teo]{Corollary}

\newtheorem{co}[teo]{Condition}

\DeclareMathOperator{\ent}{Ent}

\newcommand{\ind}{\mathbf{1}}
\newcommand{\real}{\mathbb{R}}
\newcommand{\integer}{\mathbb{Z}}
\renewcommand{\natural}{\mathbb{N}}
\renewcommand{\tilde}{\widetilde}

\newcommand{\ie}{\emph{i.e.\ }}

\def\n#1#2{\overline{#1}_{#2}}

\def\proofof#1{\noindent{\bf Proof of #1.}\ }
\def\endproofof{{\mbox{}\nolinebreak\hfill\rule{2mm}{2mm}\par\medbreak} }

\begin{document}

\title{Logarithmic Sobolev Inequality \\
  for  \\
  Zero--Range Dynamics: independence of the number of particles}
\author{Paolo Dai Pra \\
Dipartimento di Matematica Pura e Applicata \\
Universit\`{a} di Padova \\ Via Belzoni 7, 35131 Padova, Italy \\ 
e-mail: \texttt{daipra@math.unipd.it}\\~
\\Gustavo Posta \\
Dipartimento di Matematica \\ Politecnico di Milano \\ Piazza L. Da Vinci 32, 20133 Milano, Italy\\
e-mail: \texttt{guspos@mate.polimi.it}
}
\date{}

\maketitle

\begin{abstract}
We prove that the  logarithmic-Sobolev constant for Zero-Range Processes in a box of diameter $L$ may depend on $L$ but not on the number of particles. This is a first, but relevant and quite technical step, in the proof that this logarithmic-Sobolev constant grows as $L^2$, that will
be presented in a forthcoming paper (\cite{Da:Po}).
\end{abstract}



{\em Abbreviated title}: log-Sobolev inequality for Zero-Range.

\newpage

\section{Introduction}

The zero-range process is a system of interacting particles moving in a discrete lattice
$\L$, that here we will assume to be a subset of $\Z^d$. The interaction is ``zero range'', i.e.
the motion of a particle may be only affected by particles located at the same lattice site. Let
$c: \N \ra [0,+\infty)$ be a function such that $c(0) = 0$
and $c(n) >0$ for every $n>0$. In the zero-range process associated to $c(\cdot)$ particles
evolve according to the following rule:
for each site $x \in \L$, containing $\eta_x$ particles, with probability rate $c(\eta_x)$ one
particle jumps from $x$ to one of its nearest neighbors at random. Waiting jump times of different sites
are independent. If $c(n) = \l n$ then particles are independent, and evolve as simple
random walks; nonlinearity of $c(\cdot)$ is responsible for the interaction. 
Note that particles are neither created nor destroyed. When $\L$ is a finite lattice, for each $N \geq 1$,
the zero-range process in $\L$ restricted to configurations with $N$ particles
is a finite irreducible Markov chain, whose unique invariant measure $\nuln$ is proportional to 
\be{int1}
\prod_{x \in \L} \frac{1}{c(\eta_x)!},
\ee
where
\[
c(n)! = \left\{
\begin{array}{ll}
1 & \mbox{for $n=0$} \\
c(n)c(n-1) \cdots c(1) & \mbox{otherwise.}
\end{array} \right.
\]
Moreover the process is reversible with respect to $\nuln$. 

If the function $n \mapsto c(n)$ does not grow too fast in $n$, then the zero range process
can be defined in the whole lattice $\Z^d$. In this case the extremal invariant measures form a one parameter
family of uniform product measures, with marginals
\be{int2}
\mur[\eta_x = k] = \frac{1}{Z(\a(\rho))}\frac{\a(\rho)^k}{c(k)!},
\ee
where $\rho \geq 0$, $Z(\a(\rho))$ is the normalization, and $\a(\rho)$ is uniquely determined
by the condition $\mur[\eta_x] = \rho$ (we use here the notation $\mu[f]$ for $\int f d\mu$).

In this paper we are interested in the rate of relaxation to equilibrium of zero-range processes. In 
general, for conservative systems of symmetric interacting  random walks in a spatial 
region $\L$, for which the interaction between particles 
is not too strong, the relaxation time to equilibrium is expected to be of the order of the
square of the diameter of $\L$, as it is the case for independent random walks. On a rigorous 
basis, however, this result has been proved in rather few cases. For dynamics with exclusion rule and finite-range interaction
(Kawasaki dynamics) relaxation to equilibrium with rate diam$(\L)^2$ has been proved 
(see \cite{Lu:Ya, Ca:Ma:Ro}) in the high temperature regime. For models without the exclusion
rule, i.e. with a possibly unbounded particle density, available results are even weaker (see \cite{La:Se:Va} and \cite{Po}). Zero-range processes are special models without the exclusion rule.
In some respect zero-range processes may appear simpler than Kawasaki dynamics: the 
interaction is zero-range rather than finite-range and, as a consequence, invariant measures 
have a simpler form. However, they exhibit various fundamental difficulties, including: 
\bi
\item
due to unboundedness of the density of particles various arguments used for exclusion processes
fail; in principle the rate of relaxation to equilibrium may depend on the number of particles, as it
actually does in some cases;
\item
there is no small parameter in the model; one would not like to restrict to ``small'' perturbations
of a system of independent particles.
\ei
In \cite{La:Se:Va} a first result for  zero-range processes has been obtained. Let
$\EE_{\nuln}$ be the Dirichlet form associated to the zero-range process in $\L = 
[0,L]^d \cap \Z^d$ with
$N$ particles. Then, under suitable growth conditions on $c(\cdot)$ (see Section 2),
the following Poincar\'e inequality holds 
\be{gap}
\nuln[f,f] \leq CL^2 \EE_{\nuln}(f,f),
\ee
where $C$ may depend on the dimension $d$ but not on $N$ or $L$. Moreover, by suitable
test functions, the $L-$dependence in (\ref{gap}) cannot be improved, i.e. one can find a positive
constant $c>0$ and functions $f = f_{N,L}$ so that $\nuln[f,f] \geq cL^2 \EE_{\nuln}(f,f)$ for all
$L,N$.
In other terms, (\ref{gap})
says that the spectral gap of $\EE_{\nuln}$ shrinks proportionally to $\frac{1}{L^2}$, independently
of the number of particles $N$. It is well known that Poincar\'e inequality controls convergence to
equilibrium in the $L^2(\nuln)$-sense: if $(S_t^{\L})_{t \geq 0}$ is the Markov semigroup 
corresponding to the process, then for every function $f$
\be{l2conv}
\nuln\left[ \left(S_t^{\L} f - \nuln[f]\right)^2\right] \leq \exp\left(-\frac{2t}{CL^2}\right) \nuln\left[
(f-\nuln[f])^2\right].
\ee
Poincar\'e inequality is however not sufficient to control convergence in stronger norms, e.g.
in total variation, that would follow from the logarithmic-Sobolev inequality
\be{ls}
\ent_{\nuln}(f) \leq s(L,N)  \EE_{\nuln}(\sqrt{f},\sqrt{f}),
\ee
where $\ent_{\nuln}(f) = \nuln[f\log f] - \nuln[f] \log \nuln[f]$. The constant  $s(L,N)$ in (\ref{ls})
is intended to be the smallest possible, and, in principle, may depend on both $L$ and $N$. 

Our
 main aim is to prove that $s(L,N) \leq CL^2$ for some $C>0$, i.e. the logarithmic-Sobolev
constant scales as the inverse of the spectral gap. This is the first conservative system with
unbounded particle density for which this scaling is established (see comments on page 423 of
\cite{La:Ki}). It turns out that the proof of this result is very long and technical, and it roughly
consists in two parts. In the first part one needs  to show that 
\be{rough}
s(L) := \sup_{N \geq 1} s(L,N) < +\infty,
\ee
i.e. that $s(L,N)$ has an upper bound independent of $N$, while in the second part a sharp
induction in $L$ is set up to prove the actual $L^2$ dependence. For this induction to work 
one has to choose $L$ sufficiently large as a starting point, and for this $L$ an upper bound for 
$s(L,N)$
independent of $N$ has to be known in advance. Note that for models with bounded particle
density inequality (\ref{rough}) is trivial.

This paper is devoted  to the proof of (\ref{rough}), while the induction leading to the $L^2$
growth is included in \cite{Da:Po}. The proof of (\ref{rough}) is indeed very long, and relies 
on quite sharp estimate on the measure $\nuln$. After introducing the model and stating the main result in Section 2, we devote Section 3 to the presentation of the essential steps of the proof, 
leaving the (many) technical details for the remaining sections.

\section{Notations and Main result}

Throughout this paper, for a given probability space $(\O, {\cal{F}}, \mu)$ and $f: \O \ra \R$
measurable, we use the following notations for mean value and covariance:
\[
\mu[f] := \int f d\mu,  \qquad  \mu[f,g] := \mu\left[ (f-\mu[f])(g-\mu[g]) \right]
\]
and, for $f \geq 0$,
\[
\ent_{\mu}(f) := \mu[f \log f] - \mu[f] \log \mu[f],
\]
where, by convention, $0\log 0 = 0$. Similarly, for ${\cal{G}}$ a sub-$\s$-field of ${\cal{F}}$, we let
$\mu[f|{\cal{G}}]$ to denote the conditional expectation, 
\[
\mu[f,g|{\cal{G}}] := \mu[(f-\mu[f|{\cal{G}}])(g-\mu[g|{\cal{G}}])|{\cal{G}}]
\]
the conditional covariance, and
\[
\ent_{\mu}(f|{\cal{G}}) := \mu[f \log f|{\cal{G}}] - \mu[f|{\cal{G}}] \log \mu[f|{\cal{G}}]
\]
the conditional entropy.\\
If $A \subset \O$, we denote by $\ind(x \in A)$ the indicator function of $A$.
If $B\subset A$ is \emph{finite} we will write $B\subset\subset A$.
For any $x\in\real$ we will write $\lfloor x\rfloor:=\sup\{n\in\integer: n\leq x\}$ and $\lceil x\rceil:=\inf\{n\in\integer: n\geq x\}$.

Let $\L$ be a possibly infinite subset of $\Z^d$, and $\O_{\L} = \N^{\L}$ be the corresponding
configuration space, where $\N = \{0,1,2,\ldots\}$ is the set of natural numbers. Given
a configuration $\eta \in \O_{\L}$ and $x \in \L$, the natural number $\eta_x$ will
be referred to as the number of particles at $x$.
Moreover if $\Lambda^\prime\subset\Lambda$ $\eta_{\Lambda^\prime}$ will denote the restriction of $\eta$ to $\Lambda^\prime$.
For two 
elements $\sigma,\xi\in\Omega_{\Lambda}$, the operations $\sigma\pm\xi$ are defined componentwise (for the difference whenever it returns an element of $\O_{\L}$). 
In what follows, given
$x \in \L$, we make use of the special configuration $\d^x$, having one particle at $x$
and no other particle. For $f :\O_{\L} \ra \R$ and $x,y \in \L$, we let
\[
\partial_{xy}f(\eta) := f(\eta - \d^x + \d^y) - f(\eta).
\]
Consider, at a formal level, the operator
\be{not1}
\LL_{\L} f(\eta) := \sum_{x \in \L}\sum_{y \sim x} c(\eta_x) \partial_{xy}f(\eta),
\ee
where $y \sim x$ means $|x-y| = 1$, and $c:\N \ra \R_+$ is a function such that $c(0) = 0$ and
$\inf\{c(n):n>0\} >0$.
In the case of $\L$ finite, for each $N \in \N \setminus \{0\}$, $\LL_{\L}$ is the infinitesimal
generator of a irreducible Markov chain on the finite state space $\{\eta \in \O_{\L}: \ov{\eta}_{\L}
=N\}$, where
\[
\ov{\eta}_{\L} := \sum_{x \in \L} \eta_x
\]
is the total number of particles in $\L$. The unique stationary measure for this Markov 
chain is denoted by $\nuln$ and is given by
\be{not2}
\nuln[\{\eta\}] := \frac{1}{Z_{\L}^N} \prod_{x \in \L} \frac{1}{c(\eta_x)!},
\ee
where
$c(0)!:=1$, $c(k)!:=c(1)\cdot\dots\cdot c(k)$, for $k>0$, and $Z_{\L}^N$ is the normalization factor.
The measure $\nuln$ will be referred to as the {\em canonical measure}. Note that the system is
reversible for $\nuln$, i.e. $\LL_{\L}$ is self-adjoint in $L^2(\nuln)$ or, equivalently, the
{\em detailed balance condition}
\be{detbal}
c(\eta_x) \nuln[\{\eta\}] = c(\eta_y +1) \nuln[\{\eta-\d_x + \d_y\}]
\ee
holds for every $x \in \L$ and $\eta \in \O_{\L}$ such that $\eta_x >0$.

Our main result, that is stated next, will be proved under the following conditions.
\begin{co}[LG]
\label{co:LG}
  \begin{displaymath}
    \sup_{k\in\N}|c(k+1)-c(k)|:=a_1<+\infty.
  \end{displaymath}
\end{co}
As remarked in \cite{La:Se:Va} for the spectral gap, $N$-independence
 of the logarithmic-Sobolev constant
requires extra-conditions; in particular, our main result would not hold true in the case
$c(k) = c \ind(k\in \natural \setminus \{0\})$. The following condition, that is the same assumed in
\cite{La:Se:Va}, is a monotonicity requirement on $c(\cdot)$ that rules out the case above.
\begin{co}[M]
\label{co:M}
  There exists $k_0>0$ and $a_2>0$ such that $c(k)-c(j)\geq a_2$ for any $j\in \N$ and $k\geq j+k_0$.
 \end{co}
A simple but key consequence of conditions above is that there exists $A_0>0$ such that 
\be{lineargrowth}
A_0^{-1}k\leq c(k) \leq A_0 k \ \mbox{ for any $k\in\natural$}.
\ee

In what follows, we choose $\L = [0,L]^d \cap \Z^d$. In order to state our main result, we
define the Dirichlet form corresponding to $\LL_{\L}$ and $\nuln$:
\be{not3}
\EE_{\nuln}(f,g) = - \nuln[f\LL_{\L} g] = \frac{1}{2} \sum_{x \in \L}\sum_{y \sim x} \nuln\left[c(\eta_x)\partial_{xy}f(\eta)
\partial_{xy}g(\eta) \right].
\ee
\bt{maint}
Assume that conditions (LG) and (M) hold.
Then there exists a constant $C(L)>0$, that may only depend on $a_1,a_2$, the dimension 
$d$ and $L$, such that
for every choice of $N \geq 1$, $L \geq 2$ and $f: \O_{\L} \ra \R$, $f>0$, we have
\be{not4}
\ent_{\nuln}(f) \leq C(L) \EE_{\nuln}\left(\sqrt{f},\sqrt{f}\right).
\ee
\et

\section{Outline of the proof}

For simplicity, the proof will be outlined in one dimension. The essential steps for the extension to any dimension are analogous to the ones for the spectral gap, that can be 
found in \cite{La:Ki}, Appendix 3.3. However, in the most technical and original estimates contained
in this work (see Sections 7 and 8), proofs are given in a general dimension $d \geq 1$.

\subsection{Step 1: duplication}

The idea is to prove Theorem \ref{maint} by induction on $|\L|$. Suppose $|\L| = 2L$, so that
$\L = \L_1 \cup \L_2$, $|\L_1| = |\L_2| = L$, where $\L_1, \L_2$ are two disjoint  adjacent
segments in $\Z$. By a basic identity on the entropy, we have 
\be{paolo1}
\ent_{\nu_{\L}^N}(f) = \nu_{\L}^N \left[\ent_{\nu_{\L}^N[\cdot | \ov{\eta}_{\L_1}]}(f)  \right]
+\ent_{\nuln}(
\nuln[f|\ov{\eta}_{\L_1}]).
\ee
Note that $\nu_{\L}^N[\cdot | \ov{\eta}_{\L_1}] = \nu^{\ov{\eta}_{\L_1}}_{\L_1} \otimes
\nu^{N-\ov{\eta}_{\L_1}}_{\L_2}$. Thus, by the tensor property of the entropy (see \cite{toulouse},
Th. 3.2.2):
\be{paolo2}
 \nu_{\L}^N \left[\ent_{\nu_{\L}^N[\cdot | \ov{\eta}_{\L_1}]}(f) \right] \leq 
\nuln\left[ \ent_{\nu_{\L_1}^{\ov{\eta}_{\L_1}}}(f) + 
\ent_{\nu_{\L_2}^{N-\ov{\eta}_{\L_1}} }(f) \right].
\ee
Now, let $s(L,N)$ be the maximum of the logarithmic-Sobolev constant for the zero-range process
in  volumes $\L$ with $|\L| \leq L$ and less that $N$ particles, i.e. $s(L,N)$ is the smallest constant
such that
\[
\ent_{\nu_{\L}^n}(f) \leq s(L,N)\EE_{\nu_{\L}^n}(\sqrt{f},\sqrt{f}).
\]
for all $f >0$,  $|\L| \leq L$ and $n \leq N$. Then, by (\ref{paolo2}),
\begin{eqnarray}
\nu_{\L}^N \left[\ent_{\nu_{\L}^N[\cdot | \ov{\eta}_{\L_1}]}(f) \right] & \leq & s(L,N) \nuln 
\left[\EE_{\nu_{\L_1}^{\ov{\eta}_{\L_1}}}(\sqrt{f},\sqrt{f}) +
\EE_{\nu_{\L_2}^{N-\ov{\eta}_{\L_1}} }(\sqrt{f},\sqrt{f}) \right]  \nonumber \\
\addtocounter{for}{1}
& = & s(L,N)\EE_{\nuln}(\sqrt{f},\sqrt{f}) \label{paolo3}.
\end{eqnarray}
Identity (\ref{paolo1}) and inequality 
(\ref{paolo3}) suggest to estimate $s(L,N)$ by induction on $L$. The hardest thing
is to make appropriate estimates on the term $\ent_{\nuln}(
\nuln[f|\ov{\eta}_{\L_1}])$. Note that this term is the entropy of a function depending only on the number of
particles in $\L_1$.

\subsection{Step 2: logarithmic  Sobolev inequality for the distribution of the number of particles in $\boldsymbol{\L_1}$}

Let
\[
\g_{\L}^N(n) = \g(n) := \nuln[\ov{\eta}_{\L_1} =n].
\]
$\g(\cdot)$ is a probability measure on $\{0,1,\ldots,N\}$ that is reversible for the birth and death process
with generator
\be{paolo4}
\A \varphi(n) = \left[\frac{\g(n+1)}{\g(n)} \wedge 1 \right](\varphi(n+1) - \varphi(n)) + \left[\frac{\g(n-1)}{\g(n)} \wedge 1 \right](\varphi(n-1) - \varphi(n))
\ee
and Dirichlet form
\[
\dir(\varphi,\varphi) = - \langle \varphi, \A \varphi \rangle_{L^2(\g)} = \sum_{n=1}^N [\g(n) \wedge \g(n-1)](\varphi(n) - \varphi(n-1))^2.
\]
Logarithmic Sobolev inequalities for birth and death processes are studied in \cite{Mi}. The nontrivial proof that
conditions in \cite{Mi} are satisfied by $\g(n)$, leads to the following result.

\bp{ppaolo1}
The Markov chain with generator (\ref{paolo4}) has a logarithmic Sobolev constant proportional to $N$, i.e.
there exists a constant $C>0$ such that for all $\varphi \geq 0$
\[
\ent_{\g}(\varphi) \leq CN \dir(\sqrt{\varphi}, \sqrt{\varphi}).
\]
\ep

We now apply Proposition \ref{ppaolo1} to the second summand of the r.h.s. of (\ref{paolo1}), and we
obtain 
\[
\ent_{\nuln}(
\nuln[f|\ov{\eta}_{\L_1}])  \leq  CN  \sum_{n=1}^N [\g(n) \wedge \g(n-1)] \left[ \sqrt{\nuln[f|\ov{\eta}_{\L_1} = n]}
- \sqrt{\nuln[f|\ov{\eta}_{\L_1} = n-1]}\right]^2  
\]
\be{paolo5}
 \leq   CN  \sum_{n=1}^N \frac{\g(n) \wedge \g(n-1)}{\nuln[f|\ov{\eta}_{\L_1} = n] \vee \nuln[f|\ov{\eta}_{\L_1} = n-1]}
 \left(\nuln[f|\ov{\eta}_{\L_1} = n]
- \nuln[f|\ov{\eta}_{\L_1} = n-1]\right)^2,
\ee
where we have used the inequality $(\sqrt{x} - \sqrt{y})^2 \leq \frac{(x-y)^2}{x \vee y}$, $x,y>0$.

\subsection{Step 3: study of the term $\boldsymbol{\nuln[f|\ov{\eta}_{\L_1} = n]
- \nuln[f|\ov{\eta}_{\L_1} = n-1]}$}

One of the key points in the proof of Theorem \ref{maint} consists in finding the ``right'' representation for the
discrete gradient $\nuln[f|\ov{\eta}_{\L_1} = n]
- \nuln[f|\ov{\eta}_{\L_1} = n-1]$, that appears in the r.h.s. of (\ref{paolo5}). 

\bp{ppaolo2}
For every $f$ and every $n=1,2,\ldots, N$ we have
\begin{eqnarray}
\nuln[f|\ov{\eta}_{\L_1} = n]
- \nuln[f|\ov{\eta}_{\L_1} = n-1] & = &
\frac{\g(n-1)}{\g(n)} \frac{1}{nL} \left( \nuln\left[\sum_{\stackrel{x \in \L_1}{\scriptscriptstyle{y \in \L_2}}}
h(\eta_x) c(\eta_y) \partial_{yx}f \Bigg| \ov{\eta}_{\L_1} = n-1\right]\right. \nonumber \\ \addtocounter{for}{1}
& & \left. + \nuln \left[f,\sum_{\stackrel{x \in \L_1}{\scriptscriptstyle{y \in \L_2}}}
h(\eta_x)c(\eta_y) \Bigg| \ov{\eta}_{\L_1} = n-1\right]  \right) \label{paolo6},
\end{eqnarray}
where
\[
h(n) := \frac{n +1}{c(n+1)}.
\]
Moreover, by exchanging the roles of $\L_1$ and $\L_2$, the r.h.s. of (\ref{paolo6}) can be
equivalently written as, for every $n-0,1,\ldots,N-1$,
\begin{multline}
-\frac{\g(N-n)}{\g(N-n+1)} \frac{1}{(N-n+1)L} \left( \nuln\left[\sum_{\stackrel{x \in \L_1}{\scriptscriptstyle{y \in \L_2}}}
h(\eta_y) c(\eta_x) \partial_{xy}f \Bigg| \ov{\eta}_{\L_1} = n\right] \right.
\\
+ \left. \nuln \left[f,\sum_{\stackrel{x \in \L_1}{\scriptscriptstyle{y \in \L_2}}}
h(\eta_y) c(\eta_x) \Bigg| \ov{\eta}_{\L_1} = n\right]  \right). 
\addtocounter{for}{1} \label{paolo6b}
\end{multline}
\ep

The representations (\ref{paolo6}) and (\ref{paolo6b}) will be used for $n \geq \frac{N}{2}$ and $n < \frac{N}{2}$
respectively. For convenience, we rewrite (\ref{paolo6}) and (\ref{paolo6b}) as
\be{paolo6c}
\nuln[f|\ov{\eta}_{\L_1} = n]
- \nuln[f|\ov{\eta}_{\L_1} = n-1] =: A(n) + B(n),
\ee
where
\be{paolo7}
A(n) := \left\{
\begin{array}{ll}
\frac{\g(n-1)}{\g(n)} \frac{1}{nL}\nuln\left[\sum_{\stackrel{x \in \L_1}{\scriptscriptstyle{y \in \L_2}}}
h(\eta_x)c(\eta_y) \partial_{yx}f \Bigg| \ov{\eta}_{\L_1} = n-1\right] & \text{for } n \geq \frac{N}{2} \\
-\frac{\g(N-n)}{\g(N-n+1)} \frac{1}{(N-n+1)L} \nuln\left[\sum_{\stackrel{x \in \L_1}{\scriptscriptstyle{y \in \L_2}}}
h(\eta_y) c(\eta_x) \partial_{xy}f \Bigg| \ov{\eta}_{\L_1} = n\right] & \text{for } n<\frac{N}{2},
\end{array}
\right.
\ee
and
\be{paolo8}
B(n):= \left\{
\begin{array}{ll}
\frac{\g(n-1)}{\g(n)} \frac{1}{nL}\nuln \left[f,\sum_{\stackrel{x \in \L_1}{\scriptscriptstyle{y \in \L_2}}}
h(\eta_x)c(\eta_y) \Bigg| \ov{\eta}_{\L_1} = n-1\right]  & \text{for } n \geq \frac{N}{2} \\
-\frac{\g(N-n)}{\g(N-n+1)} \frac{1}{(N-n+1)L}\nuln \left[f,\sum_{\stackrel{x \in \L_1}{\scriptscriptstyle{y \in \L_2}}}
h(\eta_y) c(\eta_x) \Bigg| \ov{\eta}_{\L_1} = n\right]  & \text{for } n<\frac{N}{2}.
\end{array}
\right.
\ee

Thus, our next aim is to get estimates on the two terms in the r.h.s. of (\ref{paolo6c}). It is useful to stress that
the two terms are qualitatively different. Estimates on $A(n)$ are essentially insensitive to the precise form of $c(\cdot)$. Indeed,
the dependence of $A(n)$ on $L$ and $N$ is of the same order as in the case $c(n) = \l n$, i.e. the case of independent
particles. Quite differently, the term $B(n)$
vanishes in the case of independent particles, since, in that case, the term $\sum_{\stackrel{x \in \L_1}{\scriptscriptstyle{y \in \L_2}}}
h(\eta_x)c(\eta_y)$ is a.s. constant with respect to $ \nuln \left[\cdot |\ov{\eta}_{\L_1} = n-1\right]$. Thus, $B(n)$ 
somewhat depends on interaction between particles. Note that our model is not necessarily a ``small
perturbation'' of a system of independent particles; there is no small parameter in the model that
guarantees that $B(n)$ is small enough. Essentially all technical results of this paper are concerned with
estimating $B(n)$.

\subsection{Step 4: estimates on $\boldsymbol{A(n)}$}

The following proposition gives the key estimate on $A(n)$.
\bp{ppaolo3}
There is a constant $C>0$ such that
\begin{multline*}
  A^2(n) \leq \frac{CL^2}{N}\left(\nuln[f|\ov{\eta}_{\L_1} = n] \vee \nuln[f|\ov{\eta}_{\L_1} = n-1]\right)\\ 
\times\left[\frac{\gamma(n-1)}{\gamma(n)}\EE_{\nuln[\cdot|\ov{\eta}_{\L_1}= n-1]}(\sqrt{f},\sqrt{f})+\EE_{\nuln
[\cdot|\ov{\eta}_{\L_1}= n]}(\sqrt{f},\sqrt{f})\right].
\end{multline*}
\ep

\begin{remark}

Let us try so see where we are now. Let us ignore, for the moment the term $B(n)$, i.e. let
us pretend that $B(n) \equiv 0$. Thus, by (\ref{paolo6c}) and Proposition \ref{ppaolo3} we
would have

\begin{multline}
 \left(\nuln[f|\ov{\eta}_{\L_1} = n]
 - \nuln[f|\ov{\eta}_{\L_1} = n-1]\right)^2 
 \leq \frac{CL^2}{N}\left(\nuln[f|\ov{\eta}_{\L_1} = n] \vee \nuln[f|\ov{\eta}_{\L_1} = n-1]\right)\\ 
\times\left[\frac{\gamma(n-1)}{\gamma(n)}\EE_{\nuln[\cdot|\ov{\eta}_{\L_1}= n-1]}(\sqrt{f},\sqrt{f})+\EE_{\nuln
[\cdot|\ov{\eta}_{\L_1}= n]}(\sqrt{f},\sqrt{f})\right]
 \addtocounter{for}{1} \label{paolo9}
\end{multline}

Inserting (\ref{paolo9}) into (\ref{paolo5}) we get, for some possibly different constant $C$,
\be{paolo9b}
\ent_{\nuln}(
\nuln(f|\ov{\eta}_{\L_1})) \leq CL^2 \EE_{\nuln}(\sqrt{f},\sqrt{f}),
\ee
where we have used the obvious identity:
\be{paolo10}
\nuln\left[\EE_{\nuln[\cdot | {\cal{G}}]} (\sqrt{f},\sqrt{f}) \right] =  \EE_{\nuln}(\sqrt{f},\sqrt{f}),
\ee
that holds for any $\s$-field ${\cal{G}}$.
Inequality (\ref{paolo9b}), together with (\ref{paolo1}) and (\ref{paolo3}) yields
\be{paolo11}
s(2L,N) \leq s(L,N) + CL^2.
\ee
Thus, if we can show that
\be{paolo12}
\sup_N s(2,N) < +\infty,
\ee
(see Proposition \ref{ppaolo6} next) then we would have
\[
\sup_N s(L,N) \leq CL^2
\]
for some $C>0$, i.e. we get the exact order of growth of $s(L,N)$.
 In all this, however, we have totally 
ignored the contribution of $B(n)$.

\end{remark}

\subsection{Step 5: preliminary analysis of $\boldsymbol{B(n)}$}

We confine ourselves to the analysis of $B(n)$ for $n \geq \frac{N}{2}$, since the case $n < \frac{N}{2}$ is identical. Consider
the covariance that appears in the r.h.s. of the first formula of (\ref{paolo8}).
By elementary properties of the covariance and the fact that $\nu_{\L}^N[\cdot | \ov{\eta}_{\L_1}] = \nu^{\ov{\eta}_{\L_1}}_{\L_1} \otimes
\nu^{N-\ov{\eta}_{\L_1}}_{\L_2}$, we get
\begin{multline}
\nuln \left[\left. f,\sum_{\stackrel{x \in \L_1}{\scriptscriptstyle{y \in \L_2}}}
h(\eta_x)c(\eta_y) \right| \ov{\eta}_{\L_1} = n-1\right]
= \nu_{\L_1}^{n-1} \left[ \sum_{x \in \L_1} h(\eta_x)\nu_{\L_2}^{N-n+1}\left[f, \sum_{y \in \L_2} c(\eta_y)\right] \right]  
\\
+ \nu_{\L_1}^{n-1}\left[ \nu_{\L_2}^{N-n+1}[f], \sum_{x \in \L_1} h(\eta_x)\right] \nu_{\L_2}^{N-n+1}\left[
\sum_{y \in \L_2} c(\eta_y) \right]. 
\addtocounter{for}{1} \label{paolo13}
\end{multline}
It follows by Conditions (LG) and (M) (see  (\ref{lineargrowth})) that, for some constant $C>0$, $h(\eta_x)\leq C$
and $c(\eta_y) \leq C \eta_y$. Thus, a simple estimate on the two summands in (\ref{paolo13}),
yields, for 
some $C>0$
\begin{multline}
B^2(n) \leq \frac{\g^2(n-1)}{\g^2(n)} \left(\frac{C}{n^2}  \nu_{\L_2}^{N-n+1} \left[\nu_{\L_1}^{n-1}[f], \sum_{y \in \L_2} c(\eta_y)\right]^2 + \right. \\
+ \left.
\frac{C(N-n+1)^2}{n^2 L^2} \nu_{\L_1}^{n-1}\left[ \nu_{\L_2}^{N-n+1}[f], \sum_{x \in \L_1} h(\eta_x)\right]^2 \right). \label{paolo14}
\end{multline}
Thus, our next aim is to estimate the two covariances in (\ref{paolo14}). 

\subsection{Estimates on $\boldsymbol{B(n)}$: entropy inequality and estimates on moment generating functions}

By (\ref{paolo14}), estimating $B^2(n)$ consists in estimating two covariances. In general, covariances can
be estimated by the following {\em entropy inequality}, that holds for every probability measure
$\nu$ (see \cite{toulouse}, Section 1.2.2):
\be{paolo16}
\nu[f,g] = \nu[f(g-\nu[g])] \leq \frac{\nu[f]}{t} \log \nu\left[ e^{t(g-\nu[g])} \right] + \frac{1}{t}\ent_{\nu}(f),
\ee
where $f \geq 0$ and $t>0$ is arbitrary. Since in (\ref{paolo14}) we need to estimate the square of a covariance,
we write (\ref{paolo16}) with $-g$ in place of $g$, and obtain
\be{paolo17}
|\nu[f,g]| \leq \frac{\nu[f]}{t} \log\left( \nu\left[ e^{t(g-\nu[g])}\right] \vee \nu\left[e^{-t(g-\nu[g])}\right]\right) + \frac{1}{t}\ent_{\nu}(f).
\ee
Therefore, we first get estimates on the moment generating functions $ \nu\left[ e^{\pm t(g-\nu[g])}\right]$, and then optimize
(\ref{paolo17}) over $t>0$. 

Note that the covariances in (\ref{paolo14}) involve
 functions of 
$\eta_{\L_1}$ or $\eta_{\L_2}$. In next two propositions we write $\L$ for $\L_1$ and $\L_2$, and
denote by $N$ the number of particles in $\L$. Their proof can be found in \cite{Da:Po}.
\bp{ppaolo4}
Let $x \in \L$. Then there is a constant $A_L$ depending on
$L = |\L|$ such that for every $N> 0$ and $t \in [-1,1]$
\be{paolo18}
\nuln \left[e^{t\left(c(\eta_x) - \nuln[c(\eta_x)]\right)}\right] \leq e^{A_L N t^2 },
\ee
and
\be{paolo19}
\nuln \left[e^{tN \left(h(\eta_x)- \nuln[h(\eta_x)]\right)}\right] \leq A_L e^{A_L (N t^2 + \sqrt{N}|t|)},
\ee
\ep
Using (\ref{paolo16}), (\ref{paolo17}) and (\ref{paolo19}) and optimizing over $t>0$, we get
estimates on the covariances appearing in (\ref{paolo14}). 
\bp{ppaolo5}
There exists a constant $C_L$ depending on $L= |\L|$ such that the following inequalities
hold:
\be{paolo20}
\nuln\left[f,\sum_{x \in \L} c(\eta_x)\right]^2 \leq C_L N \nuln[f] \ent_{\nuln}(f),
\ee
\be{paolo21}
\nuln\left[f,\sum_{x \in \L} h(\eta_x)\right]^2 \leq C_L N^{-1} \nuln[f] \left[ \nuln[f] + \ent_{\nuln}(f)\right].
\ee
\ep

Inserting these new estimates in (\ref{paolo14}) we obtain, for some possibly different $C_L$,
\be{paolo22}
\begin{array}{r}
B^2(n) \leq \frac{C_L \g^2(n-1)}{\g^2(n)} \nuln[f|\ov{\eta}_{\L_1} =n-1] \Bigg( 
\frac{N-n+1}{n^2}  \nu_{\L_1}^{n-1}\left[\ent_{\nu_{\L_2}^{N-n+1}}(f)\right] 
\\
 +  \frac{(N-n+1)^2 }{n^3}\Bigg( 
\nuln[f|\ov{\eta}_{\L_1} =n-1] + \nu_{\L_2}^{N-n+1}\left[\ent_{\nu_{\L_1}^{n}}(f)\right] \Bigg)
\Bigg).
\end{array}
\ee
In order to simplify (\ref{paolo22}), we use Proposition 
\ref{pro:dec},  which gives
\be{paolo22bis}
\frac{n}{C(N-n+1)} \leq \frac{\g(n-1)}{\g(n)} \leq \frac{Cn}{N-n+1}
\ee
for some $C>0$. It follows that
\[
\frac{\g^2(n-1)}{\g^2(n)}\frac{N-n+1}{n^2} \leq \frac{C^2}{N-n+1} \leq \frac{\g(n-1)}{\g(n)} \frac{C^3}{n},
\]
and
\[
\frac{\g^2(n-1)}{\g^2(n)} \frac{(N-n+1)^2 }{n^3} \leq \frac{C^2}{n}.
\]
Thus, (\ref{paolo22}) implies, for some $C_L >0$ depending on $L$, recalling also that $n \geq \frac{N}{2}$,
\begin{multline}
N \frac{\g(n) \wedge \g(n-1)}{\nuln[f|\ov{\eta}_{\L_1} = n] \vee \nuln[f|\ov{\eta}_{\L_1} = n-1]} B^2(n) \leq C_L \g(n-1)\Big(\nuln[f|\ov{\eta}_{\L_1} =n-1]  
 \\ + 
\nu_{\L_1}^{n-1}\left[\ent_{\nu_{\L_2}^{N-n+1}}(f)\right] + \nu_{\L_2}^{N-n+1}\left[\ent_{\nu_{\L_1}^{n}}(f)\right] \Big).
\label{paolo22ter}
\end{multline}
Now, we bound the two terms 
\[
\ent_{\nu^{n-1}_{\L_1}}(f) \ \text{ and }
\ \ent_{\nu^{N-n+1}_{\L_2}}(f)
\]
 by, respectively, 
\[
s(N,L) \EE_{\nu_{\L_1}^{n-1}}(\sqrt{f},\sqrt{f}) \ 
\text{ and } \ s(N,L) \EE_{\nu_{\L_2}^{N-n+1}}(\sqrt{f},\sqrt{f}), 
\]
and insert these estimates in 
(\ref{paolo6}). What comes out is then used to estimate (\ref{paolo5}), after 
having obtained the corresponding estimates for $n<\frac{N}{2}$.  Recalling the estimates for $A(n)$, 
straightforward computations yield
\be{paolo23}
\ent_{\nuln}(
\nuln[f|\ov{\eta}_{\L_1}]) \leq C_L\EE_{\nuln}(\sqrt{f},\sqrt{f}) + C_L\nuln[f] + 
C_Ls(N,L)\EE_{\nuln}(\sqrt{f},\sqrt{f}),
\ee
for some $C_L>0$ depending on $L$. Inequality (\ref{paolo23}), together with (\ref{paolo3}), gives, for
a possibly different constant $C_L>0$
\be{paolo24}
\ent_{\nuln}(f) \leq C_L\EE_{\nuln}(\sqrt{f},\sqrt{f}) + C_L\nuln[f] + 
C_Ls(N,L)\EE_{\nuln}(\sqrt{f},\sqrt{f}).
\ee
To deal with the term $\nuln[f]$ in (\ref{paolo24}) we use the following well known
argument. Set $\ov{f} = \left(\sqrt{f} - \nuln[\sqrt{f}]\right)^2$. By Rothaus inequality (see \cite{toulouse}, Lemma 4.3.8)
\[
\ent_{\nuln}(f) \leq \ent_{\nuln}(\ov{f}) + 2 \nuln[\sqrt{f},\sqrt{f}].
\]
Using this inequality and replacing $f$ by $\ov{f}$ in (\ref{paolo24}) we get, for a different $C_L$,
\begin{eqnarray}
\ent_{\nuln}(f) & \leq & C_L\EE_{\nuln}(\sqrt{f},\sqrt{f}) + C_L \nuln[\sqrt{f},\sqrt{f}] + 
C_Ls(N,L)\EE_{\nuln}(\sqrt{f},\sqrt{f}) \nonumber \\
& \leq & D_L\EE_{\nuln}(\sqrt{f},\sqrt{f}) + D_Ls(N,L)\EE_{\nuln}(\sqrt{f},\sqrt{f})
\addtocounter{for}{1} \label{paolo25}
\end{eqnarray}
where, in the last line, we have used the Poincar\'e inequality (\ref{gap}). Therefore 
\be{paolo26}
s(2L,N) \leq D_L [s(L,N) +1],
\ee
that implies  (\ref{not4}), provided we prove the following ``basis step'' for the induction.
\bp{ppaolo6}
\[
\sup_N s(2,N) < +\infty
\]
\ep
The proof of Proposition \ref{ppaolo6} is also given in Section 8.
As we pointed out above, (\ref{paolo26}) gives no indication on how $s(L) = \sup_N s(N,L)$ grows 
with $L$. The proof of the actual $L^2$-growth is given in \cite{Da:Po}.

\section{Representation of
$\boldsymbol{\nu_\Lambda^N[f|\n{\eta}{\Lambda_1} = n]-\nu_\Lambda^N[f|\n{\eta}{\Lambda_1} = n-1]}$: proof of Proposition~\ref{ppaolo2}}

We begin with a simple consequence of reversibility

\begin{lem}
  \label{lemma:rev}
  Let $\Lambda=\Lambda_1\cup\Lambda_2$ be a partition of $\Lambda$, $N\in\natural$.
  Define $\gamma(n):=\nu_\Lambda^N[\n\eta{\Lambda_1}=n]$ and $\tau_{xy} f:=\partial_{xy}f +f$.
  Then:
  \begin{displaymath}
    \nu_{\Lambda}^{N}[f\ind(\eta_x>0)|\n\eta{\Lambda_1}=n]=
    \begin{cases}
     \frac{\gamma(n-1)}{\gamma(n)}\nu_{\Lambda}^{N}\left[\left.\frac{c(\eta_y)}{c(\eta_x+1)}\tau_{yx}f\right|\n\eta{\Lambda_1}=n-1\right] & \text{for any $x\in\Lambda_1$ and $y\in\Lambda_2$}\\
     \\
     \nu_{\Lambda}^{N}\left[\left.\frac{c(\eta_y)}{c(\eta_x+1)}\tau_{yx}f\right|\n\eta{\Lambda_1}=n\right] & \text{for any $x,y\in\Lambda_2$}\\
    \end{cases}
  \end{displaymath}
   for any $f\in L^1(\nu_\Lambda^N)$ and $n\in\{1,\dots,N\}$.
\end{lem}
\begin{proof}
  Assume that $y\in\Lambda_2$.
  We use first the detailed balance condition \eqref{detbal}, then the change of variables $\xi\mapsto\xi+\delta_y-\delta_x$ to obtain
  \begin{multline*}
    \nu_{\Lambda}^{N}\left[f\ind(\eta_x>0)|\n\eta{\Lambda_1}=n\right]=
    \frac{\sum_{\xi:\xi_x>0}\nu_\Lambda^N[\eta=\xi]f(\xi)\ind(\n\xi{\Lambda_1}=n)}{\nu_\Lambda^N[\bar\eta_{\Lambda_1}=n]}\\
    =\frac{\sum_{\xi:\xi_x>0}\nu_\Lambda^N[\eta=\xi-\delta_x+\delta_y]\frac{c(\xi_y+1)}{c(\xi_x)}f(\xi)\ind(\n\xi{\Lambda_1}=n)}{\nu_\Lambda^N[\n\eta{\Lambda_1}=n]}\\
        =
        \begin{cases}
     \frac{\gamma(n-1)}{\gamma(n)}\nu_{\Lambda}^{N}\left[\left.\frac{c(\eta_y)}{c(\eta_x+1)}\tau_{yx}f\right|\n\eta{\Lambda_1}=n-1\right] & \text{if $x\in\Lambda_1$}\\
     \\
     \nu_{\Lambda}^{N}\left[\left.\frac{c(\eta_y)}{c(\eta_x+1)}\tau_{yx}f\right|\n\eta{\Lambda_1}=n\right] & \text{if $x\in\Lambda_2$}\\
    \end{cases}
  \end{multline*}
\end{proof}

\begin{lem}
\label{pro:1}
  Let $\Lambda=\Lambda_1\cup\Lambda_2$ be a partition of $\Lambda$, $N\in\natural$.
  Define $\gamma(n):=\nu_\Lambda^N[\n\eta{\Lambda_1}=n]$ and $h:\natural\to\real$ by $h(n):=(n+1)/c(n+1)$.
  Then:
  \begin{multline}
    \label{eq:pro:1:1}
    \nu_{\Lambda}^{N}[f|\n\eta{\Lambda_1}=n]-\nu_{\Lambda}^{N}[f|\n\eta{\Lambda_1}=n-1]\\
    =\frac{\gamma(n-1)}{n\gamma(n)}
    \left(\nu_{\Lambda}^{N}\left[c(\eta_y)\sum_{x\in\Lambda_1}h(\eta_x)\partial_{yx}f\Bigg|\n\eta{\Lambda_1}=n-1\right]+\nu_{\Lambda}^{N}\left[f,c(\eta_y)\sum_{x\in\Lambda_1}h(\eta_x)\Bigg|\n\eta{\Lambda_1}=n-1\right]\right),
  \end{multline}
  for any $f\in L^1(\nu_\Lambda^N)$, $y\in\Lambda_2$ and $n\in\{1,\dots,N\}$.
\end{lem}
\begin{proof}
  We begin by writing
\[
    \nu_{\Lambda}^{N}[f|\n\eta{\Lambda_1}=n]=
     -\frac{1}{n}\sum_{x\in\Lambda_1}\nu_{\Lambda}^{N}\left[\eta_x \partial_{xy}f\big|\n\eta{\Lambda_1}=n\right]+ \frac{1}{n}\sum_{x\in\Lambda_1}\nu_{\Lambda}^{N}\left[\eta_x \tau_{xy}f\big|\n\eta{\Lambda_1}=n\right].
\]
  By Lemma~\ref{lemma:rev} we get
  \begin{multline*}
    -\nu_{\Lambda}^{N}\left[\eta_x \partial_{xy}f\big|\n\eta{\Lambda_1}=n\right]
    =\frac{\gamma(n-1)}{\gamma(n)}\nu_{\Lambda}^{N}\left[\frac{c(\eta_y)}{c(\eta_x+1)}(\eta_x+1)\partial_{yx}f\Bigg|\n\eta{\Lambda_1}=n-1\right]\\
    =\frac{\gamma(n-1)}{\gamma(n)}\nu_{\Lambda}^{N}\left[c(\eta_y)h(\eta_x)\partial_{yx}f\big|\n\eta{\Lambda_1}=n-1\right],
  \end{multline*}
  and similarly
  \begin{displaymath}
    \nu_{\Lambda}^{N}\left[\eta_x \tau_{xy}f\big|\n\eta{\Lambda_1}=n\right]
    =\frac{\gamma(n-1)}{\gamma(n)}\nu_{\Lambda}^{N}\left[c(\eta_y)h(\eta_x) f\big|\n\eta{\Lambda_1}=n-1\right].
  \end{displaymath}
  Thus
  \begin{multline*}
     \nu_{\Lambda}^{N}[f|\n\eta{\Lambda_1}=n]\\
    =\frac{\gamma(n-1)}{n\gamma(n)}
    \left(\nuln\left[c(\eta_y)\sum_{x\in\Lambda_1}h(\eta_x)\partial_{yx}f\Bigg|\n\eta{\Lambda_1}=n-1\right]+\nu_{\Lambda}^{N}\left[c(\eta_y)f \sum_{x\in\Lambda_1}h(\eta_x) \Bigg|\n\eta{\Lambda_1}=n-1\right]\right).
   \end{multline*}
   Letting $f\equiv1$ in this last formula we obtain
   \begin{displaymath}
     \frac{\gamma(n-1)}{n\gamma(n)}\nu_{\Lambda}^{N}\left[c(\eta_y) \sum_{x\in\Lambda_1}h(\eta_x) \Bigg|\n\eta{\Lambda_1}=n-1\right]=1
   \end{displaymath}
   that implies
   \begin{multline*}
     \frac{\gamma(n-1)}{n\gamma(n)}\nu_{\Lambda}^{N}\left[c(\eta_y)f \sum_{x\in\Lambda_1}h(\eta_x) \Bigg|\n\eta{\Lambda_1}=n-1\right]\\
    =\nu_{\Lambda}^{N}\left[f|\n\eta{\Lambda_1}=n-1\right] +
     \frac{\gamma(n-1)}{n\gamma(n)}\nu_{\Lambda}^{N}\left[f, c(\eta_y) \sum_{x\in\Lambda_1}h(\eta_x) \Bigg|\n\eta{\Lambda_1}=n-1\right]
    \end{multline*}
   and, therefore,
   \begin{multline*}
     \nu_{\Lambda}^{N}[f|\n\eta{\Lambda_1}=n]-\nu_{\Lambda}^{N}\left[f|\n\eta{\Lambda_1}=n-1\right] \\
     =\frac{\gamma(n-1)}{n\gamma(n)}
    \left(\nuln\left[c(\eta_y)\sum_{x\in\Lambda_1}h(\eta_x)\partial_{yx}f\Bigg|\n\eta{\Lambda_1}=n-1\right]+\nu_{\Lambda}^{N}\left[f,c(\eta_y) \sum_{x\in\Lambda_1}h(\eta_x) \Bigg|\n\eta{\Lambda_1}=n-1\right]\right).
   \end{multline*}
\end{proof}

\begin{proofof}{Proposition~\ref{ppaolo2}}
  Equation \eqref{paolo6} is obtained by \eqref{eq:pro:1:1} by averaging over $y\in\Lambda_2$.
\end{proofof}

\section{Bounds on $\boldsymbol{A(n)}$: proof of Proposition~\ref{ppaolo3}}

\begin{proofof}{Proposition~\ref{ppaolo3}}
  The Proposition is proved if we can show that for $1\leq n< N/2$
  \begin{multline*}
      A^2(n) \leq \frac{CL^2}{N-n}\left(\nuln[f|\ov{\eta}_{\L_1} = n] \vee \nuln[f|\ov{\eta}_{\L_1} = n-1]\right)\\ 
\times\left[\frac{\gamma(n-1)}{\gamma(n)}\EE_{\nuln[\cdot|\ov{\eta}_{\L_1}= n-1]}(\sqrt{f},\sqrt{f})+\EE_{\nuln
[\cdot|\ov{\eta}_{\L_1}= n]}(\sqrt{f},\sqrt{f})\right],
  \end{multline*}
  while for $n\geq N/2$
  \begin{multline}
    \label{eq:1000}
      A^2(n) \leq \frac{CL^2}{n}\left(\nuln[f|\ov{\eta}_{\L_1} = n] \vee \nuln[f|\ov{\eta}_{\L_1} = n-1]\right)\\ 
\times\left[\frac{\gamma(n-1)}{\gamma(n)}\EE_{\nuln[\cdot|\ov{\eta}_{\L_1}= n-1]}(\sqrt{f},\sqrt{f})+\EE_{\nuln
[\cdot|\ov{\eta}_{\L_1}= n]}(\sqrt{f},\sqrt{f})\right].
  \end{multline}
  We will prove only this last bound being the proof of the previous one identical.
  So assume $n\geq N/2$ and notice that $\partial_{yx}f=(\partial_{yx}\sqrt{f})(\sqrt{f}+\tau_{yx}\sqrt{f})$.
  By Cauchy-Schwartz inequality
  \begin{multline}
    \label{eq:pro:2:1}
    \left\{\frac{\gamma(n-1)}{n L\gamma(n)}\nu_\Lambda^N\Bigg[\sum_{x\in\Lambda_1,\, y\in\Lambda_2}c(\eta_y)h(\eta_x)\partial_{yx}f\Bigg|\n\eta{\Lambda_1}=n-1\Bigg]\right\}^2\\
    \leq\left[\frac{\gamma(n-1)}{n L\gamma(n)}\right]^2
    \nu_\Lambda^N\Bigg[\sum_{x\in\Lambda_1,\, y\in\Lambda_2}c(\eta_y)h(\eta_x)\left(\partial_{yx}\sqrt{f}\right)^2\Bigg|\n\eta{\Lambda_1}=n-1\Bigg]\\
   \times\nu_\Lambda^N\Bigg[\sum_{x\in\Lambda_1,\,y\in\Lambda_2}c(\eta_y)h(\eta_x)\left(\sqrt{f}+\tau_{yx}\sqrt{f}\right)^2\Bigg|\n\eta{\Lambda_1}=n-1\Bigg]\\
    \leq 2 a^2\left[\frac{\gamma(n-1)}{n L\gamma(n)}\right]^2\nu_\Lambda^N\Bigg[\sum_{x\in\Lambda_1,\, y\in\Lambda_2}c(\eta_y)\left(\partial_{yx}\sqrt{f}\right)^2\Bigg|\n\eta{\Lambda_1}=n-1\Bigg]\\
   \times\nu_\Lambda^N\Bigg[\sum_{x\in\Lambda_1,\,y\in\Lambda_2}c(\eta_y)\left(f+\tau_{yx}f\right)\Bigg|\n\eta{\Lambda_1}=n-1\Bigg]
  \end{multline}
  where $a:=\Vert h \Vert_{+\infty}$ ((\ref{lineargrowth}) implies boundedness of $h$).
  In order to bound the last factor in \eqref{eq:pro:2:1} we
  observe that by Lemma~\ref{lemma:rev}
  \begin{displaymath}
    \nu_\Lambda^N\left[c(\eta_y)\tau_{yx}f\big|\n\eta{\Lambda_1}=n-1\right]=
    \frac{\gamma(n)}{\gamma(n-1)}\nu_\Lambda^N\left[c(\eta_x)f\big|\n\eta{\Lambda_1}=n\right],
  \end{displaymath}
  so that
  \begin{multline*}
    \nu_\Lambda^N\Bigg[\sum_{x\in\Lambda_1,\,y\in\Lambda_2}c(\eta_y)\left(f+\tau_{yx}f\right)\Bigg|\n\eta{\Lambda_1}=n-1\Bigg]\\
    =\sum_{x\in\Lambda_1,\, y\in\Lambda_2}\left(\nu_\Lambda^N\left[c(\eta_y)f\big|\n\eta{\Lambda_1}=n-1\right]+
    \frac{\gamma(n)}{\gamma(n-1)}\nu_\Lambda^N\left[c(\eta_x)f\big|\n\eta{\Lambda_1}=n\right]\right)\\
  \leq a\left\{(N-n+1)L\nu_\Lambda^N\left[f\big|\n\eta{\Lambda_1}=n-1\right]+
    \frac{\gamma(n)nL}{\gamma(n-1)}\nu_\Lambda^N\left[f\big|\n\eta{\Lambda_1}=n\right]\right\},
  \end{multline*}
  where the fact $c(k)\leq ak$ (see (\ref{lineargrowth})) has been used to obtain the last line.
  This implies that
    \begin{multline*}
    \frac{\gamma(n-1)}{\gamma(n)nL}\nu_\Lambda^N\Bigg[\sum_{x\in\Lambda_1,\,y\in\Lambda_2}c(\eta_y)\left(f+\tau_{yx}f\right)\Bigg|\n\eta{\Lambda_1}=n-1\Bigg]\\
    \leq a\left\{\frac{\gamma(n-1)(N-n+1)L}{\gamma(n)nL}\nu_\Lambda^N\left[f\big|\n\eta{\Lambda_1}=n-1\right]+
    \nu_\Lambda^N\left[f\big|\n\eta{\Lambda_1}=n\right]\right\}\\
  \leq B_1\left(\nu_\Lambda^N\left[f\big|\n\eta{\Lambda_1}=n-1\right]\vee\nu_\Lambda^N\left[f\big|\n\eta{\Lambda_1}=n\right]\right)
  \end{multline*}
  where in last step we used Proposition~\ref{pro:dec}.
  By plugging this bound in \eqref{eq:pro:2:1} we get
    \begin{multline}
          \label{eq:pro:2:2}
    \left\{\frac{\gamma(n-1)}{n L\gamma(n)}\nu_\Lambda^N\Bigg[\sum_{x\in\Lambda_1,\, y\in\Lambda_2}c(\eta_y)h(\eta_x)\partial_{yx}f\Bigg|\n\eta{\Lambda_1}=n-1\Bigg]\right\}^2\\
    \leq B_2\left(\nu_\Lambda^N\left[f\big|\n\eta{\Lambda_1}=n-1\right]\vee\nu_\Lambda^N\left[f\big|\n\eta{\Lambda_1}=n\right]\right)\\
    \times\frac{\gamma(n-1)}{n L\gamma(n)}
    \nu_\Lambda^N\Bigg[\sum_{x\in\Lambda_1,\, y\in\Lambda_2}c(\eta_y)h(\eta_x)\left(\partial_{yx}\sqrt{f}\right)^2\Bigg|\n\eta{\Lambda_1}=n-1\Bigg].
  \end{multline}
  Now we have to bound the last factor in the right hand side of \eqref{eq:pro:2:2}.
  For $x,y\in\integer$ the \emph{path between $x$ and $y$} will be the sequence of nearest-neighbor integers $\gamma(x,y)=\{z_0,z_1,\dots,z_r\}$ of $\integer$ such that for $i=1,\dots,r$ $|z_{i-1}-z_{i}|=1$, for $i\neq j$ $z_i\neq j$, $z_0=x$ and $z_r=y$.
  Obviously, for $x,y \in \L$, we have $r=|\gamma(x,y)|\leq 2L$.
  Now let such a $\gamma(x,y)$ be given.
  Notice that if $\eta_y\geq1$, we can write
  \begin{displaymath}
    (\partial_{yx}\sqrt{f})(\eta)=
     \sum_{k=0}^{r-1}\left(\partial_{z_{k+1}z_k}\sqrt{f}\right)(\eta-\delta_{z_r}+\delta_{z_{k+1}}).
  \end{displaymath}
  By Jensen inequality, we obtain
  \begin{displaymath}
    c(\eta_y)\left(\partial_{yx}\sqrt{f}\right)^2(\eta)\leq
    2Lc(\eta_y)\sum_{k=0}^{r-1}\left(\partial_{z_{k+1}z_k}\sqrt{f}\right)^2(\eta-\delta_{z_r}+\delta_{z_{k+1}}),
  \end{displaymath}
  so
  \begin{multline*}
    \frac{\gamma(n-1)}{\gamma(n)nL}\nu_{\Lambda}^{N}\left[\sum_{x\in\Lambda_1,y\in\Lambda_2}c(\eta_y)h(\eta_x)\left(\partial_{yx}\sqrt{f}\right)^2\Bigg|\n\eta{\Lambda_1}=n-1\right]\\
    \leq \frac{2\gamma(n-1)}{\gamma(n)n}\sum_{x\in\Lambda_1,y\in\Lambda_2}\sum_{k=0}^{r-1}\nu_{\Lambda}^{N}\left[c(\eta_y)h(\eta_x)\left(\tau_{z_r z_{k+1}}\partial_{z_{k+1}z_k}\sqrt{f}\right)^2\Bigg|\n\eta{\Lambda_1}=n-1\right]\\
    \leq \frac{2a\gamma(n-1)}{\gamma(n)n}\sum_{x\in\Lambda_1,y\in\Lambda_2}\sum_{k=0}^{r-1}\nu_{\Lambda}^{N}\left[c(\eta_y)\left(\tau_{y z_{k+1}}\partial_{z_{k+1}z_k}\sqrt{f}\right)^2\Bigg|\n\eta{\Lambda_1}=n-1\right],
  \end{multline*}
  where again $a=\Vert h \Vert_{+\infty}$.
  By Lemma~\ref{lemma:rev}
  \begin{multline*}
    \nu_{\Lambda}^{N}\left[c(\eta_y)\left(\tau_{y z_{k+1}}\partial_{z_{k+1}z_k}\sqrt{f}\right)^2\Bigg|\n\eta{\Lambda_1}=n-1\right] \\
    =
    \begin{cases}
      \frac{\gamma(n)}{\gamma(n-1)}\nu_{\Lambda}^{N}\left[\left. c(\eta_{z_{k+1}})\left(\partial_{z_{k+1}z_k}\sqrt{f}\right)^2\right|\n\eta{\Lambda_1}=n\right] & \text{if $z_{k+1}\in\Lambda_1$}\\
      \\
      \nu_{\Lambda}^{N}\left[\left. c(\eta_{z_{k+1}})\left(\partial_{z_{k+1}z_k}\sqrt{f}\right)^2\right|\n\eta{\Lambda_1}=n-1\right] & \text{if $z_{k+1}\in\Lambda_2$}
    \end{cases}
  \end{multline*}
  which implies
  \begin{multline*}
    \frac{\gamma(n-1)}{nL\gamma(n)}\nu_{\Lambda}^{N}\left[\sum_{x\in\Lambda_1,y\in\Lambda_2}c(\eta_y)h(\eta_x)\left(\partial_{yx}\sqrt{f}\right)^2\Bigg|\n\eta{\Lambda_1}=n-1\right]\\
    \leq \frac{2a}{n}\sum_{x\in\Lambda_1,y\in\Lambda_2}\Bigg\{\sum_{k:z_{k+1}\in\Lambda_1}\nu_{\Lambda}^{N}\left[c(\eta_{z_{k+1}})\left(\partial_{z_{k+1}z_k}\sqrt{f}\right)^2\Bigg|\n\eta{\Lambda_1}=n-1\right]+\\
\frac{\gamma(n-1)}{\gamma(n)}\sum_{k:z_{k+1}\in\Lambda_2}\nu_{\Lambda}^{N}\left[c(\eta_{z_{k+1}})\left(\partial_{z_{k+1}z_k}\sqrt{f}\right)^2\Bigg|\n\eta{\Lambda_1}=n-1\right]\Bigg\} \\
    \leq\frac{2aL^2}{n}\sum_{\substack{x,y\in\Lambda \\ |x-y|=1}}\Bigg\{\nu_\Lambda^N\left[c(\eta_x)\left(\partial_{xy}\sqrt{f}\right)^2\Bigg|\n\eta{\Lambda_1}=n\right]+\\
    \frac{\gamma(n-1)}{\gamma(n)}\nu_\Lambda^N\left[c(\eta_x)\left(\partial_{xy}\sqrt{f}\right)^2\Bigg|\n\eta{\Lambda_1}=n-1\right]\Bigg\}.
  \end{multline*}
   By plugging this bound in \eqref{eq:pro:2:2} we get \eqref{eq:1000}.
\end{proofof}

\section{Rough estimates: proofs of Propositions \ref{ppaolo4} and \ref{ppaolo5}}

\proofof{Proposition \ref{ppaolo4}}
We begin with  the proof of (\ref{paolo18}). Denote by 
$\O_{\L}^N$ the set of configurations
having $N$ particles in $\L$. Then
\[
\nuln[\{\eta\}] = \frac{\ind(\eta \in \O_{\L}^N)}{Z_{\L}^N} \prod_{x \in \L} \frac{1}{c(\eta_x)!},
\]
where $Z_{\L}^N$ is a normalization factor. It is therefore easily checked that
\be{ob3}
c(\eta_x) \nuln[\{\eta\}] = \frac{Z_{\L}^{N-1}}{Z_{\L}^N} \nu_{\L}^{N-1}[\{\eta - \d_x\}]
\ee
for any $x\in\Lambda$.
Thus, for every function $f$,
\be{ob4}
\nuln[c(\eta_x) f] = \frac{Z_{\L}^{N-1}}{Z_{\L}^N} \nu_{\L}^{N-1}[\s_x f],
\ee
where $\s_x f(\eta) = f(\eta + \d_x)$. Letting $f \equiv 1$ in (\ref{ob4}) we have 
$\frac{Z_{\L}^{N-1}}{Z_{\L}^N} = \nuln[c(\eta_x)]$, and so we rewrite (\ref{ob4}) as
\be{ob5}
\nuln[c(\eta_x) f] =\nuln[c(\eta_x)] \nu_{\L}^{N-1}[\s_x f].
\ee
Now, let $t>0$, and define
\[
\varphi(t) = \nuln\left[e^{tc(\eta_x)}\right].
\]
We now estimate $\varphi(t)$ by means of the so-called Herbst argument (see
\cite{toulouse}, Section 7.4.1). By direct
computation, Jensen's inequality and (\ref{ob5}), we have
\begin{eqnarray}
t \varphi'(t) - \varphi(t) \log \varphi(t) & = & t \nuln\left[ c(\eta_x) e^{tc(\eta_x)} \right] 
- \nuln\left[e^{tc(\eta_x)}\right] \log \nuln\left[e^{tc(\eta_x)}\right] \nonumber \\
 & \leq & t \nuln\left[ c(\eta_x) e^{tc(\eta_x)} \right] - t \nuln[c(\eta_x)] \nuln\left[e^{tc(\eta_x)} \right] \nonumber \\
 & = & t \sum_{x \in \L} \nuln[c(\eta_x)] \left( \nu_{\L}^{N-1} \left[ e^{t  c(\eta_x +1)}\right]
- \nuln\left[e^{tc(\eta_x)}\right] \right) .
\addtocounter{for}{1} \label{ob6}
\end{eqnarray}
We now claim that, for every $x \in \L$ and $0 \leq t \leq 1$
\be{ob7}
\left|\nu_{\L}^{N-1} \left[ e^{t  c(\eta_x +1)}\right]
- \nuln\left[e^{tc(\eta_x)}\right] \right| \leq C_L t \nuln\left[e^{tc(\eta_x)}\right],
\ee
for some constant $C_L$ depending on $L$ but not on $N$. For the moment, let us accept 
(\ref{ob7}), and show how it is used to complete the proof. By (\ref{ob6}), (\ref{ob7}) and the fact that $\nu_\Lambda^N[c(\eta_x)]\leq CN/L$ we
get
\[
t \varphi'(t) - \varphi(t) \log \varphi(t) \leq C_L N t^2 \varphi(t),
\]
for $0 \leq t \leq 1$ and some possibly different $C_L$. Equivalently, letting $\psi(t) = 
\frac{ \log \varphi(t)}{t}$,
\be{ob8}
\psi'(t) \leq C_L N.
\ee
Observing that $\lim _{t \downarrow 0} \psi(t) = \nuln[c(\eta_x)]$, by (\ref{ob8}) and Gronwall   
lemma we have
\[
\psi(t) \leq \nuln[c(\eta_x)] + C_L N t,
\]
from which \eqref{paolo18} easily follows, for $t \geq 0$.

We now prove (\ref{ob7}). We shall use repeatedly the inequality
\be{ob9}
\left| e^x - e^y \right| \leq  |x-y| e^{|x-y|}e^x,
\ee
which holds for $x,y \in \R$. Using (\ref{ob9}) and the fact that, by condition (LG), $\|
c(\eta_x +1) - c(\eta_x) \|_{\infty} \leq a_1$, we have
\[
\left|\nu_{\L}^{N-1} \left[ e^{t  c(\eta_x +1)}\right]
- \nu_{\L}^{N-1} \left[ e^{t  c(\eta_x)}\right] \right| \leq a_1 t e^{a_1 t} \nu_{\L}^{N-1} \left[ e^{t  c(\eta_x)}\right],
\]
from which we have that (\ref{ob7}) follows if we show
\be{ob10}
\left|\nu_{\L}^{N-1} \left[ e^{t c(\eta_x)}\right]
- \nuln\left[e^{tc(\eta_x)}\right] \right| \leq C_L t \nuln\left[e^{tc(\eta_x)}\right].
\ee
At this point we use the notion of stochastic order between probability measures. For two
probabilities $\mu$ and $\nu$ on a partially ordered space $X$, we say that $\nu \prec \mu$
if $\int f d\nu \leq \int f d\mu$ for every integrable, increasing $f$. This is equivalent to the 
existence of a {\em monotone coupling} of $\nu$ and $\mu$, i.e. a probability $P$ on $X \times X$
supported on $\{(x,y) \in X \times X: x \leq y\}$, having marginals $\nu$ and $\mu$ respectively 
(see e.g. \cite{Li}, Chapter 2).

As we will see, (\ref{ob10}) would not be hard to prove if we had $\nu_{\L}^{N-1} \prec \nuln$. 
Under assumptions
(LG) and (M) a slightly weaker fact holds, namely that there is a constant $B>0$ independent of
$N$ and $L$ such that if $N \geq N'+BL$ then $\nu_{\L}^{N'} \prec \nuln$ 
(see \cite{La:Se:Va}, Lemma 4.4). In what follows, we may assume that $N>BL$. Indeed, in the case $N \leq BL$
there is no real dependence on $N$, and (\ref{paolo18}) can be proved by observing that
\[
\nuln\left[e^{t(c(\eta_x) - \nuln[c(\eta_x)])}\right] \leq 1 + \frac{t^2}{2} \nuln[c(\eta_x),c(\eta_x)]e^{2t\sup\{c(\eta_x):\eta 
\in \Omega_\Lambda^N\}} \leq 1+C_L t^2 \leq e^{C_L t^2} \leq e^{C_L N t^2}
\]
as $1 \leq N \leq BL$ ((\ref{paolo18}) is trivial for $N=0$).
For $N \geq BL$ we use the rough inequality
\begin{multline}
\left|\nu_{\L}^{N-1} \left[ e^{t c(\eta_x)}\right]
- \nuln\left[e^{tc(\eta_x)}\right] \right| \leq \left|\nu_{\L}^{N-1} \left[ e^{t c(\eta_x)}\right]
- \nu_{\L}^{N-1-BL}\left[e^{tc(\eta_x)}\right] \right| \\
+  \left|\nu_{\L}^{N} \left[ e^{t c(\eta_x)}\right]
- \nu_{\L}^{N-1-BL}\left[e^{tc(\eta_x)}\right] \right|. \label{ob11}
\end{multline}
Denote by $Q$ the probability measure on $\N^{\L} \times \N^{\L}$ that realizes a monotone coupling
of $\nu_{\L}^{N-1}$ and $\nu_{\L}^{N-1-BL}$. In other words, $Q$ has $\nu_{\L}^{N-1}$ and $\nu_{\L}^{N-1-BL}$
as marginals, and
\be{ob12}
Q[\{(\eta,\xi): \, \eta_x \geq \xi_x \ \forall x \in \L\}] =1.
\ee
Using again (\ref{ob9})
\begin{eqnarray}
\left|\nu_{\L}^{N-1} \left[ e^{t c(\eta_x)}\right]
- \nu_{\L}^{N-1-BL}\left[e^{tc(\eta_x)}\right] \right| & = & \left|Q\left[e^{tc(\eta_x)} - e^{tc(\xi_x)}\right]\right|
\nonumber \\ & 
\leq  & tQ\left[|c(\eta_x) - c(\xi_x)| e^{t|c(\eta_x) - c(\xi_x)|}e^{tc(\eta_x)} \right]. 
\addtocounter{for}{1} \label{ob13}
\end{eqnarray}
Since, $Q$-a.s., $\eta_x \geq \xi_x \ \forall x \in \L$ and $\ov{\eta}_{\L} = \ov{\xi}_{\L} + BL$, then necessarily
$\eta_x \leq \xi_x + BL \ \forall x \in \L$. Thus, it follows that, for some $C>0$,
\[
|c(\eta_x) - c(\xi_x)| \leq CL \ \ \ Q\text{-a.s.}
\]
Inserting this inequality in (\ref{ob13}), we get, for $t \leq 1$ and some $C_L >0$,
\be{ob14}
\left|\nu_{\L}^{N-1} \left[ e^{t c(\eta_x)}\right]
- \nu_{\L}^{N-1-BL}\left[e^{tc(\eta_x)}\right] \right| \leq C_L t \nu_{\L}^{N-1}  \left[ e^{t c(\eta_x)}\right].
\ee
With the same arguments, it is shown that
\be{ob15}
\left|\nu_{\L}^{N} \left[ e^{t c(\eta_x)}\right]
- \nu_{\L}^{N-1-BL}\left[e^{tc(\eta_x)}\right] \right| \leq C_L t \nu_{\L}^{N} \left[ e^{t c(\eta_x)}\right].
\ee
In order to put together (\ref{ob14}) and (\ref{ob15}), we need to show that, for $0 \leq t \leq 1$,
\be{ob16}
\nu_{\L}^{N-1}  \left[ e^{t c(\eta_x)}\right] \leq C_L \nu_{\L}^{N} \left[ e^{t c(\eta_x)}\right],
\ee
for some $L$-dependent $C_L$. By condition (LG) and (\ref{ob5}) we have
\[
\nu_{\L}^{N-1}  \left[ e^{t c(\eta_x)}\right]  \leq  e^{ta_1} \nu_{\L}^{N-1} 
 \left[ e^{t  c(\eta_x +1)}\right] 
 =  \frac{e^{ta_1}}{\nuln[c(\eta_x)]} \nuln\left[c(\eta_x) e^{tc(\eta_x)}\right] 
 \leq  C_L \nu_{\L}^{N} \left[ e^{t c(\eta_x)}\right],
\]
where, for the last step, we have used the facts that, for some $\e>0$, $\nuln[c(\eta_x)] \geq \e \frac{N}{L}$ and
$c(\eta_x) \leq \e^{-1}N$ $\nuln$-a.s. (for both inequalities we use (\ref{lineargrowth})).

The proof (\ref{paolo18}) is now completed for the case $t \geq 0$. 
For $t<0$, it is enough to observe that our 
argument 
is insensitive to replacing $c(\cdot)$ with $-c(\cdot)$.

We now prove (\ref{paolo19}) The idea is to use the fact that the tails of $N(h(\eta_x) - \nuln[h(\eta_x)])$ 
are not thicker than those of $c(\eta_x) - \nuln[c(\eta_x)]$.  Similarly to (\ref{ob3}), note that
\be{obn1}
h(\eta_x) \nuln[\{\eta\}] = \frac{Z_{\L}^{N+1}}{Z_{\L}^N} (\eta_x +1) \nu_{\L}^{N+1} [\{\eta+ \d_x\}],
\ee
which implies
\be{obn2}
\nuln[h(\eta_x) f] = \frac{Z_{\L}^{N+1}}{Z_{\L}^N} \nu_{\L}^{N+1} [\eta_x f(\eta -\d_x)].
\ee
Letting $f \equiv 1$ in (\ref{obn2}) we obtain
\[
\frac{Z_{\L}^{N+1}}{Z_{\L}^N} = \nuln[h(\eta_x)] \frac{L}{N+1}
\]
that, together with the previously obtained identity $\frac{Z_{\L}^{N-1}}{Z_{\L}^N} = \nuln[c(\eta_x)]$, 
yields
\be{obn3}
\nuln[h(\eta_x)] = \frac{N+1}{L \nu_{\L}^{N+1}[c(\eta_x)]}.
\ee
Thus
\begin{multline}
|h(\eta_x) - \nuln[h(\eta_x)]|  =  \left| \frac{\eta_x +1}{c(\eta_x +1)} -  \frac{N+1}{L \nu_{\L}^{N+1}[c(\eta_x)]} \right|
 \\
 \leq  \frac{1}{c(\eta_x +1) \nu_{\L}^{N+1}[c(\eta_x)]} \left[ (\eta_x +1) \left| c(\eta_x +1) - \nu_{\L}^{N+1}
[c(\eta_x)] \right| + c(\eta_x +1) \left| \eta_x + 1 - \frac{N+1}{L}\right| \right]  \\
  \leq  \frac{C}{N} \left[ \left| c(\eta_x +1) - \nu_{\L}^{N+1}
[c(\eta_x)] \right| + \left| \eta_x + 1 - \frac{N+1}{L}\right| \right],
\label{obn4}
\end{multline}
where, in the last step, we use the fact that $\nu_{\L}^{N+1}[c(\eta_x)] \geq \e N/L$ 
for some $\e >0$. In (\ref{obn4}) and in what follows, the $L$-dependence of constants
is omitted. For $\rho>0$, let $c(\rho)$ be obtained by linear interpolation of $c(n), n \in \N$.
Observe that, by (\ref{gap}) with $f(\eta) = \eta_x$,
\[
\nuln[\eta_x,\eta_x] \leq CL^2 \nuln[c(\eta_x)] \leq C_L N
\]
for some $C_L > 0$ that depends on $L$.
Thus, by Condition (LG)
\be{obn5}
\left| \nuln[c(\eta_x)] - c(N/L) \right| \leq c_1 \nuln\left[\left| \eta_x - \frac{N}{L} \right| \right] 
\leq c_1 \sqrt{ \nuln[\eta_x, \eta_x]}  \leq C \sqrt{N},
\ee
for some $C>0$, possibly depending on $L$. From (\ref{obn4}) and (\ref{obn5}) it
follows that, for some $C>0$
\be{obn6}
N|h(\eta_x) - \nuln[h(\eta_x)]| \leq C \left| \eta_x - \frac{N}{L} \right| + C \sqrt{N}.
\ee
Moreover, from Condition (M) it follows that there is a constant $C>0$ such that
\[
 \left| \eta_x - \frac{N}{L} \right| \leq 
C\left(\left| c(\eta_x) - c(N/L) \right|+1\right).
\]
Thus, for every $M>0$ there exists $C>0$ such that
\be{obn7}
\left| c(\eta_x) - \nuln[c(\eta_x)] \right| \leq M \Rightarrow N|h(\eta_x) - \nuln[h(\eta_x)]| \leq CM + C \sqrt{N}.
\ee
It follows that, for all $M>0$
\begin{multline}
\nuln\left[N(h(\eta_x) - \nuln[h(\eta_x)]))> CM + C \sqrt{N}\right]  \leq  \nuln\left[\left| c(\eta_x) - \nuln[c(\eta_x)] \right| > M\right]
 \\
 \leq  \nuln[ c(\eta_x) - \nuln[c(\eta_x)]  > M] +  \nuln[ \nuln[c(\eta_x)] - c(\eta_x)    > M] .
\label{obn8}
\end{multline}
Note that (\ref{obn8}) is trivially true for $M \leq 0$, so it holds for all $M \in \R$.

Now, take $t \in (0,1]$. We have
\begin{eqnarray*}
\nuln \left[ e^{tN(h(\eta_x) - \nuln[h(\eta_x)])} \right] & = & t \int_{-\infty}^{+\infty} e^{tz} \nuln[N(h(\eta_x) - \nuln[h(\eta_x)]) >z]\, dz \\
& \leq &  t \int_{-\infty}^{+\infty} e^{tz}\nuln\left[ c(\eta_x) - \nuln[c(\eta_x)]  > \frac{z}{C} - C \sqrt{N}\right]dz + \\
& & ~~~~~~~~ + 
t \int_{-\infty}^{+\infty} e^{tz}\nuln\left[\nuln[c(\eta_x)] - c(\eta_x)    > \frac{z}{C} - C \sqrt{N}\right]dz \\
& = & Ct \int_{-\infty}^{+\infty} e^{t(CM+C\sqrt{N})} \nuln[ c(\eta_x) - \nuln[c(\eta_x)]  > M]\, dM  \\
& & ~~~~~~~~ + Ct \int_{-\infty}^{+\infty} e^{t(CM+C\sqrt{N})}\nuln[ \nuln[c(\eta_x)] - c(\eta_x)    > M] \,
dM \\
& = & Ce^{Ct\sqrt{N}} \left( \nuln\left[e^{Ct(c(\eta_x) - \nuln[c(\eta_x)])}\right]
 + \nuln\left[e^{-Ct(c(\eta_x) - \nuln[c(\eta_x)])} \right] \right) \\
&  \leq  & 2C e^{Ct\sqrt{N}} e^{ANt^2},
\end{eqnarray*}
where, in the last step, we have used (\ref{paolo18}). This completes the proof for the case $t \in (0,1]$.
For $t \in [-1,0)$ we proceed similarly, after having replaced, in (\ref{obn8}), $N(h(\eta_x) - \nuln[h(\eta_x)])$
with its opposite.

\endproofof

\vspace{1cm}

\noindent
\proofof{Proposition \ref{ppaolo5}}
We first prove inequality (\ref{paolo20}). Clearly, it is enough to show that, for all $x \in \L$,
\[
\nuln\left[f,  c(\eta_x)\right]^2 \leq C_L N \nuln[f] \ent_{\nuln}(f),
\] 
for some $C_L >0$

By the entropy inequality (\ref{paolo17}) and (\ref{paolo18}) we have,
for $0<t\leq 1$:
\be{ob17}
\nuln\left[f, c(\eta_x)\right]^2 \leq 2\nuln[f]^2 A^2_L  N^2  t^2 + \frac{2}{t^2} \ent^2_{\nuln}(f).
\ee
Set 
\be{ob18}
t_*^2 = \frac{\ent_{\nuln}(f)}{N\nuln[f]}.
\ee
Suppose, first, that $t_* \leq 1$. Then if we insert $t_*$ in (\ref{ob17}) we get
\be{ob19}
\nuln\left[f,  c(\eta_x)\right]^2 \leq 2\nuln[f] A^2_L N \ent_{\nuln}(f) + 2N \nuln[f]\ent_{\nuln}(f) =:
C_L N \nuln[f]\ent_{\nuln}(f) .
\ee
In the case $t_* >1$, we easily have
\be{ob20}
\nuln\left[f,  c(\eta_x)\right]^2 \leq C \nuln[f]^2 N^2 \leq C \nuln[f]^2 N^2 {t^*}^2 = CN \nuln[f]\ent_{\nuln}(f) .
\ee
By (\ref{ob19}) and (\ref{ob20}) the conclusion follows.

Let us now prove (\ref{paolo21}). As before, by (\ref{paolo17}) and (\ref{paolo19}), for $t \in (0,1]$:
\be{ob21}
N^2 \nuln\left[f, h(\eta_x) \right]^2 \leq \frac{2\nuln[f]^2}{t^2}\left[ \log^2 A_L + A^2_L t^2 N + 
 A^2_L  N^2  t^4 \right] + \frac{2}{t^2} \ent^2_{\nuln}(f).
\ee
Here we set
\[
t_*^2 = \frac{\ent_{\nuln}(f) \vee \nuln[f]}{N\nuln[f]},
\]
and proceed as in (\ref{ob19}) and (\ref{ob20}).
\endproofof

\section{Local limit theorems}

The rest of the paper is devoted to the proofs of all technical results that have been used in previous sections.
For the sake of generality, we state and prove all results in dimension $d\geq 1$.

Using the language of statistical mechanics, having defined the canonical measure $\nuln$,
for $\rho >0$
we consider the corresponding {\em grand canonical} measure 
\be{grandcan}
\mur[\{\eta\}] :=  \frac{\a(\rho)^{\ov{\eta}_{\L}} \nu_{\L}^{\ov{\eta}_{\L}}[\{\eta\}]}{\sum_{\xi \in \O_{\L}}
\a(\rho)^{\ov{\xi}_{\L}} \nu_{\L}^{\ov{\xi}_{\L}}[\{\xi\}]} = \frac{1}{Z(\a(\rho))} \prod_{x \in \L} \frac{\a(\rho)^{\eta_x}}
{c(\eta_x)!},
\ee
where $\a(\rho)$ is chosen so that $\mur(\eta_x) = \rho$, $x \in \L$, and
$Z(\a(\rho))$ is the corresponding normalization. Clearly
$\mur$ is a product measure with marginals given by (\ref{int2}). Monotonicity and the Inverse
Function Theorem for analytic functions, guarantee that $\a(\rho)$ is well defined and it is
a analytic function of $\rho \in [0,+\infty)$. We state here without proof some direct consequences of Conditions \ref{co:LG} and \ref{co:M}. The proofs
of some of them can be found in \cite{La:Se:Va}.

\begin{pro}
  \label{lemma:M1}
  \ 
  
  \begin{enumerate}
    \item Let $\sigma^2(\rho):=\mu_\rho[\eta_x,\eta_x]$, then
  \label{pro:mu}
    \begin{equation}
      \label{eq:sigma}
    0<
    \inf_{\rho>0}\frac{\sigma^2(\rho)}{\rho}\leq
    \sup_{\rho>0}\frac{\sigma^2(\rho)}{\rho}<
    +\infty
    \end{equation}
   \item Let $\alpha(\rho)$ be the function appearing in (\ref{grandcan}); then
    \begin{equation}
      \label{eq:alpha}
    0<
    \inf_{\rho>0}\frac{\alpha(\rho)}{\rho}\leq
    \sup_{\rho>0}\frac{\alpha(\rho)}{\rho}<
    +\infty
    \end{equation}
  \end{enumerate}
\end{pro}

\subsection{Local limit theorems for the grand canonical measure}

The next two results are a form of local limit theorems for the density of $\n\eta\Lambda$ under $\mu_\rho$. Define  $p_\Lambda^\rho(n):=\mu_\rho[\n\eta\Lambda=n]$ for $\rho>0$, $n\in\natural$ and $\Lambda\subset\subset\integer^d$.
The idea is to get a Poisson approximation of $p_\Lambda^{N/|\Lambda|}(n)$ for very small values of $N/|\Lambda|$ and to use the uniform local limit theorem (see Theorem~6.1 in \cite{La:Se:Va}) for the other cases.  

\begin{lem}
  \label{lemma:poi}
  For every $N_0\in\natural\setminus\{0\}$ there exist $A_0>0$ such that 
  \begin{displaymath}
    \sup_{\substack{N\in\natural\setminus\{0\},\; n\in\natural\\  n\leq N\leq N_0}}\left|p_\Lambda^{N/|\Lambda|}(n)-\frac{N^n}{n!}e^{-N}\right|\leq
    \frac{A_0}{|\Lambda|}
  \end{displaymath}
  for any $\Lambda\subset\subset \integer^d$.
\end{lem}

\begin{proof}
  Let $\rho:=N/|\Lambda|$ and assume, without loss of generality, that $0\in\Lambda$.
  Notice that
  \begin{multline*}
    \mu_\rho\left[\n\eta\Lambda=n\right]\\
    =\mu_\rho\Big[\n\eta\Lambda=n\Big|\max_{x\in\Lambda}\eta_x\leq 1\Big]\mu_\rho\Big[\max_{x\in\Lambda}\eta_x\leq 1\Big]+
    \mu_\rho\Big[\n\eta\Lambda=n\Big|\max_{x\in\Lambda}\eta_x> 1\Big]\mu_\rho\Big[\max_{x\in\Lambda}\eta_x> 1\Big].
  \end{multline*}

  We begin by proving that
  \begin{equation}
    \label{eq:exstp1}
    \mu_\rho\Big[\max_{x\in\Lambda}\eta_x> 1\Big]=O(|\Lambda|^{-1}),
  \end{equation}
  uniformly in $0<N\leq N_0$.

  Indeed
  \begin{displaymath}
    \mu_\rho\Big[\max_{x\in\Lambda}\eta_x> 1\Big]=
     1-\left(1-\mu_\rho\left[\eta_0> 1\right]\right)^{|\Lambda|},
  \end{displaymath}
  and
  \begin{displaymath}
    \mu_\rho\left[\eta_0> 1\right]=
    \frac{1}{Z(\rho)}\sum_{k=2}^{+\infty}\frac{\alpha(\rho)^k}{c(k)!}=
    \frac{\alpha(\rho)^2}{Z(\rho)}\sum_{k=0}^{+\infty}\frac{\alpha(\rho)^k}{c(k+2)!}=
    \frac{\alpha(\rho)^2}{Z(\rho)}\sum_{k=0}^{+\infty}\frac{c(k)!}{{c(k+2)!}}\frac{\alpha(\rho)^k}{c(k)!}.
  \end{displaymath}
  Since $c(k)!/c(k+2)!$ is uniformly bounded, we have
  $\mu_\rho\left[\eta_0> 1\right]\leq B_1\alpha(\rho)^2=O(|\Lambda|^{-2})$, uniformly in $0<N\leq N_0$.
  Thus
  \begin{displaymath}
    \left(1-\mu_\rho\left[\eta_0> 1\right]\right)^{|\Lambda|}=
    \left(1-O(|\Lambda|^{-2})\right)^{|\Lambda|}=O(|\Lambda|^{-1}),
  \end{displaymath}
  which establishes \eqref{eq:exstp1}.

  Now let
  \begin{displaymath}
    \tilde{\rho}:=\sum_{k=0}^{+\infty} k \mu_\rho\Big[\eta_0=k\Big|\eta_0\leq 1\Big].
  \end{displaymath}
  
  A trivial calculation shows that
  \begin{displaymath}
    \tilde{\rho}=
    \frac{\mu_\rho[\eta_0\ind(\eta_0\leq 1)]}{\mu_\rho[\eta_0\leq 1]}=
    \frac{\mu_\rho[\eta_0]-\mu_\rho[\eta_0\ind(\eta_0>1)]}{\mu_\rho[\eta_0\leq 1]}\\
    =\frac{\rho-\mu_\rho[\eta_0\ind(\eta_0>1)]}{\mu_\rho[\eta_0\leq 1]}.
  \end{displaymath}
  Moreover
  \begin{multline*}
    \mu_\rho[\eta_0\ind(\eta_0>1)]=
    \frac{1}{Z(\rho)}\sum_{k=2}^{+\infty}\frac{k\alpha(\rho)^k}{c(k)!}=
    \frac{\alpha(\rho)^2}{Z(\rho)}\sum_{k=0}^{+\infty}\frac{(k+2)\alpha(\rho)^k}{c(k+2)!}\\
    \leq\frac{B_2\alpha(\rho)^2}{Z(\rho)}\sum_{k=0}^{+\infty}\frac{(k+2)\alpha(\rho)^k}{c(k)!}=
    B_2\alpha(\rho)^2(\rho+2)=
    O(|\Lambda|^{-2}),
  \end{multline*}
  and finally
  \begin{equation}
    \label{eq:exstp2}
    \tilde\rho
    =\frac{\rho+O(|\Lambda|^{-2})}{1+O(|\Lambda|^{-2})}=\rho+O(|\Lambda|^{-2}).
  \end{equation}

  Observe that for any $n\in\{0,\dots,|\Lambda|\}$, we have
  \begin{equation}
    \label{eq:poi1}
    \mu_\rho\Big[\n\eta\Lambda=n\Big|\max_{x\in\Lambda}\eta_x\leq 1\Big]=
    \binom{|\Lambda|}{n}\tilde\rho(1-\tilde\rho)^{|\Lambda|-n}.
  \end{equation}

  This comes from the fact that the random variables $\left\{\eta_x: x\in\Lambda\right\}$, under the probability measure $\mu_\rho[\cdot|\max_{x\in\Lambda}\eta_x\leq 1]$, are Bernoulli independent random variables with mean $\tilde\rho$.
  The remaining part of the proof follows the classical argument of approximation of the binomial distribution with the Poisson distribution.
  Using \eqref{eq:exstp2} and \eqref{eq:poi1}, after some simple calculations we get 
  \begin{multline*}
    \mu_\rho\Big[\n\eta\Lambda=n\Big|\max_{x\in\Lambda}\eta_x\leq 1\Big]=
    \binom{|\Lambda|}{n}\tilde\rho(1-\tilde\rho)^{|\Lambda|-n}\\
    =\frac{|\Lambda|!}{n!(|\Lambda|-n)!}\left[\frac{N}{|\Lambda|}+O(|\Lambda|^{-2})\right]^n\left[1-\frac{N}{|\Lambda|}+O(|\Lambda|^{-2})\right]^{|\Lambda|-n}\\
    =\frac{1}{n!}\left[N+O(|\Lambda|^{-1})\right]^n\left[1-\frac{N}{|\Lambda|}+O(|\Lambda|^{-2})\right]^{|\Lambda|}\\
    \times\frac{|\Lambda|(|\Lambda|-1)\cdot\dots\cdot(|\Lambda|-n+1)}{|\Lambda|^n}\left[1-\frac{N}{|\Lambda|}+O(|\Lambda|^{-2})\right]^{-n}=
    \frac{N^n}{n!}e^{-N}+O(|\Lambda|^{-1})
  \end{multline*}
  uniformly in $0<N\leq N_0$ and $0\leq n <N$.
  This proves that there exist positive constants $v_1$ and $B_3$ such that if $\Lambda\subset\subset\integer^d$ is such that $|\Lambda|>v_1$ then
    \begin{displaymath}
    \left|p_\Lambda^{N/|\Lambda|}(n)-\frac{N^n}{n!}e^{-N}\right|\leq
    \frac{B_3}{|\Lambda|}
  \end{displaymath}
    uniformly in $N\in\natural\setminus\{0\}$ with $N\leq N_0$ and $n\in\natural$ with $n\leq N$.
    The general case follows easily because the set of $n\in\natural$, $N\in\natural\setminus\{0\}$ and $\Lambda\subset\subset\integer^d$ such that $n\leq N\leq N_0$, $|\Lambda|\leq v_1$ and $0\in\Lambda$ is finite.

\end{proof}

\begin{pro}
  \label{pro:LLT}
   For any $\rho_0>0$, there exist finite positive constants $A_0$, $n_0$ and $v_0$ such that:
  \begin{enumerate}
  \item 
  \begin{equation}
    \label{eq:LLT:1}
   \sup_{n\in\natural}\left|\sqrt{\sigma^2(\rho)|\Lambda|}p_{\Lambda}^\rho(n)-\frac{1}{\sqrt{2 \pi}}e^{-\frac{(n-\rho|\Lambda|)^2}{2\sigma^2(\rho)|\Lambda|}}\right|\leq
      \frac{A_0}{\sqrt{\sigma^2(\rho)|\Lambda|}}
  \end{equation}
  for any $\rho\leq\rho_0$ and any $\Lambda\subset\subset\integer^d$ such that $\sigma^2(\rho)|\Lambda|\geq n_0$;
  \item
  \begin{equation}
    \label{eq:LLT:2}
   \sup_{\substack{\rho>\rho_0 \\ n\in\natural}}\left|\sqrt{\sigma^2(\rho)|\Lambda|}p_{\Lambda}^\rho(n)-\frac{1}{\sqrt{2 \pi}}e^{-\frac{(n-\rho|\Lambda|)^2}{2\sigma^2(\rho)|\Lambda|}}\right|\leq
      \frac{A_0}{\sqrt{|\Lambda|}}
  \end{equation}
  for any $\Lambda\subset\subset\integer^d$ such that $|\Lambda|\geq v_0$.    
  \end{enumerate}
\end{pro}

\begin{proof}
  This is a special case of the Local Limit Theorem for $\mu_\rho$ (see Theorem~6.1 in \cite{La:Se:Va}).
\end{proof}

We conclude this section with a bound on the tail of $p_\Lambda^\rho$ which will be used in the regimes not covered by Lemma~\ref{lemma:poi} or Proposition~\ref{pro:LLT}.

\begin{lem}
  \label{lemma:dec}
  There exists a positive constant $A_0$ such that
  \begin{displaymath}
    \frac{\rho|\Lambda|}{A_0(n+1)}\leq
    \frac{p_\Lambda^\rho(n+1)}{p_\Lambda^\rho(n)}\leq
    \frac{A_0\rho|\Lambda|}{n+1}
  \end{displaymath}
  for any $n\in\natural\setminus\{0\}$, $\Lambda\subset\subset\integer^d$ and $\rho>0$.
\end{lem}

\begin{proof}
  Notice that
  \begin{displaymath}
    p_\Lambda^\rho(n+1)=
    \mu_\rho[\n\eta\Lambda=n+1]=
    \frac{1}{n+1}\sum_{x\in\Lambda}\mu_\rho[\eta_x,\ind(\n\eta\Lambda=n+1)]
  \end{displaymath}
  and, by the change of variable $\xi:=\sigma-\delta^x$:
  \begin{multline*}
    \mu_\rho[\eta_x\ind(\n\eta\Lambda=n+1)]=
    \sum_{\sigma\in\Omega_\Lambda}\mu_\rho[\eta_\Lambda=\sigma]\sigma_x\ind(\n\sigma\Lambda=n+1)=
    \sum_{\sigma\in\Omega_\Lambda}\mu_\rho[\eta_\Lambda=\sigma]\sigma_x\ind(\n\sigma\Lambda=n+1,\sigma_x>0)\\
    =\sum_{\xi\in\Omega_\Lambda}\mu_\rho[\eta_\Lambda=\xi+\delta^x](\xi_x+1)\ind(\n\xi\Lambda=n)=
    \sum_{\xi\in\Omega_\Lambda}\mu_\rho[\eta_\Lambda=\xi]\frac{\mu_\rho[\eta_\Lambda=\xi+\delta^x]}{\mu_\rho[\eta_\Lambda=\xi]}(\xi_x+1)\ind(\n\xi\Lambda=n)\\
    =\sum_{\xi\in\Omega_\Lambda}\mu_\rho[\eta_\Lambda=\xi]\frac{\alpha(\rho)(\xi_x+1)}{c(\xi_x+1)}\ind(\n\xi\Lambda=n)=
    \alpha(\rho)\mu_\rho\left[\frac{\eta_x+1}{c(\eta_x+1)},\n\eta\Lambda=n \right].
  \end{multline*}
  This means  that
  \begin{equation}
    \label{eq:dec:1}
    p_\Lambda^\rho(n+1)=
    \frac{\alpha(\rho)}{n+1}\sum_{x\in\Lambda}\mu_\rho\left[\frac{\eta_x+1}{c(\eta_x+1)}\ind(\n\eta\Lambda=n) \right].
  \end{equation}
  By (\ref{lineargrowth}) we know that there exists a positive constant $B_0$ such that
  \begin{displaymath}
    B_0^{-1}\leq
    \frac{c(k)}{k}\leq
    B_0
    \qquad\text{for any $k\in\natural\setminus\{0\}$.}
  \end{displaymath}
  Thus, by plugging these bounds in \eqref{eq:dec:1}, we get
  \begin{displaymath}
    \frac{\alpha(\rho)|\Lambda|p_\Lambda^\rho(n)}{B_0(n+1)}\leq
    p_\Lambda^\rho(n+1)\leq
    \frac{B_0\alpha(\rho)|\Lambda|p_\Lambda^\rho(n)}{n+1},
  \end{displaymath}
  from which the conclusion follows.
\end{proof}

\subsection{Gaussian estimates for the canonical measure}
\label{sec:gau}

In this section we will prove some Gaussian bounds on $\nu_\Lambda^N[\n\eta{\Lambda^\prime}=\cdot]$, when the volumes $|\Lambda|$ and $|\Lambda^\prime|$ are of the same order (typically it will be $|\Lambda^\prime|/|\Lambda|=1/2$).
These bounds are volume dependent (see Proposition~\ref{pro:gaul} below) and so are of limited utility.
However they will be used to prove Proposition~\ref{pro:cen} which, in turn, is used in Section~\ref{sec:odlsi} to regularize $\nu_\Lambda^N[\n\eta{\Lambda^\prime}=\cdot]$ and prove Gaussian uniform estimates on it. 

Assume $\Lambda^\prime\subset\Lambda\subset\subset\integer^d$, $N\in\natural\setminus\{0\}$ and consider the probability measure on $\{0,1,\dots,N\}$

\begin{displaymath}
  \nu_\Lambda^N[\n\eta{\Lambda^\prime}=n]=
  \frac{p_{\Lambda^\prime}^\rho(n)p_{\Lambda\setminus\Lambda^\prime}^\rho(N-n)}{p_\Lambda^\rho(N)}.
\end{displaymath}
Notice that $\nu_\Lambda^N[\n\eta{\Lambda^\prime}=\cdot]$ does not depend on $\rho>0$ (and depends on $\Lambda$ and $\Lambda^\prime$ only through $|\Lambda|$ and $|\Lambda^\prime|$).
Indeed a simple computation shows that
\begin{displaymath}
  \nu_\Lambda^N[\n\eta\Lambda=n]=
  \frac{\sum_{\sigma\in\Omega_{\Lambda^\prime}}\frac{\ind(\n\sigma{\Lambda^\prime}=n)}{\prod_{x\in\Lambda^\prime}c(\sigma_x)!}\sum_{\xi\in\Omega_{\Lambda\setminus\Lambda^\prime}}\frac{\ind(\n\xi{\Lambda\setminus\Lambda^\prime}=N-n)}{\prod_{y\in\Lambda\setminus\Lambda^\prime}c(\xi_y)!}}{\sum_{n=0}^N \sum_{\sigma\in\Omega_{\Lambda^\prime}}\frac{\ind(\n\sigma{\Lambda^\prime}=n)}{\prod_{x\in\Lambda^\prime}c(\sigma_x)!}\sum_{\xi\in\Omega_{\Lambda\setminus\Lambda^\prime}}\frac{\ind(\n\xi{\Lambda\setminus\Lambda^\prime}=N-n)}{\prod_{y\in\Lambda\setminus\Lambda^\prime}c(\xi_y)!}},
\end{displaymath}
for any $n\in\{0,1,\dots,N\}$.
A particular case is when $\Lambda^\prime\subset\Lambda\subset\subset\integer^d$ is such that $|\Lambda^\prime|=|\Lambda\setminus\Lambda^\prime|$.
In this case we define

\begin{equation}
  \label{eq:gamma}
  \gamma_\Lambda^N(n):=
  \nu_\Lambda^N[\n\eta{\Lambda^\prime}=n]=
  \frac{p_{\Lambda^\prime}^\rho(n)p_{\Lambda\setminus\Lambda^\prime}^\rho(N-n)}{p_\Lambda^\rho(N)}=
  \frac{p_{\Lambda^\prime}^\rho(n)p_{\Lambda^\prime}^\rho(N-n)}{p_\Lambda^\rho(N)}.
\end{equation}
We begin by proving a simple result on the decay of the tails of $\gamma_\Lambda^N$.

\begin{pro}
  \label{pro:dec}
  There exist a positive constant $A_0$ such that for every $\Lambda\subset\subset\integer^d$ such that $|\Lambda|\geq 2$ and every $N\in\natural\setminus\{0\}$
  \begin{equation}
    \label{eq:dec:0}
    \frac{N-n}{A_0(n+1)}\leq
    \frac{\gamma_\Lambda^N(n+1)}{\gamma_\Lambda^N(n)}\leq
    \frac{A_0(N-n)}{(n+1)},
  \end{equation}
  for any $n\in\{0,\dots,N-1\}$.
\end{pro}

\begin{proof}
  Let $\Lambda^\prime\subset\Lambda$ be a non empty lattice set.
  Then by \eqref{eq:gamma}
  \begin{displaymath}
    \frac{\gamma_\Lambda^N(n+1)}{\gamma_\Lambda^N(n)}=
    \frac{p_{\Lambda^\prime}^\rho(n+1)p_{\Lambda\setminus\Lambda^\prime}^\rho(N-n-1)}{p_{\Lambda^\prime}^\rho(n)p_{\Lambda\setminus\Lambda^\prime}^\rho(N-n)},
  \end{displaymath}
  and \eqref{eq:dec:0} follows from Lemma~\ref{lemma:dec}.
\end{proof}

Next we prove Gaussian bounds on $\gamma_\Lambda^N$.
This is a very technical argument.
We begin with the case $|\Lambda|=2$.

\begin{lem}
  \label{lemma:gau2}
  Assume that $|\Lambda|=2$, and define $\bar N:=\lceil N/2 \rceil$ for $N\in\natural\setminus\{0\}$.
  There exist a positive constant $A_0$ such that
  \begin{equation}
    \label{eq:gau2}
    \frac{1}{A_0\sqrt{\bar N}}e^{-\frac{A_0(n-\bar N)^2}{\bar N}}\leq
    \gamma_\Lambda^N(n)\leq
    \frac{A_0}{\sqrt{\bar N}}e^{-\frac{(n-\bar N)^2}{A_0\bar N}},
  \end{equation}
  uniformly in $N\in\natural\setminus\{0\}$ and $n\in\{0,1,\dots,N\}$.
\end{lem}

\begin{proof}
  We split the proof in several steps for clarity purpose.
  
\paragraph{Step 1.}

 There exists $A_1>0$ such that for any $N\in\natural\setminus\{0\}$,
  \begin{equation}
   \label{eq:gau2:stp1:1}
   \log\gamma_\Lambda^N(n-1)-\log\gamma_\Lambda^N(n)\leq
   -\frac{\bar N-n}{A_1\bar N}+\frac{A_1}{N}
 \end{equation}
  for any $n\in\{1,\dots,\bar N-1\}$
 and
 \begin{equation}
   \label{eq:gau2:stp1}
   \log\gamma_\Lambda^N(n+1)-\log\gamma_\Lambda^N(n)\leq
   -\frac{n-\bar N}{A_1\bar N}+\frac{A_1}{N}
 \end{equation}
  for any $n\in\{\bar N,\dots,N-1\}$.

  \subparagraph{Proof of Step 1.}

  Assume $n\in\{\bar N,\dots,N-1\}$, the other case will follow by a symmetry argument since $\gamma_\Lambda^N(n)=\gamma_\Lambda^N(N-n)$.
  By \eqref{eq:gamma} we obtain, in the case $|\Lambda|=2$, that
  \begin{displaymath}
    \frac{\gamma_\Lambda^N(n+1)}{\gamma_\Lambda^N(n)}=
    \frac{c(N-n)}{c(n+1)}.
  \end{displaymath}
  If $n+1-(N-n)\leq k_0$, where $k_0$ is the constant which appears in Condition~\ref{co:M},
then
  by Condition~\ref{co:LG} and (\ref{lineargrowth}) there exist a constant $B_1>0$ such that
  \begin{displaymath}
    \frac{c(N-n)}{c(n+1)}=
    1-\frac{c(n+1)-c(N-n)}{c(n+1)}\leq
    1+\frac{|c(n+1)-c(N-n)|}{c(n+1)|}\leq
    1+\frac{a_1k_0B_1}{n+1}
  \end{displaymath}
  and, since $\log (1+x)\leq x$ and $n\geq\bar N$, there exist a constant $B_2>0$ such that
  \begin{displaymath}
    \log\gamma_\Lambda^N(n+1)-\log\gamma_\Lambda^N(n)\leq
    \frac{a_1k_0B_1}{n+1}\leq
    \frac{B_2}{N}.
  \end{displaymath}
  Since $n+1-(N-n)\leq k_0$ we have that $n-\bar N\leq (k_0-1)/2$, which implies $(n-\bar N)/\bar N\leq (k_0-1)/N$.
  Thus
  \begin{displaymath}
    \frac{B_2}{N}\leq
    \frac{B_2}{N}+\frac{k_0-1}{N}-\frac{n-\bar N}{\bar N}\leq
    \frac{B_3}{N}-\frac{n-\bar N}{\bar N}
  \end{displaymath}
  and so
  \begin{displaymath}
    \log\gamma_\Lambda^N(n+1)-\log\gamma_\Lambda^N(n)\leq
    \frac{B_3}{N}-\frac{n-\bar N}{\bar N}
  \end{displaymath}
  in this case.
  Assume now that $n+1-(N-n)> k_0$;
  more precisely assume that $rk_0<n+1-(N-n)\leq (r+1)k_0$ for some $r\in\natural\setminus\{0\}$.
  Then
  \begin{multline}
    \label{eq:gau2:1}
    \log\gamma_\Lambda^N(n+1)-\log\gamma_\Lambda^N(n)=
    \log c(N-n)-\log c(n+1)\\
    =\sum_{s=1}^r\left[\log c(N-n+k_0(s-1))-\log c(N-n+k_0s)\right]+\log c(n-N+k_0r)-\log c(n+1).
  \end{multline}
  For any $s\in\{1,2,\dots,r\}$, by Condition~\ref{co:M} and the fact that $\log(1-x)\leq -x$, we obtain
  \begin{multline*}
    \log c(N-n+k_0(s-1))-\log c(N-n+k_0s)=
    \log\left[1-\frac{c(N-n+k_0s)-c(N-n+k_0(s-1))}{c(N-n+k_0s)}\right]\\
    \leq-\frac{c(N-n+k_0s)-c(N-n+k_0(s-1))}{c(N-n+k_0s)}\leq
    -\frac{a_2}{c(N-n+k_0s)}.
  \end{multline*}
   By (\ref{lineargrowth}) we know that $c(N-n+k_0s)\leq a_1 (N-n+k_0s)$ and because $N-n+k_0s\leq N-n+k_0r\leq n+1$ we get
  \begin{equation}
    \label{eq:gau2:2}
    \log c(N-n+k_0(s-1))-\log c(N-n+k_0s)\leq 
    -\frac{a_2}{a_1(n+1)}.
  \end{equation}
  Furthermore by Condition~\ref{co:LG} and the fact that $\log(1-x)\leq -x$ we obtain
  \begin{multline*}
    \log c(N-n+k_0r)-\log c(n+1)=
    \log\left[1-\frac{c(n+1)-c(N-n+k_0r)}{c(n+1)}\right]\\
    \leq-\frac{c(n+1)-c(N-n+k_0r)}{c(n+1)}\leq
    \frac{|c(n+1)-c(N-n+k_0r)|}{c(n+1)}\leq
    \frac{a_1}{c(n+1)},
  \end{multline*}
  and by (\ref{lineargrowth}) we have that there exists $B_4>0$ such that $c(n+1)\geq B_4^{-1}(n+1)$.
  Thus, since $n+1>\bar N$, we have
  \begin{equation}
    \label{eq:gau2:3}
    \log c(N-n+k_0r)-\log c(n+1)\leq
    \frac{2a_1B_4}{N}.
  \end{equation}
  By plugging the bounds \eqref{eq:gau2:3} and \eqref{eq:gau2:2} into \eqref{eq:gau2:1} we obtain
  \begin{equation}
    \label{eq:gau2:4}
    \log\gamma_\Lambda^N(n+1)-\log\gamma_\Lambda^N(n)\leq
    -\frac{ra_2}{a_1(n+1)}+\frac{2a_1B_4}{N}.
  \end{equation}
  Recalling that $(r+1)k_0\geq n+1-(N-n)=2n-N+1$ and $k_0\geq1$, we have $rk_0\geq n+1-(N-n)=2n-N$ and
  \begin{displaymath}
    \frac{r}{n+1}\geq
    \frac{2n-N}{k_0(n+1)}\geq
    \frac{2n-N}{k_0N}\geq
    \frac{n-\bar N}{k_0\bar N}
  \end{displaymath}
  which together with \eqref{eq:gau2:4} completes the proof of \eqref{eq:gau2:stp1}.
  
  \paragraph{Step 2.}

There exists $A_2>0$ such that if $N\in\natural\setminus\{0\}$ then,
    \begin{equation}
   \label{eq:gau2:stp2:1}
   \log\gamma_\Lambda^N(n-1)-\log\gamma_\Lambda^N(n)\geq
   -A_2\frac{\bar N-n}{\bar N}-\frac{A_2}{N}
 \end{equation}
  for any $n\in\{\lceil  N/4 \rceil,\dots,\bar N-1\}$
 and
 \begin{equation}
   \label{eq:gau2:stp2:2}
   \log\gamma_\Lambda^N(n+1)-\log\gamma_\Lambda^N(n)\geq
   -A_2\frac{n-\bar N}{\bar N}-\frac{A_2}{N}
 \end{equation}
  for any $n\in\{\bar N,\dots,\lfloor 3N/4\rfloor\}$.

  \subparagraph{Proof of Step 2.}

    Assume $n\in\{\bar N,\dots,\lfloor 3N/4\rfloor\}$, the other case again  will follow by symmetry as in \textbf{Step 1}.
    Notice that
  \begin{multline}
    \label{eq:gau2:stp2:3}
    \log\gamma_\Lambda^N(n+1)-\log\gamma_\Lambda^N(n)
    =\log c(N-n) -\log c(n+1)\\
    =\sum_{s=1}^{2n-N+1}[\log c(N-n + s-1)-\log c(N-n+s)].
  \end{multline}
  Moreover
  \begin{multline*}
    \log c(N-n + s-1)-\log c(N-n+s)=
    \log \left[1-\frac{c(N-n+s)-c(N-n+s-1)}{c(N-n+s)}\right]\\
    \geq \log\left[1-\frac{|c(N-n+s)-c(N-n+s-1)|}{c(N-n+s)}\right];
  \end{multline*}
  by Condition~\ref{co:LG} and (\ref{lineargrowth}) there exist $B_1>0$ such that
  \begin{displaymath}
    \frac{|c(N-n+s)-c(N-n+s-1)|}{c(N-n+s)}\leq
    \frac{B_1}{N-n+s}\leq
    \frac{B_1}{N-n}\leq
    \frac{4B_1}{ N},
  \end{displaymath}
  where we used the fact that $n\leq 3N/4$;
  thus
  \begin{multline*}
    \log c(N-n + s-1)-\log c(N-n+s)\\
    \geq\log\left[1-\frac{|c(N-n+s)-c(N-n+s-1)|}{c(N-n+s)}\right]\geq
    \log\left(1-\frac{4B_1}{ N}\right).
  \end{multline*}
  Since $\log(1-x)\geq -2x$, for $x\in[0,1/2]$, then 
  \begin{displaymath}
    \log\left(1-\frac{4B_1}{ N}\right)\geq
    -{8B_1}{ N}
  \end{displaymath}
  for $N\geq n_1:=8B_1$.
  Thus
  \begin{displaymath}
    \log c(N-n + s-1)-\log c(N-n+s)\geq
    -\frac{8B_1}{N}.
  \end{displaymath}
  By plugging this bound into \eqref{eq:gau2:stp2:3} we obtain
  \begin{displaymath}
    \log\gamma_\Lambda^N(n+1)-\log\gamma_\Lambda^N(n)\geq
    -8B_1\frac{2n-N+1}{ N}\geq
    -16B_1\frac{n-\bar N}{\bar N}-\frac{8B_1}{\bar N},
  \end{displaymath}
  which completes the proof of \eqref{eq:gau2:stp2:2} in the case $N\geq n_1$.
  The general case is obtained by a a finiteness argument.
  
  \paragraph{Step 3.}

  There exist $A_3>0$ such that for any $N\in\natural\setminus\{0\}$,
  \begin{equation}
    \label{eq:gau2:stp3}
    \frac{\gamma_\Lambda^N(n)}{\gamma_\Lambda^N(\bar N)}\leq
    A_3e^{-\frac{(n-\bar N)^2}{A_3\bar N}},
  \end{equation}
  for any $n\in\{0,\dots,N\}$.

  \subparagraph{Proof of Step 3.}

  Assume that $n\in\{\bar N+1,\dots,N\}$.
  By \eqref{eq:gau2:stp1} there exists $B_1>0$ such that 
  \begin{equation}
    \label{eq:gau2:stp3:1}
    \log \gamma_\Lambda^N(n)-\log \gamma_\Lambda^N(\bar N)=
    \sum_{k=\bar N}^{n-1} \left[\log \gamma_\Lambda^N(k+1)-\log\gamma_\Lambda^N(k)\right]\leq
    \sum_{k=\bar N}^{n-1}\left(\frac{B_1}{N}-\frac{k-\bar N}{B_1\bar N}\right).
  \end{equation}
  Using the fact that $n-\bar N\leq N/2$, and some elementary computation we obtain
  \begin{multline*}
    \sum_{k=\bar N}^{n-1}\left(\frac{B_1}{N}-\frac{k-\bar N}{B_1\bar N}\right)\leq
    \frac{B_1(n-\bar N)}{N}-\frac{1}{B_1\bar N}\sum_{k=1}^{n-\bar N-1}k\leq
    \frac{B_1}{2}-\frac{(n-\bar N)(n-\bar N-1)}{B_1 N}\\
    =\frac{B_1}{2}-\frac{(n-\bar N)^2}{B_1 N}+\frac{n-\bar N}{B_1 N}\leq
    B_1-\frac{(n-\bar N)^2}{B_1 N},
  \end{multline*}
  which plugged into \eqref{eq:gau2:stp3:1} implies \eqref{eq:gau2:stp3} for $n\in\{\bar N+1,\dots,N\}$.
  The case $n\in\{0,\dots,\bar N-1\}$ is again obtained by symmetry, while the case $n=\bar N$ is trivial.

  \paragraph{Step 4.}

    There exist $N_2>0$ and $A_4>0$ such that if $N\in\natural\setminus\{0\}$ then
  \begin{equation}
    \label{eq:gau2:stp4}
    \frac{\gamma_\Lambda^N(n)}{\gamma_\Lambda^N(\bar N)}\geq
    \frac{1}{A_4}e^{-A_4\frac{(n-\bar N)^2}{\bar N}},
  \end{equation}
  for any $n\in\{\lceil N/4\rceil,\dots,\lfloor3N/4\rfloor\}$.
    
  \subparagraph{Proof of Step 4.}

  Assume that $n\in\{\bar N+1,\dots,\lfloor 3N/4\rfloor\}$. 
  By \eqref{eq:gau2:stp2:2} there exists $B_1>0$ such that 
  \begin{equation}
    \label{eq:gau2:stp4:1}
    \log \gamma_\Lambda^N(n)-\log \gamma_\Lambda^N(\bar N)=
    \sum_{k=\bar N}^{n-1} \left[\log \gamma_\Lambda^N(k+1)-\log\gamma_\Lambda^N(k)\right]\geq
    -B_1\sum_{k=\bar N}^{n-1}\left(\frac{1}{N}+\frac{k-\bar N}{\bar N}\right).
  \end{equation}
  Using the fact that $0\leq n-\bar N\leq N/2$, and some elementary computation we obtain
  \begin{multline*}
    \sum_{k=\bar N}^{n-1}\left(\frac{1}{N}+\frac{k-\bar N}{\bar N}\right)\leq
    \frac{n-\bar N}{N}+\frac{1}{\bar N}\sum_{k=1}^{n-\bar N-1}k\leq
    \frac{1}{2}+\frac{(n-\bar N)(n-\bar N-1)}{ N}\\
    =\frac{1}{2}+\frac{(n-\bar N)^2}{N}-\frac{n-\bar N}{N}\leq
    \frac{1}{2}+\frac{(n-\bar N)^2}{N},
  \end{multline*}
  which plugged into \eqref{eq:gau2:stp4:1} implies \eqref{eq:gau2:stp4} for $n\in\{\bar N+1,\dots,\lfloor 3N/4\rfloor\}$.
  The general case can be obtained again by symmetry.

  \paragraph{Step 5.}

    There exists $A_5>0$ such that if $N\in\natural\setminus\{0\}$ then
  \begin{equation}
    \label{eq:gau2:stp5}
    \frac{\gamma_\Lambda^N(n)}{\gamma_\Lambda^N(\bar N)}\geq
    \frac{1}{A_5}e^{-A_5\frac{(n-\bar N)^2}{\bar N}},
  \end{equation}
  for any $n\in\{0,\dots,\lceil N/4\rceil-1\}\cup\{\lfloor 3N/4\rfloor+1,\dots,N\}$.
    
  \subparagraph{Proof of Step 5.}

  Observe that $n\in\{0,\dots,\lceil N/4\rceil-1\}\cup\{\lfloor 3N/4\rfloor+1,\dots,N\}$ implies $|n-\bar N|>N/8$.
  Thus to proof \eqref{eq:gau2:stp5} we have only to show that the ratio $\gamma_\Lambda^N(n)/\gamma_\Lambda^N(\bar N)$ is bounded from below by a negative exponential of $N$.
  Assume that $n\in\{\lfloor3N/4\rfloor+1,\dots,N\}$.
  Then
  \begin{equation}
    \label{eq:gau2:stp5:0}
    \frac{\gamma_\Lambda^N(n)}{\gamma_\Lambda^N(\bar N)}=
    \frac{\gamma_\Lambda^N(\lfloor 3N/4\rfloor)}{\gamma_\Lambda^N(\bar N)}\prod_{k=\lfloor 3N/4\rfloor}^{n-1}\frac{\gamma_\Lambda^N(k+1)}{\gamma_\Lambda^N(k)}.
  \end{equation}
  By \eqref{eq:gau2:stp4} there exists $B_1>0$ such that
  \begin{displaymath}
    \frac{\gamma_\Lambda^N(\lfloor 3N/4\rfloor)}{\gamma_\Lambda^N(\bar N)}\geq
    \frac{1}{B_1}e^{-\frac{B_1(\lfloor 3N/4\rfloor-\bar N)^2}{\bar N}}.
  \end{displaymath}
  But $|\lfloor 3N/4\rfloor-\bar N|\leq N$, thus
  \begin{equation}
    \label{eq:gau2:stp5:1}
    \frac{\gamma_\Lambda^N(\lfloor 3N/4\rfloor)}{\gamma_\Lambda^N(\bar N)}\geq
    \frac{1}{B_1}e^{-B_1N}.
  \end{equation}
  In order to bound the product factor in the right hand side of \eqref{eq:gau2:stp5:0}, notice that by Proposition~\ref{pro:dec} there exists $B_2>0$ such that
  \begin{displaymath}
    \frac{\gamma_\Lambda^N(k+1)}{\gamma_\Lambda^N(k)}\geq
    \frac{N-k}{B_2(k+1)},
  \end{displaymath}
  for any $k\in\{1,\dots,N-1\}$, thus
  \begin{multline*}
    \sum_{k=\lfloor 3N/4\rfloor}^{n-1}\left[\log\gamma_\Lambda^N(k)-\log\gamma_\Lambda^N(k+1)\right]\\
    \leq (n-\lfloor 3N/4\rfloor)\log B_2+\sum_{k=\lfloor 3N/4\rfloor}^{n-1}\left[\log(k+1)-\log(N-k)\right]\\
    \leq N\log B_2+\sum_{k=\lfloor 3N/4\rfloor}^{n-1}\frac{2k+1-N}{N-k}\leq
     N\log B_2+(n-\lfloor 3N/4\rfloor)\frac{2(n-1)+1-N}{N-3N/4}\\
    \leq N\log B_2+N\frac{N-1}{N-3N/4}\leq
    N(\log B_2+4)
  \end{multline*}
  which implies
  \begin{displaymath}
    \prod_{k=\lfloor 3N/4\rfloor}^{n-1}\frac{\gamma_\Lambda^N(k+1)}{\gamma_\Lambda^N(k)}\geq
    e^{-(\log B_2+4)N}.
  \end{displaymath}
  By plugging this bound and \eqref{eq:gau2:stp5:1} into \eqref{eq:gau2:stp5:0}, we get \eqref{eq:gau2:stp5}.

  \paragraph{Step 6.}

  There exists $A_6>0$ such that if $N\in\natural\setminus\{0\}$, then,
  \begin{equation}
    \label{eq:gau2:stp6}
    \frac{1}{A_6\sqrt{\bar N}}\leq
    \gamma_\Lambda^N(\bar N)\leq
    \frac{A_6}{\sqrt{\bar N}}.
  \end{equation}

  \subparagraph{Proof of Step 6.}

  By \eqref{eq:gau2:stp3}, \eqref{eq:gau2:stp4} and \eqref{eq:gau2:stp5} we obtain $B_1>0$ such that for any $N\in\natural\setminus\{0\}$
  \begin{displaymath}
    \frac{1}{B_1}e^{-\frac{B_1(n-\bar N)^2}{\bar N}}\leq
    \frac{\gamma_\Lambda^N(n)}{\gamma_\Lambda^N(\bar N)}\leq
    B_1e^{-\frac{(n-\bar N)^2}{B_1\bar N}}
  \end{displaymath}
  for any $n\in\{0,\dots,N\}$.
  By summing for $n\in\{0,\dots,N\}$ we have
   \begin{displaymath}
     \sum_{n=0}^N\frac{1}{B_1}e^{-\frac{B_1(n-\bar N)^2}{\bar N}}\leq
    \frac{1}{\gamma_\Lambda^N(\bar N)}\leq
     B_1\sum_{n=0}^Ne^{-\frac{(n-\bar N)^2}{B_1\bar N}},
   \end{displaymath}
   Thus we are done if we can show that, for any $B_2>0$, there exists $B_3>0$ such that
   \begin{equation}
     \label{eq:gau2:stp6:1}
     \frac{\sqrt{\bar N}}{B_3}\leq
     \sum_{n=0}^Ne^{-\frac{(n-\bar N)^2}{B_2\bar N}}\leq
     B_3\sqrt{\bar N}
   \end{equation}
   for any $N\in\natural\setminus\{0\}$. 
   An elementary computation gives
   \begin{multline*}
     \sum_{n=0}^{N}e^{-B_2\frac{(n-\bar N)^2}{\bar N}}\leq
     \sum_{n=0}^{\bar N}e^{-B_2\frac{(n-\bar N)^2}{\bar N}}+\sum_{n=\bar N}^{N}e^{-B_2\frac{(n-\bar N)^2}{\bar N}}\\
     \leq\int_{-\infty}^{\bar N}e^{-B_2\frac{(x-\bar N)^2}{\bar N}}\; dx+\int_{\bar N}^{+\infty}e^{-B_2\frac{(x-\bar N)^2}{\bar N}}\; dx +1=
     1+\sqrt{\frac{\pi\bar N}{B_2}}.
   \end{multline*}
   that proves the upper bound in \eqref{eq:gau2:stp6:1}.
   The proof of the lower bound is similar.

   \paragraph{Conclusion}

   By  \eqref{eq:gau2:stp3}, \eqref{eq:gau2:stp4}, \eqref{eq:gau2:stp5} and \eqref{eq:gau2:stp6} we get \eqref{eq:gau2}.
 \end{proof}

We now prove an iterative procedure that allow us to extend the Gaussian bounds of Lemma~\ref{lemma:gau2}, from $|\Lambda|=2$ to a generic $\Lambda\subset\subset\integer^d$.

\begin{lem}
  \label{lemma:prec}
  Let $\Lambda^\prime\subset\Lambda\subset\subset\integer^d$ be such that $2\leq|\Lambda^\prime|<|\Lambda|$ and $N\in\natural\setminus\{0\}$;
  assume $\Lambda^\prime=\Lambda_1^\prime\cup\Lambda_2^\prime$ with $\Lambda_1^\prime\cap\Lambda_2^\prime=\emptyset$, $\Lambda_i^\prime\not=\emptyset$ for any $i\in\{1,2\}$ and $\Lambda^{\prime\prime}:=\Lambda\setminus\Lambda^\prime=\Lambda_1^{\prime\prime}\cup\Lambda_2^{\prime\prime}$ with $\Lambda_1^{\prime\prime}\cap\Lambda_2^{\prime\prime}=\emptyset$;
  define $\Lambda_i:=\Lambda_i^\prime\cup\Lambda_i^{\prime\prime}$ for any $i\in\{1,2\}$.
  Then
  \begin{equation}
    \label{eq:lemma:prec}
    \nu_\Lambda^N[\n\eta{\Lambda^\prime}=n]=
    \sum_{k=0}^N\nu_{\Lambda}^N[\n\eta{\Lambda_1}=k]\sum_{h=0\vee[n-(N-k)]}^{k\wedge n}\nu_{\Lambda_1}^k[\n\eta{\Lambda_1^\prime}=h]\nu_{\Lambda_2}^{N-k}[\n\eta{\Lambda_2^\prime}=n-h]
  \end{equation}
  for any $n\in\{0,1,\dots,N\}$.
\end{lem}

\begin{proof}
  For $N\in\natural\setminus\{0\}$ and $n\in\{0,1,\dots,N\}$, by definition of conditional probability
  \begin{multline}
    \label{eq:lemma:prec:1}
    \nu_\Lambda^N[\n\eta{\Lambda^\prime}=n]=
    \nu_\Lambda^N[\n\eta{\Lambda_1^\prime\cup\Lambda_2^\prime}=n]=
    \sum_{k=0}^N\nu_\Lambda^N[\n\eta{\Lambda_1^\prime\cup\Lambda_2^\prime}=n,\ \n\eta{\Lambda_1}=k ]\\
    =\sum_{k=0}^N\nu_\Lambda^N[\n\eta{\Lambda_1^\prime\cup\Lambda_2^\prime}=n| \n\eta{\Lambda_1}=k ]\nu_\Lambda^N[ \n\eta{\Lambda_1}=k ].
  \end{multline}
  But
  \begin{multline}
    \label{eq:lemma:prec:2}
    \nu_\Lambda^N[\n\eta{\Lambda_1^\prime\cup\Lambda_2^\prime}=n| \n\eta{\Lambda_1}=k ]=
    \sum_{h=0}^{k}\nu_\Lambda^N[\n\eta{\Lambda_1^\prime}=h,\ \n\eta{\Lambda_2^\prime}=n-h| \n\eta{\Lambda_1}=k ]\\
    =\sum_{h=0\vee[n-(N-k)]}^{k\wedge n}\nu_\Lambda^N[\n\eta{\Lambda_1^\prime}=h,\ \n\eta{\Lambda_2^\prime}=n-h| \n\eta{\Lambda_1}=k ]
  \end{multline}
  and
  \begin{displaymath}
    \nu_\Lambda^N[\n\eta{\Lambda_1^\prime}=h,\ \n\eta{\Lambda_2^\prime}=n-h| \n\eta{\Lambda_1}=k ]=
    \nu_{\Lambda_1}^k[\n\eta{\Lambda_1^\prime}=h]\nu_{\Lambda_2}^{N-k}[\n\eta{\Lambda_2^\prime}=n-h].
  \end{displaymath}
  By plugging this identity into \eqref{eq:lemma:prec:2} and \eqref{eq:lemma:prec:2} into \eqref{eq:lemma:prec:1} we get \eqref{eq:lemma:prec}.
\end{proof}

\begin{pro}
  \label{pro:gaul}
    Let $\Lambda^\prime\subset\Lambda\subset\subset\integer^d$ be such that $|\Lambda^\prime|=|\Lambda\setminus\Lambda^\prime|$ and define, for any $N\in\natural\setminus\{0\}$, $\bar N:=\lceil N/2 \rceil$;
  then there exist a positive constant $A_0$
    such that for any $N\in\natural\setminus\{0\}$
  \begin{equation}
    \label{eq:pro:gaul}
    \frac{1}{A_0\sqrt{\bar N}}e^{-\frac{A_0(n-\bar N)^2}{\bar N}}\leq
    \gamma_\Lambda^N(n)\leq
    \frac{A_0}{\sqrt{\bar N}}e^{-\frac{(n-\bar N)^2}{A_0\bar N}},
  \end{equation}
  for any $n\in\{0,1,\dots,N\}$.
\end{pro}

\begin{proof}
  We will prove the result by induction on $|\Lambda|$.
  The case $|\Lambda|=2$ has been proved in Lemma~\ref{lemma:gau2}, so assume that \eqref{eq:pro:gaul} has been proved for any $\Lambda\subset\subset\integer^d$ such that $|\Lambda|=2v$ and $v\leq v_0$ where $v,v_0\in\natural\setminus\{0\}$.
  We now show that \eqref{eq:pro:gaul} holds for any $\Lambda\subset\subset\integer^d$ such that $|\Lambda|=2(v_0+1)$.
  
  Let $\Lambda\subset\subset\integer^d$ be such that $|\Lambda|=2(v_0+1)$.
  Then take $\Lambda^\prime\subset\Lambda$ such that $|\Lambda^\prime|=v_0+1$, furthermore there exist $\Lambda_1^\prime,\Lambda_2^\prime\subset\Lambda^\prime$ and $\Lambda_1^{\prime\prime},\Lambda_2^{\prime\prime}\subset\Lambda\setminus\Lambda^\prime$ such that the following holds:
  $\Lambda^\prime=\Lambda_1^\prime\cup\Lambda_2^\prime$, $\Lambda_1^\prime\cap\Lambda_2^\prime=\emptyset$, $\Lambda_i^\prime\not=\emptyset$ for any $i\in\{1,2\}$;
  $\Lambda\setminus\Lambda^\prime=\Lambda_1^{\prime\prime}\cup\Lambda_2^{\prime\prime}$ with $\Lambda_1^{\prime\prime}\cap\Lambda_2^{\prime\prime}=\emptyset$;
  $|\Lambda_1|$ and $|\Lambda_2|$ are even, where $\Lambda_i:=\Lambda_i^\prime\cup\Lambda_i^{\prime\prime}$ for any $i\in\{1,2\}$.
  For instance it suffices to take disjoint nonempty lattice sets $\Lambda_1^\prime, \Lambda_2^\prime, \Lambda_1^{\prime\prime}, \Lambda_2^{\prime\prime}\subset\Lambda$ such that $\Lambda^\prime=\Lambda_1^\prime\cup\Lambda_2^\prime$ and $\Lambda\setminus\Lambda^\prime=\Lambda_1^{\prime\prime}\cup\Lambda_2^{\prime\prime}$, with $|\Lambda_1^\prime|=|\Lambda_2^\prime|=|\Lambda_1^{\prime\prime}|=|\Lambda_2^{\prime\prime}|=(v_0+1)/2$ if $v_0$ is odd, and $|\Lambda_1^\prime|=|\Lambda_1^{\prime\prime}|=v_0/2$, $|\Lambda_2^\prime|=|\Lambda_2^{\prime\prime}|=(v_0/2)+1$ if $v_0$ is even.
  Thus by Lemma~\ref{lemma:prec}
  \begin{equation}
    \label{eq:pro:gaul:1}
        \nu_\Lambda^N[\n\eta{\Lambda^\prime}=n]=
    \sum_{k=0}^N\nu_{\Lambda}^N[\n\eta{\Lambda_1}=k]\sum_{h=0\vee[n-(N-k)]}^{k\wedge n}\nu_{\Lambda_1}^k[\n\eta{\Lambda_1^\prime}=h]\nu_{\Lambda_2}^{N-k}[\n\eta{\Lambda_2^\prime}=n-h],
  \end{equation}
  where $\Lambda_i:=\Lambda_i^\prime\cup\Lambda_i^{\prime\prime}$, for $i\in\{1,2\}$.
  Since $|\Lambda_i|\leq 2v_0$,
  and rewriting \eqref{eq:pro:gaul:1} as
  \begin{equation}
    \label{eq:pro:gaul:2}
            \gamma_\Lambda^N(n)=
    \sum_{k=0}^N\gamma_{\Lambda}^N(k)\sum_{h=0\vee[n-(N-k)]}^{k\wedge n}\gamma_{\Lambda_1}^k(h)\gamma_{\Lambda_2}^{N-k}(n-h),
  \end{equation}
  we can use the inductive assumption \eqref{eq:pro:gaul} on $\gamma_{\Lambda_1}^k(\cdot)=\nu_{\Lambda_1}^k[\n\eta{\Lambda_1^\prime}=\cdot]$ and $\gamma_{\Lambda_2}^{N-k}(\cdot)=\nu_{\Lambda_2}^{N-k}[\n\eta{\Lambda_2^\prime}=\cdot]$ to bound $\gamma_\Lambda^N$.

    First of all notice that because $\gamma_{\Lambda_1}^0(0)=\gamma_{\Lambda_2}^0(0)=1$ then
  \begin{multline}
    \label{eq:pro:gaul:4}
    \sum_{k=0}^N\gamma_{\Lambda}^N(k)\sum_{h=0\vee[n-(N-k)]}^{k\wedge n}\gamma_{\Lambda_1}^k(h)\gamma_{\Lambda_2}^{N-k}(n-h)\\
    =\gamma_{\Lambda}^N(0)\gamma_{\Lambda_2}^{N}(n)+
    \gamma_{\Lambda}^N(N)\gamma_{\Lambda_1}^{N}(n)+
    \sum_{k=1}^{N-1}\gamma_{\Lambda}^N(k)\sum_{h=0\vee[n-(N-k)]}^{k\wedge n}\gamma_{\Lambda_1}^k(h)\gamma_{\Lambda_2}^{N-k}(n-h).
  \end{multline}
  Then, by \eqref{eq:pro:gaul:2} and the inductive assumption there exists $B_1>0$ such that
  \begin{multline*}
    \sum_{h=0\vee[n-(N-k)]}^{k\wedge n}\gamma_{\Lambda_1}^k(h)\gamma_{\Lambda_2}^{N-k}(n-h)\\
    \leq B_1^2\sum_{h=0\vee[n-(N-k)]}^{k\wedge n}\frac{1}{\sqrt{\bar k(\bar N -\bar k)}}\exp\left(-\left\{\frac{(h-\bar k)^2}{B_1\bar k}+\frac{[n-h-(\bar N-\bar k)]^2}{B_1(\bar N-\bar k)}\right\}\right)\\
    \leq B_1^2\sum_{h=-\infty}^{+\infty}\frac{1}{\sqrt{\bar k(\bar N -\bar k)}}\exp\left(-\left\{\frac{(h-\bar k)^2}{B_1\bar k}+\frac{[n-h-(\bar N-\bar k)]^2}{B_1(\bar N-\bar k)}\right\}\right)
  \end{multline*}
  for any $k\in\{1,\dots,N-1\}$.
  It is elementary to show, by comparison with the similar property of the normal density, that there exists  $B_2>0$ such that
  \begin{displaymath}
    B_1^2\sum_{h=-\infty}^{+\infty}\frac{1}{\sqrt{\bar k(\bar N -\bar k)}}\exp\left(-\left\{\frac{(h-\bar k)^2}{B_1\bar k}+\frac{[n-h-(\bar N-\bar k)]^2}{B_1(\bar N-\bar k)}\right\}\right)\leq
    \frac{B_2}{\sqrt{\bar N}}\exp\left({-\frac{(n-\bar N)^2}{B_2\bar N}}\right).
  \end{displaymath}
  Thus by \eqref{eq:pro:gaul:4} and again the inductive assumption
  \begin{multline*}
    \sum_{k=0}^N\gamma_{\Lambda}^N(k)\sum_{h=0\vee[n-(N-k)]}^{k\wedge n}\gamma_{\Lambda_1}^k(h)\gamma_{\Lambda_2}^{N-k}(n-h)\\
    \leq\gamma_{\Lambda}^N(0)\gamma_{\Lambda_2}^{N}(n)+
    \gamma_{\Lambda}^N(N)\gamma_{\Lambda_1}^{N}(n)+
    \frac{B_2}{\sqrt{\bar N}}e^{-\frac{(n-\bar N)^2}{B_2\bar N}}\sum_{k=1}^{N-1}\gamma_{\Lambda}^N(k)\leq
    \frac{3B_2}{\sqrt{\bar N}}e^{-\frac{(n-\bar N)^2}{B_2\bar N}}\leq
    \frac{3B_2}{\sqrt{\bar N}}e^{-\frac{(n-\bar N)^2}{3B_2\bar N}},
  \end{multline*}
  which completes the upper bound in \eqref{eq:pro:gaul}.

   Notice that if $n\in\{0,\dots,\lceil N/4\rceil-1\}\cup\{\lfloor 3N/4\rfloor+1,\dots,N\}$ then $|n-\bar N|> \bar N/8$.
   Thus in this case in order to prove the lower bound in \eqref{eq:pro:gaul} we have to show that $\gamma_\Lambda^N(n)$ can be bounded from below by a negative exponential of $N$.
  This is easy to show; in fact  by the inductive assumption there exists $B_1>0$ such that $\gamma_{\Lambda_i}^s(t)\geq e^{-B_1s}$ for any $i\in\{1,2\}$, $s\in\natural\setminus\{0\}$ and $t\in\{1,\dots,s\}$.
  So, by \eqref{eq:pro:gaul:4},
  \begin{multline*}
    \gamma_{\Lambda}^N(n)=
    \gamma_{\Lambda}^N(0)\gamma_{\Lambda_2}^{N}(n)+
    \gamma_{\Lambda}^N(N)\gamma_{\Lambda_1}^{N}(n)+
    \sum_{k=1}^{N-1}\gamma_{\Lambda}^N(k)\sum_{h=0\vee[n-(N-k)]}^{k\wedge n}\gamma_{\Lambda_1}^k(h)\gamma_{\Lambda_2}^{N-k}(n-h)\\
    \geq\gamma_{\Lambda}^N(0)e^{-B_1N}+\gamma_{\Lambda}^N(N)e^{-B_1N}+e^{-B_1N}\sum_{k=1}^{N-1}\gamma_{\Lambda}^N(k)=
    e^{-B_1N},
  \end{multline*}
  for any $N\in\natural\setminus\{0\}$ and $n\in\{0,\dots,N\}$.

  Assume now that $n\in\{\lceil N/4\rceil,\dots,\lfloor 3N/4\rfloor\}$.
  By symmetry, without loss of generality, we may assume
  $n\in\{\bar N,\dots,\lfloor 3N/4\rfloor\}$.
  In this case, by \eqref{eq:pro:gaul:2}, 
  \begin{equation}
    \label{eq:pro:gaul:6}
    \gamma_\Lambda^N(n)\geq
    \sum_{k=\lceil N/4\rceil}^{\lfloor 3N/4\rfloor}\gamma_{\Lambda}^N(k)\sum_{h=0\vee[n-(N-k)]}^{k\wedge n}\gamma_{\Lambda_1}^k(h)\gamma_{\Lambda_2}^{N-k}(n-h).
  \end{equation}
  By the inductive assumption there exists $B_2>0$ such that
  \begin{equation}
    \label{eq:pro:gaul:7}
    \gamma_{\Lambda_1}^k(h)\gamma_{\Lambda_2}^{N-k}(n-h)\geq
    \frac{1}{B_2^2\sqrt{ \frac{k}{2}\left(\frac{N- k}{2}\right)}}\exp\left(-B_2\left\{\frac{(h- k/2)^2}{k/2}+\frac{[n-h-(N- k)/2]^2}{(N -k)/2}\right\}\right),
  \end{equation}
  and an elementary calculation shows that
  \begin{multline}
    \label{eq:pro:gaul:8}
    \frac{1}{B_2^2\sqrt{ \frac{k}{2}\left(\frac{N- k}{2}\right)}}\exp\left(-B_2\left\{\frac{(h- k/2)^2}{k/2}+\frac{[n-h-(N- k)/2]^2}{(N -k)/2}\right\}\right)\\
        =\frac{1}{B_2\sqrt{\frac{N}{2}}}\exp\left(-B_2\frac{(n- N/2)^2}{ N/2}\right)
        \frac{1}{B_2\sqrt{\frac{(N-k)k}{2N}}}\exp\left(-B_2\frac{(h-kn/N)^2}{(N-k)k/(2N)}\right).
  \end{multline}
  Thus we have only to show that there exists $B_3>0$ such that
  \begin{equation}
    \label{eq:pro:gaul:9}
    \sum_{h=0\vee[n-(N-k)]}^{k\wedge n}\frac{1}{B_2\sqrt{\frac{(N-k)k}{2N}}}\exp\left(-B_2\frac{(h-kn/N)^2}{(N-k)k/(2N)}\right)\geq \frac{1}{B_3}
  \end{equation}
  for any $N\in\natural\setminus\{0\}$, and $n\in\{\bar N,\dots,\lfloor 3N/4\rfloor\}$.
  In fact assume that \eqref{eq:pro:gaul:9} holds;
  then by \eqref{eq:pro:gaul:8}, \eqref{eq:pro:gaul:7} and \eqref{eq:pro:gaul:6} there exists $B_4>0$ such that
  \begin{multline*}
    \gamma_\Lambda^N(n)\geq
    \frac{1}{B_4\sqrt{\bar N}}e^{-B_4\frac{(n- \bar N)^2}{\bar N}}\sum_{k=\lceil N/4\rceil}^{\lfloor 3N/4\rfloor}\gamma_\Lambda^N(k)\\
    \geq\frac{1}{B_4\sqrt{\bar N}}e^{-B_4\frac{(n- \bar N)^2}{\bar N}}\left[1-\sum_{k=0}^{\lceil N/4\rceil-1}\gamma_\Lambda^N(k)+\sum_{\lfloor 3N/4\rfloor+1}^{N}\gamma_\Lambda^N(k)\right].
  \end{multline*}
  But we just proved that there exists $B_1>0$ such that $\gamma_\Lambda^N(n)\geq e^{-B_1N}$ for any $N\in\natural\setminus\{0\}$ and $n\in\{0,\dots,N\}$.
  Thus
  \begin{displaymath}
    \sum_{k=0}^{\lceil N/4\rceil-1}\gamma_\Lambda^N(k)+\sum_{\lfloor 3N/4\rfloor+1}^{N}\gamma_\Lambda^N(k)\leq
    \frac{Ne^{-B_1N}}{2}\leq
    \frac{1}{2}
  \end{displaymath}
  for $N$ large enough.
  This proves that
    \begin{displaymath}
    \gamma_\Lambda^N(n)\geq
    \frac{1}{2B_4\sqrt{\bar N}}e^{-2B_4\frac{(n- \bar N)^2}{\bar N}},
  \end{displaymath}
  for any $N\in\natural\setminus\{0\}$ large enough and $n\in\{\bar N,\dots,\lfloor 3N/4\rfloor\}$;
  the case of small $N$ is controlled by a finiteness argument.

  It remains to prove \eqref{eq:pro:gaul:9}, which is an elementary but tedious calculation:
  first of all observe that if $N\in\natural\setminus\{0\}$, $k\in\{\lceil N/4\rceil,\dots,\lfloor 3N/4\rfloor\}$ and $n\in\{\bar N,\dots,\lfloor 3N/4\rfloor\}$, then
  \begin{displaymath}
    (k\wedge n) -\frac{kn}{N}=
    kn\left(\frac{k\wedge n}{kn} -\frac{1}{N}\right)=
    kn\left(\frac{1}{k\vee n} -\frac{1}{N}\right)\geq
    \frac{N^2}{8}\left(\frac{4}{3N} -\frac{1}{N}\right)=
    \frac{N}{24}.
  \end{displaymath}
  Furthermore, since
  \begin{displaymath}
    n-(N-k)=
    n+k-N\leq
    \frac{3}{2}N-N=
    \frac{N}{2},
  \end{displaymath}
  then
  \begin{displaymath}
    0\vee[n-(N-k)]-\frac{kn}{N}=
    kn\left\{\frac{0\vee[n-(N-k)]}{kn}-\frac{1}{N}\right\}\leq
    \frac{N^2}{8}\left(\frac{\frac{N}{2}}{\frac{9N^2}{16}}-\frac{1}{N}\right)=
    -\frac{N}{72}.
  \end{displaymath}
  This implies
  \begin{multline*}
    \sum_{h=0\vee[n-(N-k)]}^{k\wedge n}\frac{1}{B_2\sqrt{(N-k)k/(2N)}}\exp\left(-B_2\frac{(h-kn/N)^2}{(N-k)k/(2N)}\right)\\
    \geq\sum_{h=\lceil kn/N\rceil }^{k\wedge n}\frac{1}{B_2\sqrt{(N-k)k/(2N)}}\exp\left(-B_2\frac{(h-kn/N)^2}{(N-k)k/(2N)}\right)\\
    \geq \frac{1}{B_2\sqrt{(N-k)k/(2N)}}\int_{\lceil kn/N\rceil }^{k\wedge n}\exp\left(-B_2\frac{(h-kn/N)^2}{(N-k)k/(2N)}\right)\;dh\\
    =\sqrt{\frac{\pi}{B_2}}\left\{\Phi\left(\frac{\sqrt{2B_2}\left[(k\wedge n)-kn/N\right]}{\sqrt{(N-k)k/(2N)}}\right)-\Phi\left(\frac{\sqrt{2B_2}\left(\lceil kn/N\rceil-kn/N\right)}{\sqrt{(N-k)k/(2N)}}\right)\right\},
  \end{multline*}
  where
  \begin{displaymath}
    \Phi(x):=\frac{1}{\sqrt{2\pi}}\int_{-\infty}^{x}e^{-\frac{z^2}{2}}\; dz.
  \end{displaymath}
  Notice that, since $k\in\{\lceil N/4\rceil,\dots,\lfloor 3N/4\rfloor\}$, it is easy to show that
  \begin{displaymath}
    \frac{\left[(k\wedge n)-kn/N\right]}{\sqrt{\frac{(N-k)k}{2N}}}\geq
    \frac{\sqrt{N}}{6\sqrt{2}}
    \qquad\text{and}\qquad
    \frac{\left(\lceil kn/N\rceil-kn/N\right)}{\sqrt{(N-k)k/(2N)}}\leq
  \frac{1}{16\sqrt{3}\sqrt{N}}.
  \end{displaymath}
  Thus
  \begin{multline*}
    \Phi\left(\frac{\sqrt{2B_2}\left[(k\wedge n)-kn/N\right]}{\sqrt{(N-k)k/(2N)}}\right)-\Phi\left(\frac{\sqrt{2B_2}\left(\lceil kn/N\rceil-kn/N\right)}{\sqrt{(N-k)k/(2N)}}\right)\\
    \geq\Phi\left(\frac{\sqrt{N}}{6\sqrt{2}}\right)-\Phi\left(\frac{1}{16\sqrt{3}\sqrt{N}}\right)
    \geq\Phi\left(\frac{1}{6\sqrt{2}}\right)-\Phi\left(\frac{1}{16\sqrt{3}}\right):=
    \frac{1}{B_3},
  \end{multline*}
    which proves \eqref{eq:pro:gaul:9} and completes the proof of the lower bound in \eqref{eq:pro:gaul}.
\end{proof}

\subsection{Further estimates on the grand canonical measure}

The next result is a uniform estimate on $\mu_\rho[\n\eta\Lambda=\rho|\Lambda|]$, see \eqref{eq:cen} below, and will be used in Section~\ref{sec:odlsi}.
If $c(k)=k$ the result can be obtained elementarily from the Stirling Formula.
The general case is more difficult. 

\begin{pro}
  \label{pro:cen}
  \begin{equation}
    \label{eq:cen}
    0<
    \inf_{\substack{\Lambda\subset\subset\integer^d\\ N\in\natural\setminus\{0\}}}\sqrt{\sigma^2(N/|\Lambda|)|\Lambda|}p_\Lambda^{N/|\Lambda|}(N)\leq
        \sup_{\substack{\Lambda\subset\subset\integer^d\\ N\in\natural\setminus\{0\}}}\sqrt{\sigma^2(N/|\Lambda|)|\Lambda|}p_\Lambda^{N/|\Lambda|}(N)<
        +\infty
  \end{equation}
\end{pro}

\begin{proof}
  Fix $\rho_0>0$ and consider five different cases.

    \paragraph{Small density case.}
      There exists $N_0>0$ such that
  \begin{displaymath}
    0<
    \inf_{\substack{N_0\leq N\leq\rho_0|\Lambda|\\ \Lambda\subset\subset\integer^d}}\sqrt{\sigma^2(N/|\Lambda|)|\Lambda|}p_\Lambda^{N/|\Lambda|}(N)\leq
        \sup_{\substack{N_0\leq N\leq\rho_0|\Lambda|\\ \Lambda\subset\subset\integer^d}}\sqrt{\sigma^2(N/|\Lambda|)|\Lambda|}p_\Lambda^{N/|\Lambda|}(N)<
        +\infty
  \end{displaymath}

\subparagraph{Proof of small density case.}

  By point 1. of Proposition~\ref{pro:LLT} there exist positive
  constants $B_0$ and $n_0$ such that:
   \begin{displaymath}
   \left|\sqrt{\sigma^2(N/|\Lambda|)|\Lambda|}p_{\Lambda}^{N/|\Lambda|}(N)-\frac{1}{\sqrt{2 \pi}}\right|\leq
      \frac{B_0}{\sqrt{\sigma^2(N/|\Lambda|)|\Lambda|}}
  \end{displaymath}
  uniformly in $N\in\natural\setminus\{0\}$ and $\Lambda\subset\subset\integer^d$ such that $N/|\Lambda|\leq\rho_0$ and $\sigma^2(N/|\Lambda|)|\Lambda|\geq n_0$.
  Because of \eqref{eq:sigma} there exists $B_1>0$ such that $\sigma^2(N/|\Lambda|)|\Lambda|\geq B_1^{-1} N$
  for any $N\in\natural\setminus\{0\}$ and $\Lambda\subset\subset\integer^d$.
  So take $N_1:=  n_0B_1$, then, for any $N\geq N_1$, $\sigma^2(N/|\Lambda|)|\Lambda|\geq n_0$, and therefore we have.
     \begin{displaymath}
   \left|\sqrt{\sigma^2(N/|\Lambda|)|\Lambda|}p_{\Lambda}^{N/|\Lambda|}(N)-\frac{1}{\sqrt{2 \pi}}\right|\leq
      \frac{B_0}{\sqrt{\sigma^2(N/|\Lambda|)|\Lambda|}}\leq
      \frac{B_0\sqrt{B_1}}{\sqrt{N}}
  \end{displaymath}
  where $\Lambda\subset\subset\integer^d$ is such that $N/|\Lambda|\leq \rho_0$.
  Taking $N$ large enough, namely $N\geq N_0:= \lceil N_1\vee (8\pi B_0^2 B_1)\rceil$ we have:
     \begin{displaymath}
   \left|\sqrt{\sigma^2(N/|\Lambda|)|\Lambda|}p_{\Lambda}^{N/|\Lambda|}(N)-\frac{1}{\sqrt{2 \pi}}\right|\leq
      \frac{1}{2\sqrt{2\pi}}
  \end{displaymath}
  which means
       \begin{displaymath}
      \frac{1}{2\sqrt{2\pi}}\leq
   \sqrt{\sigma^2(N/|\Lambda|)|\Lambda|}p_{\Lambda}^{N/|\Lambda|}(N)\leq
      \frac{3}{2\sqrt{2\pi}}
  \end{displaymath}
  for any $N\geq N_0$ and any $\Lambda\subset\subset\integer^d$ such that $N/|\Lambda|\leq \rho_0$.

  \paragraph{Very small density case.}

  For any fixed $N_0\in\natural\setminus\{0\}$
  \begin{displaymath}
    0<
    \inf_{\substack{0< N\leq N_0 \\ \Lambda\subset\subset\integer^d}}\sqrt{\sigma^2(N/|\Lambda|)|\Lambda|}p_\Lambda^{N/|\Lambda|}(N)\leq
        \sup_{\substack{0< N\leq N_0 \\ \Lambda\subset\subset\integer^d}}\sqrt{\sigma^2(N/|\Lambda|)|\Lambda|}p_\Lambda^{N/|\Lambda|}(N)<
        +\infty
  \end{displaymath}

  \subparagraph{Proof of very small density case.}

  By Lemma~\ref{lemma:poi} there exists $B_2>0$ such that for any $\Lambda\subset\subset\integer^d$
  \begin{displaymath}
        \left|\sqrt{\sigma^2(N/|\Lambda|)|\Lambda|}p_\Lambda^{N/|\Lambda|}(N)-\frac{\sqrt{\sigma^2(N/|\Lambda|)|\Lambda|}N^N}{N!}e^{-N}\right|\leq
    B_2\sqrt{\frac{\sigma^2(N/|\Lambda|)}{|\Lambda|}}
  \end{displaymath}
  uniformly in $N\in\natural\setminus\{0\}$ with $N\leq N_0$.
  Because of \eqref{eq:sigma} there exists $B_1>0$ such that $B_1^{-1}N\leq\sigma^2(N/|\Lambda|)|\Lambda|\leq B_1N$, so we may write
   \begin{displaymath}
    \frac{N^{N+1/2}}{\sqrt{B_1}N!}e^{-N}-\frac{\sqrt{N B_1}B_2}{|\Lambda|}\leq
        \sqrt{\sigma^2(N/|\Lambda|)|\Lambda|}p_\Lambda^{N/|\Lambda|}(N)\leq
    \frac{\sqrt{B_1}N^{N+1/2}}{N!}e^{-N}+\frac{\sqrt{N B_1}B_2}{|\Lambda|}.
  \end{displaymath}
  Now taking $v_1$ large enough
  it easy to show that
    \begin{displaymath}
    0<
    \inf_{\substack{0< N\leq N_0 \\ \Lambda\subset\subset\integer^d,\; |\Lambda|>v_1}}\sqrt{\sigma^2(N/|\Lambda|)|\Lambda|}p_\Lambda^{N/|\Lambda|}(N)\leq
        \sup_{\substack{0< N\leq N_0 \\ \Lambda\subset\subset\integer^d,\; |\Lambda|>v_1}}\sqrt{\sigma^2(N/|\Lambda|)|\Lambda|}p_\Lambda^{N/|\Lambda|}(N)<
        +\infty.
  \end{displaymath}
  The general case follows again by a finiteness argument.

  \paragraph{Normal density case.}
      There exist $v_0>0$ such that
  \begin{displaymath}
    0<
    \inf_{\substack{v_0<|\Lambda|<\rho_0^{-1}N\\ \Lambda\subset\subset\integer^d}}\sqrt{\sigma^2(N/|\Lambda|)|\Lambda|}p_\Lambda^{N/|\Lambda|}(N)\leq
        \sup_{\substack{v_0<|\Lambda|<\rho_0^{-1}N\\ \Lambda\subset\subset\integer^d}}\sqrt{\sigma^2(N/|\Lambda|)|\Lambda|}p_\Lambda^{N/|\Lambda|}(N)<
        +\infty
  \end{displaymath}

  \subparagraph{Proof of normal density case.}

  By point 2. of Proposition~\ref{pro:LLT} there exist positive
  constants $B_3$ and $v_4$ such that:
   \begin{displaymath}
   \left|\sqrt{\sigma^2(N/|\Lambda|)|\Lambda|}p_{\Lambda}^{N/|\Lambda|}(N)-\frac{1}{\sqrt{2 \pi}}\right|\leq
      \frac{B_3}{\sqrt{|\Lambda|}}
  \end{displaymath}
  uniformly in $N\in\natural\setminus\{0\}$ and $\Lambda\subset\subset\integer^d$ such that $N/|\Lambda|>\rho_0$.
  Now it is easy to show that there exists $v_0>0$ such that for any $\Lambda\subset\subset\integer^d$ with $|\Lambda|>v_0$ and any $N\in\natural\setminus\{0\}$ such that $N> \rho_0|\Lambda|$
  \begin{displaymath}
      \frac{1}{2\sqrt{2\pi}}\leq
   \sqrt{\sigma^2(N/|\Lambda|)|\Lambda|}p_{\Lambda}^{N/|\Lambda|}(N)\leq
      \frac{3}{2\sqrt{2\pi}}.
  \end{displaymath}

  \paragraph{Large density case.}

  For any fixed $v_0\in\natural\setminus\{0\}$
  \begin{equation}
    \label{eq:vldc}
    0<
    \inf_{\substack{N\in\natural\setminus\{0\},\; \Lambda\subset\integer^d\\  0<|\Lambda|\leq v_0 }}\sqrt{\sigma^2(N/|\Lambda|)|\Lambda|}p_\Lambda^{N/|\Lambda|}(N)\leq
        \sup_{\substack{N\in\natural\setminus\{0\},\; \Lambda\subset\integer^d\\  0<|\Lambda|\leq v_0 }}\sqrt{\sigma^2(N/|\Lambda|)|\Lambda|}p_\Lambda^{N/|\Lambda|}(N)<
        +\infty
  \end{equation}

  \subparagraph{Proof of large density case.}
  In this case no local limit theorem is applicable.
  Instead we will use the Gaussian estimates of Section~\ref{sec:gau}.
    For any fixed $\rho_0>0$, by point 2. of Proposition~\ref{pro:LLT} there exist positive
  constants $B_1$ and $v_1$ such that:
   \begin{displaymath}
   \left|\sqrt{\sigma^2(N/|\Lambda|)|\Lambda|}p_{\Lambda}^{N/|\Lambda|}(N)-\frac{1}{\sqrt{2 \pi}}\right|\leq
      \frac{B_1}{\sqrt{|\Lambda|}}
  \end{displaymath}
  uniformly in $N\in\natural\setminus\{0\}$ and $\Lambda\subset\subset\integer^d$ such that $N/|\Lambda|>\rho_0$.
  Taking $|\Lambda|$ large enough, namely $|\Lambda|\geq v_2:= 8\pi B_1^2$ we have:
     \begin{displaymath}
   \left|\sqrt{\sigma^2(N/|\Lambda|)|\Lambda|}p_{\Lambda}^{N/|\Lambda|}(N)-\frac{1}{\sqrt{2 \pi}}\right|\leq
      \frac{1}{2\sqrt{2\pi}}
  \end{displaymath}
  which means
       \begin{equation}
         \label{eq:vldc:1}
      \frac{1}{2\sqrt{2\pi}}\leq
   \sqrt{\sigma^2(N/|\Lambda|)|\Lambda|}p_{\Lambda}^{N/|\Lambda|}(N)\leq
      \frac{3}{2\sqrt{2\pi}}
  \end{equation}
  for any $\Lambda\subset\subset\integer^d$ such that $|\Lambda|\geq v_2$ and any $N\in\natural\setminus\{0\}$ such that $N> \rho_0|\Lambda|$.
  Now take $\Lambda^\prime\subset\Lambda\subset\subset\integer^d$ such that $v_2/2\leq |\Lambda^\prime|< v_2$ and $|\Lambda\setminus\Lambda^\prime|=|\Lambda^\prime|$;
  notice that in this case  $|\Lambda|=2|\Lambda^\prime|$.
  By Lemma~\ref{pro:gaul} there exists a $B_2>0$, possibly depending on $v_2$, such that
  \begin{displaymath}
    \frac{1}{B_2\sqrt{\bar N}}\leq
    \gamma_\Lambda^N(\bar N)\leq
    \frac{B_2}{\sqrt{\bar N}},
  \end{displaymath}
  uniformly in $N\in\natural\setminus\{0\}$.
  By definition of $\gamma_\Lambda^N$
    \begin{displaymath}
    \frac{1}{B_2\sqrt{\bar N}}\leq
    \frac{p_{\Lambda^\prime}^{\bar N/|\Lambda^\prime|}(\bar N)^2}{p_\Lambda^{\bar N/|\Lambda^\prime|}( N)}\leq
    \frac{B_2}{\sqrt{\bar N}}
  \end{displaymath}
  \ie
    \begin{displaymath}
    \frac{p_\Lambda^{\bar N/|\Lambda^\prime|}( N)}{B_2\sqrt{\bar N}}\leq
    p_{\Lambda^\prime}^{\bar N/|\Lambda^\prime|}(\bar N)^2\leq
    \frac{B_2p_\Lambda^{\bar N/|\Lambda^\prime|}( N)}{\sqrt{\bar N}}
  \end{displaymath}
  uniformly in $N\in\natural\setminus\{0\}$;
  but $|\Lambda|=2|\Lambda^\prime|\geq v_2$ and $\bar N/|\Lambda^\prime|\geq N/2|\Lambda^\prime|=N/|\Lambda|>\rho_0$, and therefore we can use \eqref{eq:vldc:1} to get
    \begin{displaymath}
    \frac{1}{B_2\sqrt{2\pi}\sqrt{\sigma^2(N/|\Lambda|)|\Lambda|\bar N}}\leq
    p_{\Lambda^\prime}^{\bar N/|\Lambda^\prime|}(\bar N)^2\leq
    \frac{3B_2}{2\sqrt{2\pi}\sqrt{\sigma^2(N/|\Lambda|)|\Lambda|\bar N}},
  \end{displaymath}
  and by using \eqref{eq:sigma} we conclude that there exist a $B_3$, possibly depending on $v_2$, such that
    \begin{displaymath}
    \frac{1}{B_3\sqrt{N}}\leq
    p_{\Lambda^\prime}^{\bar N/|\Lambda^\prime|}(\bar N)\leq
    \frac{B_3}{\sqrt{N}},
  \end{displaymath}
  for any $N\in\natural\setminus\{0\}$ and $\Lambda^\prime\subset\integer^d$ such that $v_2/2\leq |\Lambda^\prime|< v_2$.
  From this inequality, using again \eqref{eq:sigma}, we get immediately a constant $B_4>0$, possibly depending on $v_2$, such that
  \begin{displaymath}
      \frac{1}{B_4}\leq
   \sqrt{\sigma^2(N/|\Lambda|)|\Lambda|}p_{\Lambda}^{N/|\Lambda|}(N)\leq
      B_4
  \end{displaymath}
 for any $\Lambda\subset\subset\integer^d$ such that $|\Lambda|\geq v_2/2$ and any $N\in\natural\setminus\{0\}$ such that $N> \rho_0|\Lambda|$.
 By iterating the above procedure a finite number of times we get \eqref{eq:vldc}.
\end{proof}

\section{One dimensional L.S.I.: proof of Proposition~\ref{ppaolo1}}
\label{sec:odlsi}

In this section we prove Proposition~\ref{ppaolo1}, a logarithmic Sobolev inequality for a particular one dimensional birth and death process.
Furthermore we will see that Proposition~\ref{ppaolo1} implies Proposition~\ref{ppaolo6}.

\subsection{A general result}

Next result is Proposition~A.5. in \cite{Ca:Ma:Ro}.
We report it only for completeness.

Let $\{\gamma^N:N\in\natural\setminus\{0\}\}$ be a family of positive probability on $\integer$.
Assume that for any $N\in\natural\setminus\{0\}$ the probability $\gamma^N$ is supported on $\{0,\dots,N\}$.
It is elementary to check that  $\gamma^N$ is reversible with respect to the continuous time birth and death process with rates
\begin{displaymath}
  a^N(n):=\frac{\gamma^N(n+1)}{\gamma^N(n)}\wedge 1\qquad\text{for birth}
  \qquad
  b^N(n):=\frac{\gamma^N(n-1)}{\gamma^N(n)}\wedge 1\qquad\text{for death.}
\end{displaymath}
for any $N\in\natural\setminus\{0\}$.
The Dirichlet form of this Markov process may be written as
\begin{displaymath}
  \sum_{n=1}^N[\gamma^N(n)\wedge\gamma^N(n-1)][f(n-1)-f(n)]^2.
\end{displaymath}

\begin{pro}
  \label{pro:cmr}
  Assume that there exist  a positive constant $A_0$ such that for any $N\in\natural\setminus\{0\}$ we can find $\bar N\in\{0,\dots,N\}$
  such that $A_0^{-1}\bar N\leq N-\bar N\leq A_0\bar N$ and
  \begin{align}
    \label{eq:cmr:1}
    \frac{\gamma^N(n+1)}{\gamma^N(n)} &\leq e^{-\frac{n-\bar N}{A_0 \bar N}}\qquad\qquad\qquad\ \;\text{for all  $n\in\{\bar N+1,\dots,N\}$} \\
        \label{eq:cmr:2}
    \frac{\gamma^N(n-1)}{\gamma^N(n)} &\leq e^{-\frac{\bar N-n}{A_0 \bar N}}\qquad\qquad\qquad\ \;\text{for all  $n\in\{0,\dots,\bar N-1\}$} \\
        \label{eq:cmr:3}
    \gamma^N(n)  &\geq \frac{1}{A_0\sqrt{\bar N}}e^{-\frac{A_0(n-\bar N)^2}{\bar N}}\qquad\text{for all $n\in\{0,\dots,N\}$}.
  \end{align}
  Then there exists a positive constant $A_1$ such that for any positive function $f$ on $\{0,\dots,N\}$
  \begin{displaymath}
    \ent_{\gamma^N}(f)\leq A_1 N \sum_{n=1}^N[\gamma^N(n)\wedge\gamma^N(n-1)]
[\sqrt{f(n-1)}-\sqrt{f(n)}]^2
  \end{displaymath}
  for any  $N\in\natural\setminus\{0\}$.
\end{pro}

\begin{proof}
  The proof follows closely the proof of Proposition~A.5 in \cite{Ca:Ma:Ro}.
  By Proposition~A.1 of \cite{Ca:Ma:Ro} we have to bound from above $B_0(N):=B_0^-(N)\vee B_0^+(N)$ uniformly in $N\in\natural\setminus\{0\}$, where
  \begin{align*}
    B_0^-(N)&:=
    \sup_{n\in\{0,\dots,\bar N-1\}}\left(\sum_{k=0}^n\gamma^N(k)\right)\log\left(\frac{1}{\sum_{k=0}^n\gamma^N(k)}\right)\left(\sum_{k=n}^{\bar N-1}\frac{1}{\gamma^N(k)\wedge \gamma^N(k+1)}\right)\\
    B_0^+(N)&:=
    \sup_{n\in\{\bar N+1,\dots,N\}}\left(\sum_{k=n}^N\gamma^N(k)\right)\log\left(\frac{1}{\sum_{k=n}^N\gamma^N(k)}\right)\left(\sum_{k=\bar N+1}^{n}\frac{1}{\gamma^N(k)\wedge \gamma^N(k-1)}\right),
  \end{align*}
  and take $A_1:=\sup_NB_0(N)/N$.
  We divide these bounds in several steps.

  \paragraph{Step 1.}
  There exists $B_1>0$ such that for any $N\in\natural\setminus\{0\}$
  \begin{align}
    \label{eq:cmr:4}
    \sum_{k=n}^N\gamma^N(k)&\leq
    \frac{B_1\bar N}{n-\bar N}\gamma^N(n)\qquad\text{for any $n\in\{\bar N+1,\dots,N\}$}\\
    \label{eq:cmr:5}
    \sum_{k=0}^n\gamma^N(k)&\leq
    \frac{B_1\bar N}{\bar N-n}\gamma^N(n)\qquad\text{for any $n\in\{0,\dots,\bar N-1\}$}.
  \end{align}

  \subparagraph{Proof of step 1.}

  For $n\in\{\bar N+1,\dots,N\}$, by a simple telescopic argument and \eqref{eq:cmr:1} we obtain:
  \begin{multline}
    \label{eq:cmr:6}
    \sum_{k=n}^N\frac{\gamma^N(k)}{\gamma^N(n)}=
    1+\sum_{k=n+1}^N\prod_{h=n}^{k-1}\frac{\gamma^N(h+1)}{\gamma^N(h)}\leq
    1+\sum_{k=n+1}^Ne^{-\sum_{h=n}^{k-1}\frac{h-\bar N}{A_0\bar N}}\\
    \leq1+\sum_{k=n+1}^Ne^{-\frac{(k-n)(n-\bar N)}{A_0\bar N}}\leq
    \sum_{k=0}^{+\infty}e^{-\frac{k(n-\bar N)}{A_0\bar N}}=
    \frac{1}{1-e^{-\frac{n-\bar N}{A_0\bar N}}};
  \end{multline}
  $n\in\{\bar N+1,\dots,N\}$.
  Because $N-\bar N\leq A_0\bar N$, it easy to check that
  \begin{displaymath}
    \frac{1}{1-e^{-\frac{n-\bar N}{A_0\bar N}}}\leq
    \frac{2A_0\bar N}{n-\bar N}.
  \end{displaymath}
  By plugging this bound into \eqref{eq:cmr:6} we obtain \eqref{eq:cmr:4}.
  Using the same argument and \eqref{eq:cmr:2} we obtain \eqref{eq:cmr:5}.

\subparagraph{Step 2.}
  There exists $B_2>0$ such that for any $N\in\natural\setminus\{0\}$
  \begin{align}
    \label{eq:cmr:7}
    \sum_{k=\bar N+1}^n\frac{1}{\gamma^N(k)}&\leq
    \frac{B_2\bar N}{n-\bar N}\frac{1}{\gamma^N(n)}\qquad\text{for any $n\in\{\bar N+1,\dots,N\}$}\\
    \label{eq:cmr:8}
    \sum_{k=n}^{\bar N-1}\frac{1}{\gamma^N(k)}&\leq
    \frac{B_2\bar N}{\bar N-n}\frac{1}{\gamma^N(n)}\qquad\text{for any $n\in\{0,\dots,\bar N-1\}$}.
  \end{align}

  \subparagraph{Proof of step 2.}

  For $n\in\{\bar N+1,\dots,N\}$, by the same telescopic argument used in step 1. and  \eqref{eq:cmr:1} we obtain:
  \begin{multline*}
    \sum_{k=\bar N+1}^n\frac{\gamma^N(n)}{\gamma^N(k)}=
    \sum_{k=\bar N+1}^{n-1}\prod_{h=k}^{n-1}\frac{\gamma^N(h+1)}{\gamma^N(h)}+1\leq
    \sum_{k=\bar N+1}^{n-1}e^{-\sum_{h=k}^{n-1}\frac{h-\bar N}{A_0\bar N}}+1\\
    \leq\sum_{k=\bar N+1}^{n-1}e^{-\frac{(n-k)(k-\bar N)}{A_0\bar N}}+1\leq
    \sum_{k=\bar N+1}^{n-1}e^{-\frac{(k-\bar N)}{A_0\bar N}}+1\leq
    \sum_{k=0}^{+\infty}e^{-\frac{k(n-\bar N)}{A_0\bar N}}=
    \frac{1}{1-e^{-\frac{n-\bar N}{A_0\bar N}}},
  \end{multline*}
  and we can conclude as in step 1.

  \paragraph{Step 3.}

  There exist $B_3>0$ such that for any $n\in\{0,\dots,N\}$ with $|n-\bar N|\geq\sqrt{\bar N}$
  \begin{align}
    \label{eq:cmr:9}
    \sum_{k=n}^N\gamma^N(k) &\geq
    \frac{1}{B_3}e^{-\frac{B_3(n-\bar N)^2}{\bar N}},\\
    \label{eq:cmr:900}
    \sum_{k=0}^n\gamma^N(k) &\geq
    \frac{1}{B_3}e^{-\frac{B_3(n-\bar N)^2}{\bar N}}
  \end{align}
  for any $N\in\natural\setminus\{0\}$.

  \subparagraph{Proof of step 3.}

  Assume that $n\in\{0,\dots,N\}$ is such that $n\geq\bar N+\sqrt{\bar N}$, the case $n\leq\bar N-\sqrt{\bar N}$ is similar.
  Because $A_0^{-1}\bar N\leq N-\bar N\leq A_0\bar N$, it is easy to find $n_0\in\natural\setminus\{0\}$ depending only on $A_0$ such that $\bar N+\lceil\sqrt{\bar N}\rceil\leq N-\lceil\sqrt{\bar N}\rceil$ for any $N\geq n_0$.
  We can assume $N\geq n_0$ because for the case $N\in\{1,\dots,n_0\}$  \eqref{eq:cmr:9} follows by a finiteness argument.
  Assume that $n\in\{\bar N+\lceil\sqrt{\bar N}\rceil,\dots,N-\lceil\sqrt{\bar N}\rceil\}$; by condition \eqref{eq:cmr:3} and some simple bounds,
  \begin{multline*}
    \sum_{k=n}^N\gamma_\Lambda^N(k)\geq
    \frac{1}{A_0\sqrt{\bar N}}\sum_{k=n}^Ne^{-\frac{A_0(k-\bar N)^2}{\bar N}}\geq
    \frac{1}{A_0\sqrt{\bar N}}\sum_{k=n}^{n+\lceil\sqrt{\bar N}\rceil}e^{-\frac{A_0(k-\bar N)^2}{\bar N}}\\
   \geq\frac{\lceil\sqrt{\bar N}\rceil+1}{A_0\sqrt{\bar N}}e^{-\frac{A_0(n+\lceil\sqrt{\bar N}\rceil-\bar N)^2}{\bar N}}\geq
   \frac{1}{A_0}e^{-\frac{2A_0(n-\bar N)^2+4A_0\bar N}{\bar N}},
  \end{multline*}
  from which \eqref{eq:cmr:9} follows.
  
  Assume now that $n\in\{N-\lceil\sqrt{\bar N}\rceil+1,\dots, N\}$.
  Then, again by condition \eqref{eq:cmr:3},
    \begin{multline}
      \label{eq:cmr:10}
    \sum_{k=n}^N\gamma_\Lambda^N(k)\geq
    \frac{1}{A_0\sqrt{\bar N}}\sum_{k=n}^Ne^{-\frac{A_0(k-\bar N)^2}{\bar N}}\geq
    \frac{1}{A_0\sqrt{\bar N}}e^{-\frac{A_0(n-\bar N)^2}{\bar N}}\geq
    \frac{1}{A_0\sqrt{\bar N}}e^{-\frac{A_0(n-\bar N)^2}{\bar N}}\\
    \geq\frac{e^{\frac{A_0(N-\lceil\sqrt{\bar N}\rceil-\bar N)^2}{\bar N}}}{\sqrt{A_0\bar N}}e^{-\frac{2A_0(n-\bar N)^2}{\bar N}}.
  \end{multline}
  Note that
  \begin{displaymath}
    (N-\lceil\sqrt{\bar N}\rceil-\bar N)^2\geq
    (N-\bar N)^2-3(N-\bar N)\sqrt{\bar N}\geq
    \frac{\bar N^2}{A_0^2}-3A_0\bar N^{3/2}\geq
    \frac{\bar N^2}{2A_0^2}
  \end{displaymath}
  by taking $N\geq n_0$ and $n_0$ large enough.
  By plugging this bound into \eqref{eq:cmr:10} we get \eqref{eq:cmr:9} in this case also.
  
  \bigskip
  
  We can now conclude the proof of the proposition.
  We will bound from above
  \begin{displaymath}
    B_0^+(N)=
    \sup_{n\in\{\bar N+1,\dots,N\}}\left(\sum_{k=n}^N\gamma^N(k)\right)\log\left(\frac{1}{\sum_{k=n}^N\gamma^N(k)}\right)\left(\sum_{k=\bar N+1}^{n}\frac{1}{\gamma^N(k)\wedge \gamma^N(k+1)}\right),
  \end{displaymath}
  the bound of $B_0^-(N)$ being similar.
  If $n\in\{\bar N+1,\dots,\bar N+\lfloor\sqrt{\bar N}\rfloor\}$ by \eqref{eq:cmr:3} we obtain
  \begin{multline*}
    \sum_{k=\bar N+1}^{n}\frac{1}{\gamma^N(k)\wedge \gamma^N(k-1)}\leq
    \sum_{k=\bar N+1}^{n}\frac{1}{\gamma^N(k)}\leq
    A_0\sqrt{\bar N}\sum_{k=\bar N+1}^{n}e^{\frac{A_0(k-\bar N)^2}{\bar N}}\\
    \leq  A_0\sqrt{\bar N}(n-\bar N)e^{\frac{A_0(n-\bar N)^2}{\bar N}}
    \leq  A_0\bar Ne^{A_0},  
  \end{multline*}
  while because $\sum_{k=n}^N\gamma^N(k)\leq 1$ and $-x\log x\leq1$ for any $x\in[0,1]$,
  \begin{displaymath}
    \left(\sum_{k=n}^N\gamma^N(k)\right)\log\left(\frac{1}{\sum_{k=n}^N\gamma^N(k)}\right)\leq 1.
  \end{displaymath}
  This implies
    \begin{equation}
      \label{eq:cmr:11}
    \left(\sum_{k=n}^N\gamma^N(k)\right)\log\left(\frac{1}{\sum_{k=n}^N\gamma^N(k)}\right)\left(\sum_{k=\bar N+1}^{\bar N}\frac{1}{\gamma^N(k)\wedge \gamma^N(k+1)}\right)\leq
    A_0\bar Ne^{A_0},
  \end{equation}
  for any $n\in\{\bar N+1,\dots,\bar N+\lfloor\sqrt{\bar N}\rfloor\}$.

  Assume now $n\in\{\bar N+\lfloor\sqrt{\bar N}\rfloor+1,\dots,N\}$. 
  By \eqref{eq:cmr:4}, \eqref{eq:cmr:9} and \eqref{eq:cmr:7} we obtain
  \begin{multline*}
    \left(\sum_{k=n}^N\gamma^N(k)\right)\log\left(\frac{1}{\sum_{k=n}^N\gamma^N(k)}\right)\left(\sum_{k=\bar N+1}^{n}\frac{1}{\gamma^N(k)\wedge \gamma^N(k+1)}\right)\\
    \leq\left(\sum_{k=n}^N\gamma^N(k)\right)\log\left(\frac{1}{\sum_{k=n}^N\gamma^N(k)}\right)\left(\sum_{k=\bar N+1}^{n}\frac{1}{\gamma^N(k)}\right)\\
   \leq B_1B_2\left(\frac{\bar N}{n-\bar N}\right)^2\left[\log B_3+\frac{(n-\bar N)^2}{\bar N}\right]\leq
   B_1B_2A_0^2(\log B_3+A_0\bar N).
  \end{multline*}
  By the previous bound and \eqref{eq:cmr:11} we have that $B_0^+(N)\leq A_1 N$ for any $N\in\natural\setminus\{0\}$.
\end{proof}

\subsection{Proof of Proposition~\ref{ppaolo1}}

Recall that for any $\Lambda^\prime\subset\Lambda\subset\subset\integer^d$ such that $|\Lambda|=2|\Lambda^\prime|$ and any $N\in\natural\setminus\{0\}$ we defined $\gamma_\Lambda^N=\nu_\Lambda^N(\n\eta{\Lambda^\prime}=\cdot)$.
We have to show that any $f:\{0,\dots,N\}\to\real_+$
  \begin{equation}
    \label{pro:ols}
    \ent_{\gamma_\Lambda^N}(f)\leq
    A_0N\sum_{n=1}^N[\gamma_\Lambda^N(n)\wedge\gamma_\Lambda^N(n-1)]
[\sqrt{f(n-1)}-\sqrt{f(n)}]^2.
  \end{equation}

\noindent
The assumptions of Proposition~\ref{pro:cmr} are too strong to be fulfilled by $\gamma^N=\gamma_\Lambda^N$.
In particular, while assumption \eqref{eq:cmr:3} holds, as we will see \textit{a fortiori}, assumptions \eqref{eq:cmr:1} and \eqref{eq:cmr:2} may not be fulfilled in general (the case $|\Lambda|=2$ may be instructive to see this). 
  So following \cite{Ca:Ma:Ro}, for any $\epsilon\in(0,1/4)$, we consider a regularization $\tilde\gamma_\Lambda^{N,\epsilon}$ of $\gamma_\Lambda^N$.
  This regularization will be equivalent to $\gamma_\Lambda^N$ for any $\epsilon\in(0,1/4)$, so by the comparison criterion 
in Theorem 3.4.3 of \cite{toulouse} we can replace $\gamma_\Lambda^N$ with $\tilde\gamma_\Lambda^{N,\epsilon}$ in \eqref{pro:ols}.
  Finally we will prove that there exists $\epsilon_0\in(0,1/4)$ such that $\gamma^N:=\tilde\gamma_\Lambda^{N,\epsilon_0}$ fulfills the hypothesis of Proposition~\ref{pro:cmr}.

    Assume $N\in\natural\setminus\{0\}$, $n\in\{0,\dots,N\}$ and $\Lambda\subset\subset\integer^d$ and
  define $I_{N,\epsilon}:=[\epsilon N,(1-\epsilon)N]\cap\integer$, $\bar N:=\lceil N/2\rceil$,
  \begin{displaymath}
    H_\Lambda^N(n):=
    \log\left[\frac{p_{\Lambda^\prime}^{n/|\Lambda^\prime|}(n)p_{\Lambda^\prime}^{(N-n)/|\Lambda^\prime|}(N-n)}{p_{\Lambda^\prime}^{\bar N/|\Lambda^\prime|}(n)p_{\Lambda^\prime}^{\bar N/|\Lambda^\prime|}(N-n)}\right],
    \qquad\qquad
    Z_\Lambda^{N,\epsilon}:=
    \frac{\sum_{k\in I_{N,\epsilon}}e^{-H_\Lambda^N(k)}}{\sum_{k\in I_{N,\epsilon}}\gamma_\Lambda^N(k)}
  \end{displaymath}
  and
  \begin{displaymath}
    \tilde\gamma_\Lambda^{N,\epsilon}(n):=
    \begin{cases}
      \frac{1}{Z_\Lambda^{N,\epsilon}}e^{-H_\Lambda^N(n)} & \text{if $n\in I_{N,\epsilon}$,}\\
      \gamma_\Lambda^N(n) & \text{if $n\not\in I_{N,\epsilon}$.}
    \end{cases}
  \end{displaymath}
  It easy to check that $\tilde\gamma_\Lambda^{N,\epsilon}$ is a probability density supported on $\{0,\dots,N\}$.
  Now we show that $\tilde\gamma_\Lambda^{N,\epsilon}$ is equivalent to $\gamma_\Lambda^N$.

  \begin{lem}
    \label{lemma:eqm}
  For any fixed $\epsilon\in(0,1/4)$ there exists a positive constant $A_0$ such that
    \begin{equation}
      \label{eq:ols:1}
    \frac{1}{A_0}\leq
    \inf_{\substack{N\in\natural\setminus\{0\},\,\Lambda\subset\subset\integer^d\\ n\in\{0,\dots,N\}}}\frac{\gamma_\Lambda^N(n)}{\tilde\gamma_\Lambda^{N,\epsilon}(n)}\leq
    \sup_{\substack{N\in\natural\setminus\{0\},\,\Lambda\subset\subset\integer^d\\ n\in\{0,\dots,N\}}}\frac{\gamma_\Lambda^N(n)}{\tilde\gamma_\Lambda^{N,\epsilon}(n)}\leq
    A_0
  \end{equation}    
  \end{lem}

\begin{proof}
  Fix $\epsilon\in(0,1/4)$, since $\gamma_\Lambda^N(n)/\tilde\gamma_\Lambda^{N,\epsilon}(n)=1$ for any $n\not\in I_{N,\epsilon}$, we have to bound the ratio  $\gamma_\Lambda^N(n)/\tilde\gamma_\Lambda^{N,\epsilon}(n)$ for $n\in I_{N,\epsilon}$ only.
  So assume that $n\in I_{N,\epsilon}$ and define
  \begin{displaymath}
    \varphi_\Lambda^N(n):=
    \frac{p_{\Lambda^\prime}^{n/|\Lambda^\prime|}(n)p_{\Lambda^\prime}^{(N-n)/|\Lambda^\prime|}(N-n)}{p_{\Lambda^\prime}^{n/|\Lambda^\prime|}(n)},
  \end{displaymath}
  so that
  \begin{displaymath}
    \frac{\gamma_\Lambda^N(n)}{\tilde\gamma_\Lambda^{N,\epsilon}(n)}
    =Z_\Lambda^{N,\epsilon}\varphi_\Lambda^N(n).
  \end{displaymath}
  But
  \begin{displaymath}
    Z_\Lambda^{N,\epsilon}=
   \frac{\sum_{k\in I_{N,\epsilon}}e^{-H_\Lambda^N}(k)}{\sum_{k\in I_{N,\epsilon}}\gamma_\Lambda^N(k)}=
   \frac{\sum_{k\in I_{N,\epsilon}}\frac{\gamma_\Lambda^N(k)}{\varphi_\Lambda^N(k)}}{\sum_{k\in I_{N,\epsilon}}\gamma_\Lambda^N(k)},
  \end{displaymath}
  which implies
  \begin{displaymath}
   \frac{\gamma_\Lambda^N(n)}{\tilde\gamma_\Lambda^{N,\epsilon}(n)}=
   \frac{\sum_{k\in I_{N,\epsilon}}\gamma_\Lambda^N(k)\frac{\varphi_\Lambda^N(n)}{\varphi_\Lambda^N(k)}}{\sum_{k\in I_{N,\epsilon}}\gamma_\Lambda^N(k)}.
  \end{displaymath}
  Therefore
  \begin{displaymath}
    \min_{n,m\in I_{N,\epsilon}}\frac{\varphi_\Lambda^N(n)}{\varphi_\Lambda^N(m)}\leq
    \frac{\gamma_\Lambda^N(n)}{\tilde\gamma_\Lambda^{N,\epsilon}(n)}\leq
    \max_{n,m\in I_{N,\epsilon}}\frac{\varphi_\Lambda^N(n)}{\varphi_\Lambda^N(m)},
  \end{displaymath}
  for any $n\in I_{N,\epsilon}$, $N\in\natural\setminus\{0\}$ and $\Lambda\subset\subset\integer^d$ with $|\Lambda|$ even.
  Furthermore notice that  
  \begin{displaymath}
    \frac{\varphi_\Lambda^N(n)}{\varphi_\Lambda^N(m)}=
    \frac{p_{\Lambda^\prime}^{n/|\Lambda^\prime|}(n)p_{\Lambda^\prime}^{(N-n)/|\Lambda^\prime|}(N-n)}{p_{\Lambda^\prime}^{m/|\Lambda^\prime|}(n)p_{\Lambda^\prime}^{(N-m)/|\Lambda^\prime|}(N-m)}
  \end{displaymath}
  By Proposition~\ref{pro:cen} there exists $C_1>0$
  \begin{multline*}
        \frac{1}{C_1}\sqrt{\frac{\sigma^2(m/|\Lambda^\prime|)\sigma^2((N-m)/|\Lambda^\prime|)}{\sigma^2(n/|\Lambda^\prime|)\sigma^2((N-n)/|\Lambda^\prime|)}}
       \leq \frac{p_{\Lambda^\prime}^{n/|\Lambda^\prime|}(n)p_{\Lambda^\prime}^{(N-n)/|\Lambda^\prime|}(N-n)}{p_{\Lambda^\prime}^{m/|\Lambda^\prime|}(n)p_{\Lambda^\prime}^{(N-m)/|\Lambda^\prime|}(N-m)}\\
       \leq {C_1}\sqrt{\frac{\sigma^2(m/|\Lambda^\prime|)\sigma^2((N-m)/|\Lambda^\prime|)}{\sigma^2(n/|\Lambda^\prime|)\sigma^2((N-n)/|\Lambda^\prime|)}},
  \end{multline*}
  for any $n\in \{0,\dots, N\}$, $N\in\natural\setminus\{0\}$, $\Lambda^\prime\subset\subset\integer^d$.
  So \eqref{eq:ols:1} is proved if we can show that
  \begin{displaymath}
    0<
    \inf_{\substack{n,m\in I_{N,\epsilon} \\ N,v\in\natural\setminus\{0\}}}\frac{\sigma^2(m/v)\sigma^2((N-m)/v)}{\sigma^2(n/v)\sigma^2((N-n)/v)}\leq
    \sup_{\substack{n,m\in I_{N,\epsilon} \\ N,v\in\natural\setminus\{0\}}}\frac{\sigma^2(m/v)\sigma^2((N-m)/v)}{\sigma^2(n/v)\sigma^2((N-n)/v)}<
    +\infty.
  \end{displaymath}
  Indeed, by \eqref{eq:sigma}, there exists a constant $C_2>0$ such that
  \begin{displaymath}
    \frac{m(N-m)}{C_2n(N-n)}\leq
    \frac{\sigma^2(m/v)\sigma^2((N-m)/v)}{\sigma^2(n/v)\sigma^2((N-n)/v)}\leq
    \frac{C_2m(N-m)}{n(N-n)}
  \end{displaymath}
  for any $n,m\in\{0,\dots,N\}$, $N,v\in\natural\setminus\{0\}$ and a trivial calculation shows that
  \begin{displaymath}
    4\epsilon(1-\epsilon)\leq
    \frac{m(N-m)}{n(N-n)}\leq
    \frac{1}{4\epsilon(1-\epsilon)}
  \end{displaymath}
  for any $n,m\in I_{N,\epsilon}$.
\end{proof}

The rest of the paper is devoted to prove that conditions \eqref{eq:cmr:1}, \eqref{eq:cmr:2} and \eqref{eq:cmr:3} hold for $\tilde\gamma_\Lambda^{N,\epsilon_0}$ and some $\epsilon_0\in(0,1/4)$. 
We begin by showing the exponential decay of the tails of $\tilde\gamma_\Lambda^{N,\epsilon_0}$.
\begin{lem}
  \label{lemma:cod}
  There exists $\epsilon_0\in(0,1/4)$ such that
  \begin{align}
    \frac{\tilde\gamma_\Lambda^{N,\epsilon_0}(n+1)}{\tilde\gamma_\Lambda^{N,\epsilon_0}(n)} &\leq \frac{1}{2}\qquad\text{for any } n\in [(1-\epsilon_0)N,N-1]\cap\integer; \\
        \frac{\tilde\gamma_\Lambda^{N,\epsilon_0}(n-1)}{\tilde\gamma_\Lambda^{N,\epsilon_0}(n)} &\leq \frac{1}{2}\qquad\text{for any } n\in [1,\epsilon_0 N]\cap\integer;
  \end{align}
  for any $N\in\natural\setminus\{0\}$ and any $\Lambda\subset\subset\integer^d$ with $|\Lambda|\geq 2$.
\end{lem}

\begin{proof}
  Because $\gamma_\Lambda^{N,\epsilon_0}(n)=\gamma_\Lambda^N(n)$ for $n\not\in I_{N,\epsilon_0}$ the lemma is a trivial consequence of Proposition~\ref{pro:dec}. 
\end{proof}

\begin{lem}
  \label{lemma:uc}
  For any $\epsilon\in(0,1/4)$ there exists a positive constant $A_0$ such that for any $N\in\natural\setminus\{0\}$, $\Lambda\subset\subset\integer^d$
  \begin{align}
    \label{eq:uc:1}
    \frac{n-\bar N}{A_0\bar N}&\leq H_\Lambda^N(n+1)- H_\Lambda^N(n)\leq A_0\frac{n-\bar N}{\bar N}\qquad\text{for any $n\in I_{N,\epsilon}$ such that $n> \bar N$,}\\
    \label{eq:uc:2}
    \frac{n-\bar N}{A_0\bar N}&\leq H_\Lambda^N(n-1)- H_\Lambda^N(n)\leq A_0\frac{n-\bar N}{\bar N}\qquad\text{for any $n\in I_{N,\epsilon}$ such that $n< \bar N$.}
  \end{align}
\end{lem}

\begin{proof}
  Define $v:=|\Lambda^\prime|=|\Lambda|/2$.
  An elementary computation gives
  \begin{multline*}
    H_\Lambda^N(n):=
    \log\left[\frac{p_{\Lambda^\prime}^{n/|\Lambda^\prime|}(n)p_{\Lambda^\prime}^{(N-n)/|\Lambda^\prime|}(N-n)}{p_{\Lambda^\prime}^{\bar N/|\Lambda^\prime|}(n)p_{\Lambda^\prime}^{\bar N/|\Lambda^\prime|}(N-n)}\right]=
    n\log\alpha(n/v)+(N-n)\log\alpha((N-n)/v)\\
    +2v\log Z(\bar N/v)-v\log Z(n/v) -v \log Z((N-n)/v)-N\log \alpha(\bar N/v).
  \end{multline*}
  This formula shows that it is possible to extend the function $H_\Lambda^N:\{0,\dots,N\}\to \real$ to a real function, $H_\Lambda^N:[0,N]\to \real$ by defining 
  \begin{multline*}
    H_\Lambda^N(x):=
    x\log\alpha(x/v)+(N-x)\log\alpha((N-x)/v)\\
    +2v\log Z(\bar N/v)-v\log Z(x/v) -v \log Z((N-x)/v)-N\log \alpha(\bar N/v)
  \end{multline*}
  for any $x\in[0,N]$.
  A direct calculation gives
  \begin{multline*}
    \frac{dH_\Lambda^N}{dx}=
    \log\alpha(x/v)+(x/v)[\alpha^\prime(x/v)/\alpha(x/v)]-\log\alpha((N-x)/v)\\
    -[(N-x)/v][\alpha^\prime((N-x)/v)/\alpha((N-x)/v)]-Z^\prime(x/v)/Z(x/v) +Z^\prime((N-x)/v)/Z((N-x)/v).
  \end{multline*}
  Since
  \begin{displaymath}
    Z^\prime(\rho)=
    \frac{d}{d\rho}\sum_{k=0}^{+\infty}\frac{\alpha(\rho)^{k}}{c(k)!}=
    \alpha^\prime(\rho)\sum_{k=1}^{+\infty}\frac{k\alpha(\rho)^{k-1}}{c(k)!}=
    \frac{\alpha^\prime(\rho)}{\alpha(\rho)}\sum_{k=0}^{+\infty}\frac{k\alpha(\rho)^{k}}{c(k)!}=
    \frac{\alpha^\prime(\rho)\rho Z(\rho)}{\alpha(\rho)},
  \end{displaymath}
  which, by equation (1.3) of \cite{La:Se:Va} (namely $\alpha^\prime(\rho)=\alpha(\rho)\sigma^2(\rho)$), implies
  \begin{displaymath}
    \frac{Z^\prime(\rho)}{Z(\rho)}=
    \frac{\alpha^\prime(\rho)\rho}{\alpha(\rho)}=
    \frac{\rho}{\sigma^2(\rho)},
  \end{displaymath}
  we obtain
  \begin{displaymath}
    \frac{dH_\Lambda^N}{dx}=\log\alpha(x/v)-\log\alpha((N-x)/v).
  \end{displaymath}
  Differentiating this identity and using again equation (1.3) of \cite{La:Se:Va}, we obtain
  \begin{displaymath}
    \frac{d^2H_\Lambda^N}{dx^2}=
    \frac{1}{\sigma^2(x/v)v}+\frac{1}{\sigma^2((N-x)/v)v},
  \end{displaymath}
  which, again with \eqref{eq:sigma}, implies that there exist a positive constant $B_1$ such that for any $x\in(0,N)$
  \begin{displaymath}
    \frac{1}{B_1}\left(\frac{1}{x}+\frac{1}{N-x}\right)\leq
      \frac{d^2H_\Lambda^N}{dx^2}\leq
      B_1\left(\frac{1}{x}+\frac{1}{N-x}\right).
  \end{displaymath}
  Assume now that $x\in [\bar N, (1-\epsilon)N]$.
  Then the previous inequality yields
  \begin{displaymath}
        \frac{1}{B_1}\left(\frac{1}{N(1-\epsilon)}+\frac{1}{\bar N}\right)\leq
      \frac{d^2H_\Lambda^N}{dx^2}\leq
      B_1\left(\frac{1}{\bar N}+\frac{1}{\epsilon N}\right)
  \end{displaymath}
  it follows that, if $\epsilon>0$ is small enough, there exists a positive constant $B_2(\epsilon)$ such that
  \begin{displaymath}
            \frac{1}{B_2\bar N}\leq
      \frac{d^2H_\Lambda^N}{dx^2}\leq
      \frac{B_2}{\bar N}
  \end{displaymath}
  for any $N\in\natural\setminus\{0\}$, $\Lambda\subset\subset\integer^d$.
  Integrating this inequality with respect to $x$ from $\bar N$ to $y\in[\bar N, (1-\epsilon)N]$, we obtain
  \begin{displaymath}
                \frac{y-\bar N}{B_2\bar N}\leq
      \frac{dH_\Lambda^N}{dy}\leq
      B_2\frac{y-\bar N}{\bar N}
  \end{displaymath}
  Now take $n\in I_{N,\epsilon}$ such that $n>\bar N$ and integrate the previous inequality with respect to $y$ from $n$ to $n+1$ to obtain
  \begin{displaymath}
      \frac{n-\bar N}{B_2\bar N}\leq
      H_\Lambda^N(n+1)-H_\Lambda^N(n)\leq
      B_2\frac{n+1-\bar N}{\bar N}.
  \end{displaymath}
  Equation \eqref{eq:uc:1} now follows because  $(n+1-\bar N)/\bar N\leq 2(n-\bar N)/\bar N$ for any $n>\bar N$.
  Equation \eqref{eq:uc:2} can be obtained in a similar way.
\end{proof}

By Lemma~\ref{lemma:cod} and Lemma~\ref{lemma:uc} follows easily

\begin{cor}
  \label{cor:ed}
  Define $\bar N:=\lceil N/2 \rceil$.
  There exists $\epsilon_0\in(0,1/4)$ and a positive constant $A_0$ such that for any $N\in\natural\setminus\{0\}$ and $\Lambda\subset\subset\integer^d$
  \begin{align}
    \label{eq:ed:1}
    \frac{\tilde\gamma_\Lambda^{N,\epsilon_0}(n+1)}{\tilde\gamma_\Lambda^{N,\epsilon_0}(n)} &\leq e^{-\frac{n-\bar N}{A_0 \bar N}}\qquad\qquad\qquad\ \;\text{for all  $n\in\{\bar N+1,\dots,N\}$}\\
    \label{eq:ed:2}
    \frac{\tilde\gamma_\Lambda^{N,\epsilon_0}(n-1)}{\tilde\gamma_\Lambda^{N,\epsilon_0}(n)} &\leq e^{-\frac{\bar N-n}{A_0 \bar N}}\qquad\qquad\qquad\ \;\text{for all  $n\in\{0,\dots,\bar N-1\}$.}
  \end{align}
\end{cor}

In order to use Proposition~\ref{pro:cmr} it remains to prove the Gaussian lower bound \eqref{eq:cmr:3} for $\tilde\gamma_\Lambda^{N,\epsilon_0}$.
In fact the lower bound obtained in Proposition~\ref{pro:gaul} is not useful, because it is not uniform in $|\Lambda|$. 
The proof is almost the same as the proof of Lemma~\ref{lemma:gau2}

\begin{lem}
  \label{lemma:gau}
  Define $\bar N:=\lceil N/2 \rceil$ for $N\in\natural\setminus\{0\}$.
  There exist $\epsilon_0\in(0,1/4)$ and a positive constant $A_0$ such that
  \begin{equation}
    \label{eq:gau}
    \frac{1}{A_0\sqrt{\bar N}}e^{-\frac{A_0(n-\bar N)^2}{\bar N}}\leq
    \tilde\gamma_\Lambda^{N,\epsilon_0}(n)\leq
    \frac{A_0}{\sqrt{\bar N}}e^{-\frac{(n-\bar N)^2}{A_0\bar N}},
  \end{equation}
  uniformly in $N\in\natural\setminus\{0\}$ and $n\in\{0,1,\dots,N\}$.
\end{lem}

\begin{proof}
  We split the proof in several steps for clarity purpose.

  \paragraph{Step 1.}

  There exist $A_1>0$ such that for any $N\in\natural\setminus\{0\}$,
  \begin{equation}
    \label{eq:gau:stp3}
    \frac{\tilde\gamma_\Lambda^{N,\epsilon_0}(n)}{\tilde\gamma_\Lambda^{N,\epsilon_0}(\bar N)}\leq
    A_1e^{-\frac{(n-\bar N)^2}{A_1\bar N}},
  \end{equation}
  for any $n\in\{0,\dots,N\}$.

  \subparagraph{Proof of Step 1.}

  Assume that $n\in\{\bar N+1,\dots,N\}$, by \eqref{eq:ed:1} we get $B_1>0$ such that 
  \begin{equation}
    \label{eq:gau:stp3:1}
    \log \tilde\gamma_\Lambda^{N,\epsilon_0}(n)-\log \tilde\gamma_\Lambda^{N,\epsilon_0}(\bar N)=
    \sum_{k=\bar N}^{n-1} \left[\log \tilde\gamma_\Lambda^{N,\epsilon_0}(k+1)-\log\tilde\gamma_\Lambda^{N,\epsilon_0}(k)\right]\leq
    -\sum_{k=\bar N}^{n-1}\frac{k-\bar N}{B_1\bar N}.
  \end{equation}
  Using the fact that $n-\bar N\leq N/2$, and some elementary computation we obtain
  \begin{displaymath}
    \sum_{k=\bar N}^{n-1}\frac{k-\bar N}{B_1\bar N}\geq
    \frac{1}{B_1\bar N}\sum_{k=1}^{n-\bar N-1}k\geq
    \frac{(n-\bar N)(n-\bar N-1)}{B_1 N}=
    \frac{(n-\bar N)^2}{B_1 N}-\frac{n-\bar N}{B_1 N}\geq
    \frac{(n-\bar N)^2}{B_1 N}-\frac{1}{2B_1},
  \end{displaymath}
  which plugged into \eqref{eq:gau:stp3:1} implies \eqref{eq:gau:stp3} for $n\in\{\bar N+1,\dots,N-1\}$.
  The case $n\in\{0,\dots,\bar N\}$ is similar while that case $n=\bar N$ is trivial.

  \paragraph{Step 2.}

    There exists $A_2>0$ such that if $N\in\natural\setminus\{0\}$ then
  \begin{equation}
    \label{eq:gau:stp4}
    \frac{\tilde\gamma_\Lambda^{N,\epsilon_0}(n)}{\tilde\gamma_\Lambda^{N,\epsilon_0}(\bar N)}\geq
    \frac{1}{A_2}e^{-A_2\frac{(n-\bar N)^2}{\bar N}},
  \end{equation}
  for any $n\in I_{N,\epsilon_0}$.
    
  \subparagraph{Proof of Step 2.}

  Fix $n\in I_{N,\epsilon_0}$ and assume that $n >\bar N$. 
  By \eqref{eq:uc:1} we get $B_1>0$ such that 
  \begin{equation}
    \label{eq:gau:stp4:1}
    \log \tilde\gamma_\Lambda^{N,\epsilon_0}(n)-\log \tilde\gamma_\Lambda^{N,\epsilon_0}(\bar N)=
    \sum_{k=\bar N}^{n-1} \left[\log \tilde\gamma_\Lambda^{N,\epsilon_0}(k+1)-\log\tilde\gamma_\Lambda^{N,\epsilon_0}(k)\right]\geq
    -B_1\sum_{k=\bar N}^{n-1}\frac{k-\bar N}{\bar N}
  \end{equation}
  but, using the fact that $0\leq n-\bar N\leq N/2$, and some elementary computation we obtain
  \begin{displaymath}
    \sum_{k=\bar N}^{n-1}\frac{k-\bar N}{\bar N}\leq
    \frac{1}{\bar N}\sum_{k=1}^{n-\bar N-1}k\leq
    \frac{(n-\bar N)(n-\bar N-1)}{ N}=
    \frac{(n-\bar N)^2}{N}-\frac{n-\bar N}{N}\leq
    \frac{(n-\bar N)^2}{N},
  \end{displaymath}
  which plugged into \eqref{eq:gau:stp4:1} implies \eqref{eq:gau:stp4} for $n\geq\bar N$.
  The case $n< \bar N$ is similar while the case $n=\bar N$ is trivial.

  \paragraph{Step 3.}

    There exists $A_3>0$ such that if $N\in\natural\setminus\{0\}$ then
  \begin{equation}
    \label{eq:gau:stp5}
    \frac{\tilde\gamma_\Lambda^{N,\epsilon_0}(n)}{\tilde\gamma_\Lambda^{N,\epsilon_0}(\bar N)}\geq
    \frac{1}{A_3}e^{-A_3\frac{(n-\bar N)^2}{\bar N}},
  \end{equation}
  for any $n\in\{0,\dots,N\}\setminus I_{N,\epsilon_0}$.
    
  \subparagraph{Proof of Step 3.}

  Observe that $n\in\{0,\dots,N\}\setminus I_{N,\epsilon_0}$ and $\epsilon_0\in(0,1/4)$ implies $|n-\bar N|>N/8$.
  Thus to prove \eqref{eq:gau:stp5} we have only to show that the ratio $\tilde\gamma_\Lambda^{N,\epsilon_0}(n)/\tilde\gamma_\Lambda^{N,\epsilon_0}(\bar N)$ is bounded from below by a negative exponential of $N$.
  Assume that $n\in\{\lceil(1-\epsilon_0)N\rceil+1,\dots,N\}$;
  then
  \begin{equation}
    \label{eq:gau:stp5:0}
    \frac{\tilde\gamma_\Lambda^{N,\epsilon_0}(n)}{\tilde\gamma_\Lambda^{N,\epsilon_0}(\bar N)}=
    \frac{\tilde\gamma_\Lambda^{N,\epsilon_0}(\lfloor(1-\epsilon_0)N\rfloor)}{\tilde\gamma_\Lambda^{N,\epsilon_0}(\bar N)}\prod_{k=\lfloor(1-\epsilon_0)N\rfloor}^{n-1}\frac{\tilde\gamma_\Lambda^{N,\epsilon_0}(k+1)}{\tilde\gamma_\Lambda^{N,\epsilon_0}(k)}.
  \end{equation}
  By \eqref{eq:gau:stp4} there exists $B_1>0$ such that
  \begin{displaymath}
    \frac{\tilde\gamma_\Lambda^{N,\epsilon_0}(\lfloor(1-\epsilon_0)N\rfloor)}{\tilde\gamma_\Lambda^{N,\epsilon_0}(\bar N)}\geq
    \frac{1}{B_1}e^{-\frac{B_1(\lfloor(1-\epsilon_0)N\rfloor-\bar N)^2}{\bar N}}.
  \end{displaymath}
  But
  \begin{displaymath}
    \left|\lfloor(1-\epsilon_0)N\rfloor-\bar N\right|=
    \lfloor(1-\epsilon_0)N\rfloor-\bar N\leq
    (1-\epsilon_0)N-\bar N\leq
    \bar N,
  \end{displaymath}
  so
  \begin{displaymath}
    \frac{\tilde\gamma_\Lambda^{N,\epsilon_0}(\lfloor(1-\epsilon_0)N\rfloor)}{\tilde\gamma_\Lambda^{N,\epsilon_0}(\bar N)}\geq
    \frac{1}{B_1}e^{-B_1\bar N}
  \end{displaymath}
  and because $\bar N\leq N$,
  \begin{equation}
    \label{eq:gau:stp5:1}
    \frac{\tilde\gamma_\Lambda^{N,\epsilon_0}(\lfloor(1-\epsilon_0)N\rfloor)}{\tilde\gamma_\Lambda^{N,\epsilon_0}(\bar N)}\geq
    \frac{1}{B_1}e^{-B_1N}.
  \end{equation}
  In order to bound the product factor in the right hand side of \eqref{eq:gau:stp5:0}, notice that by Proposition~\ref{pro:dec} there exists $B_2>0$ such that
  \begin{displaymath}
    \frac{\tilde\gamma_\Lambda^{N,\epsilon_0}(k+1)}{\tilde\gamma_\Lambda^{N,\epsilon_0}(k)}\geq
    \frac{N-k}{B_2(k+1)},
  \end{displaymath}
  for any $k\in\{1,\dots,N-1\}$.
  Thus
  \begin{multline*}
    \sum_{k=\lfloor(1-\epsilon_0)N\rfloor}^{n-1}\left[\log\tilde\gamma_\Lambda^{N,\epsilon_0}(k)-\log\tilde\gamma_\Lambda^{N,\epsilon_0}(k+1)\right]\\
    \leq (n-\lfloor(1-\epsilon_0)N\rfloor)\log B_2+\sum_{k=\lfloor(1-\epsilon_0)N\rfloor}^{n-1}\left[\log(k+1)-\log(N-k)\right]\\
    \leq N\log B_2+\sum_{k=\lfloor(1-\epsilon_0)N\rfloor}^{n-1}\left[k+1-(N-k)\right]\max_{x\in(N-k,k+1)}\frac{d}{dx}\log x\\
    \leq N\log B_2+\sum_{k=\lfloor(1-\epsilon_0)N\rfloor}^{n-1}\frac{2k+1-N}{N-k}\leq
     N\log B_2+(n-\lfloor(1-\epsilon_0)N\rfloor)\frac{2(n-1)+1-N}{N-(1-\epsilon_0)N}\\
    \leq N\log B_2+N\frac{N-1}{N-3N/4}\leq
    N(\log B_2+\epsilon_0^{-1})
  \end{multline*}
  which implies
  \begin{displaymath}
    \prod_{k=\lfloor(1-\epsilon_0)N\rfloor}^{n-1}\frac{\tilde\gamma_\Lambda^{N,\epsilon_0}(k+1)}{\tilde\gamma_\Lambda^{N,\epsilon_0}(k)}\geq
    e^{-(\log B_2+\epsilon_0^{-1})N}.
  \end{displaymath}
  By plugging this bound and \eqref{eq:gau:stp5:1} into \eqref{eq:gau:stp5:0}, we get \eqref{eq:gau:stp5}.

  \paragraph{Step 4.}

  There exists $A_4>0$ such that if $N\in\natural\setminus\{0\}$ and $\Lambda\subset\subset\integer^d$, then,
  \begin{equation}
    \label{eq:gau:stp6}
    \frac{1}{A_4\sqrt{\bar N}}\leq
    \tilde\gamma_\Lambda^{N,\epsilon_0}(\bar N)\leq
    \frac{A_4}{\sqrt{\bar N}}.
  \end{equation}

  \subparagraph{Proof of Step 4.}

  By \eqref{eq:gau:stp3}, \eqref{eq:gau:stp4} and \eqref{eq:gau:stp5} we obtain $B_1>0$ such that for any $N\in\natural\setminus\{0\}$
  \begin{displaymath}
    \frac{1}{B_1}e^{-\frac{B_1(n-\bar N)^2}{\bar N}}\leq
    \frac{\tilde\gamma_\Lambda^{N,\epsilon_0}(n)}{\tilde\gamma_\Lambda^{N,\epsilon_0}(\bar N)}\leq
    B_1e^{-\frac{(n-\bar N)^2}{B_1\bar N}}
  \end{displaymath}
  for any $n\in\{0,\dots,N\}$.
  By summing for $n\in\{0,\dots,N\}$ we have
   \begin{displaymath}
     \sum_{n=0}^N\frac{1}{B_1}e^{-\frac{B_1(n-\bar N)^2}{\bar N}}\leq
    \frac{1}{\tilde\gamma_\Lambda^{N,\epsilon_0}(\bar N)}\leq
     B_1\sum_{n=0}^Ne^{-\frac{(n-\bar N)^2}{B_1\bar N}}.
   \end{displaymath}
   Thus we are done if we can show that, for any $B_2>0$, there exists $B_3>0$ such that
   \begin{equation}
     \label{eq:gau:stp6:1}
     \frac{\sqrt{\bar N}}{B_3}\leq
     \sum_{n=0}^Ne^{-\frac{(n-\bar N)^2}{B_2\bar N}}\leq
     B_3\sqrt{\bar N}
   \end{equation}
   for any $N\in\natural\setminus\{0\}$. 
   An elementary computation gives
   \begin{multline*}
     \sum_{n=0}^{N}e^{-B_2\frac{(n-\bar N)^2}{\bar N}}\leq
     \sum_{n=0}^{\bar N}e^{-B_2\frac{(n-\bar N)^2}{\bar N}}+\sum_{n=\bar N}^{N}e^{-B_2\frac{(n-\bar N)^2}{\bar N}}\\
     \leq\int_{-\infty}^{\bar N}e^{-B_2\frac{(x-\bar N)^2}{\bar N}}\; dx+\int_{\bar N}^{+\infty}e^{-B_2\frac{(x-\bar N)^2}{\bar N}}\; dx +1=
     1+\sqrt{\frac{\pi\bar N}{B_2}}.
   \end{multline*}
   This proves the upper bound in \eqref{eq:gau:stp6:1}, the proof of the lower one is similar.

   \paragraph{Conclusion}

   By  \eqref{eq:gau:stp3}, \eqref{eq:gau:stp4}, \eqref{eq:gau:stp5} and \eqref{eq:gau:stp6} we get \eqref{eq:gau}.
 \end{proof}

We can now prove Proposition~\ref{ppaolo1}

\medskip

\begin{proofof}{Proposition~\ref{ppaolo1}}
  By Lemma~\ref{lemma:eqm} $\gamma_\Lambda^N$ and $\tilde\gamma_\Lambda^{N,\epsilon}$ are equivalent measures for any $N\in\natural\setminus\{0\}$, $\Lambda\subset\subset\integer^d$ and $\epsilon\in(0,1/4)$ so by the comparison criterion in Theorem 3.4.3 of \cite{toulouse}, in the proof we can replace $\gamma_\Lambda^N$ with $\tilde\gamma_\Lambda^{N,\epsilon}$.
  By Corollary~\ref{cor:ed} there exists $\epsilon_0\in(0,1/4)$ such that $\tilde\gamma_\Lambda^{N,\epsilon}$ satisfy  conditions \eqref{eq:cmr:1} and \eqref{eq:cmr:2} of Proposition~\ref{pro:cmr} uniformly in $\Lambda\subset\subset\integer^d$.
  Furthermore by Lemma~\ref{lemma:gau} condition \eqref{eq:cmr:3} too holds for $\tilde\gamma_\Lambda^{N,\epsilon}$ and we can apply Proposition~\ref{pro:cmr} to $\tilde\gamma_\Lambda^{N,\epsilon}$.
 \end{proofof}

Finally, we prove Proposition \ref{ppaolo6}

\begin{proofof}{Proposition~\ref{ppaolo6}}
  We may assume  $d=1$ and $\Lambda=\{0,1\}$.
  In this case $\nu_\Lambda^N(\eta_0=\xi_0,\eta_1=\xi_1)=\ind(\xi_1=N-\xi_0)\gamma_\Lambda^N(\xi_0)$.
  Thus $\ent_{\nu_\Lambda^N}(f)=\ent_{\gamma_\Lambda^N}\left(\tilde f\right)$ for any $f:\Omega_\Lambda\to\real_+$, where $\tilde f(n):=f(n,N-n)$, for any $n\in\natural$.
  For the same reason
  \begin{multline*}
    \mathcal{E}_{\nu_\Lambda^N}(\sqrt{f},\sqrt{f})\\
    =\frac{1}{2}\sum_{n=0}^N\gamma_\Lambda^N(n)\bigg\{c(n)\left[\sqrt{f}(n-1,N-n+1)- \sqrt{f}(n,N-n)\right]^2
      \\+c(N-n)\left[\sqrt{f}(n+1,N-n-1)-\sqrt{f}(n,N-n)\right]^2\bigg\}\\
    =\frac{1}{2}\sum_{n=0}^N\gamma_\Lambda^N(n)\left\{c(n)\left[\sqrt{\tilde f}(n-1)-\sqrt{\tilde f}(n)\right]^2+c(N-n)\left[\sqrt{\tilde f}(n+1)-\sqrt{\tilde f}(n)\right]^2\right\}\\
    =\sum_{n=1}^N\gamma_\Lambda^N(n)c(n)\left[\sqrt{\tilde f}(n-1)-\sqrt{\tilde f}(n)\right]^2
    :=\mathcal{\tilde D}\left(\sqrt{\tilde f},\sqrt{\tilde f}\right).
  \end{multline*}
  Recall now the definition of the Dirichlet form $\mathcal{D}(\varphi,\varphi)$ defined just before Proposition~\ref{ppaolo1}, namely
  \begin{displaymath}
    \mathcal{D}\left(\varphi,\varphi\right)
    =\sum_{n=1}^N[\gamma_\Lambda^N(n)\wedge\gamma_\Lambda^N(n-1)]\left[\varphi(n-1)-\varphi(n)\right]^2.
  \end{displaymath}
  We claim that there exists a constant $B_0>0$ such that for any $\varphi:\natural\to\real$ and any $N\in\natural\setminus\{0\}$ we have
  \begin{equation}
    \label{eq:dom}
    N\mathcal{D}(\varphi,\varphi)
    \leq B_0\mathcal{\tilde D}(\varphi,\varphi).
  \end{equation}
  In this case we will have by Proposition~\ref{ppaolo1}
  \begin{displaymath}
    \ent_{\nu_\Lambda^N}(f)
    =\ent_{\gamma_\Lambda^N}\left(\tilde f\right)
    \leq B_1 N\mathcal{D}\left(\sqrt{\tilde f},\sqrt{\tilde f}\right)
    \leq B_2 \mathcal{\tilde D}\left(\sqrt{\tilde f},\sqrt{\tilde f}\right)
    =B_2 \mathcal{E}_{\nu_\Lambda^N}\left( \sqrt{f}, \sqrt{f}\right),
  \end{displaymath}
  which is \eqref{not4} in the present case.
  Thus we have to verify \eqref{eq:dom}.
  Observe that
  \begin{multline*}
    \mathcal{\tilde D}(\varphi,\varphi)
    =\sum_{n=1}^N\gamma_\Lambda^N(n)c(n)\left[\varphi(n-1)-\varphi(n)\right]^2\\
    =\sum_{n=1}^N\frac{\gamma_\Lambda^N(n)c(n)}{\gamma_\Lambda^N(n)\wedge\gamma_\Lambda^N(n-1)}[\gamma_\Lambda^N(n)\wedge\gamma_\Lambda^N(n-1)]\left[\varphi(n-1)-\varphi(n)\right]^2\\
    =\sum_{n=1}^Nc(n)\left[1\vee\frac{\gamma_\Lambda^N(n)}{\gamma_\Lambda^N(n-1)}\right][\gamma_\Lambda^N(n)\wedge\gamma_\Lambda^N(n-1)]\left[\varphi(n-1)-\varphi(n)\right]^2.
  \end{multline*}
  Since, by \eqref{eq:gamma}, $\gamma_\Lambda^N(n)/\gamma_\Lambda^N(n-1)=c(N-n+1)/c(n)$  we get
  \begin{displaymath}
     \mathcal{\tilde D}(\varphi,\varphi)
    =\sum_{n=1}^N\left[c(n)\vee c(N-n+1)\right][\gamma_\Lambda^N(n)\wedge\gamma_\Lambda^N(n-1)]\left[\varphi(n-1)-\varphi(n)\right]^2.
  \end{displaymath}
  Then \eqref{eq:dom} follows by observing that by (\ref{lineargrowth}) there exists a constant $B_3>0$ such that $c(n)\vee c(N-n+1)\geq B_3^{-1}[n\vee (N-n+1)]\geq B_3^{-1}N/2$, for any $n\in\{1,\dots,N\}$.
\end{proofof}

\section*{Acknowledgments}
The authors wish to thank Lorenzo Bertini and Claudio Landim for very useful and interesting conversations on this work and Anna Maria Paganoni for her help in the first stage of the work.

\end{document}